\documentclass[reqno,letter]{amsart}
\usepackage{amsmath,amssymb,epsfig}
\usepackage{float}
\usepackage[pdftex,pdftitle={Adaptive Schemes for Reaction-Diffusion Systems},colorlinks=true,breaklinks=true,linkcolor=blue]{hyperref}

\allowdisplaybreaks

\newcommand{\eps}{\varepsilon}
\renewcommand{\epsilon}{\varepsilon}

\newcommand{\RR}{\mathbb{R}}

\newcommand{\ww}{\mathbf{w}}
\newcommand{\n}{\mathbf{n}}

\newcommand{\x}{\mathbf{x}}

\newcounter{dctr}[section]
\numberwithin{equation}{section}

\newtheorem{Alg}{Algorithm}[section]

\title[Adaptive Schemes for Reaction-Diffusion Systems]{Adaptive multiresolution schemes
with local time stepping for two-dimensional degenerate  reaction-diffusion systems}
\author[Bendahmane]{Mostafa Bendahmane$^{\mathrm{a}}$}
\author[B\"urger]{Raimund B\"urger$^{\mathrm{a}}$}
\author[Ruiz]{Ricardo Ruiz Baier$^{\mathrm{a}}$}
\author[Schneider]{Kai Schneider$^{\mathrm{b}}$}

\thanks{$^{\mathrm{a}}$Departamento de Ingenier\'{\i}a Matem\'{a}tica,
 Facultad de Ciencias F\'{\i}sicas y Matem\'{a}ticas,
 Universidad de Concepci\'{o}n, Casilla 160-C,  Concepci\'{o}n,
 Chile.
E-mail: {\tt mostafab@ing-mat.udec.cl}, {\tt rburger@ing-mat.udec.cl},
{\tt rruiz@ing-mat.udec.cl}}

\thanks{$^{\mathrm{b}}$Centre de Math\'ematiques et d'Informatique,
Universit\'e de Provence, 39 rue Joliot-Curie, 13453 Marseille cedex 13,
France.
E-mail: {\tt kschneid@cmi.univ-mrs.fr}}

\date{\today}

\begin{document}

\begin{abstract}
We present a fully adaptive multiresolution scheme  for spatially
two-dimensional, possibly degenerate reaction-diffusion systems,
focusing on combustion models and models of pattern formation and
chemotaxis in mathematical biology. Solutions of these equations
in these  applications exhibit steep gradients, and in the degenerate
case, sharp fronts and discontinuities. This calls for a concentration of
 computational effort in zones of strong variation.

The  multiresolution scheme is based on finite volume discretizations
 with explicit time stepping.  The multiresolution representation
 of the solution is stored in a graded tree (``quadtree''), whose
 leaves  are the non-uniform finite volumes on the borders of
 which the numerical divergence is evaluated.  By a  thresholding
procedure, namely the elimination of leaves that are  smaller
 than a threshold value, substantial data compression
and CPU time reduction is attained.  The threshold value is chosen
optimally, in the sense that the total error of the adaptive scheme
is of the same slope as that of the reference finite volume scheme.

Since chemical reactions involve a large range of temporal scales,
  but are spatially well localized (especially in the combustion
  model),  a locally varying adaptive time stepping strategy is applied.
  For scalar equations, this strategy  has the advantage
that consistence with a CFL condition is always enforced. Numerical experiments with five
 different scenarios, in part with local time stepping, illustrate the
effectiveness of the adaptive multiresolution method.
 It turns out that local time stepping accelerates the
adaptive multiresolution method by a factor of two, while the error
remains controlled.
\end{abstract}


\maketitle
\section{Introduction}
Multiresolution techniques were first introduced by Harten
\cite{harten:1995} to improve the performance of schemes for
 one-dimensional conservation laws. Later on, these
original ideas were  extended to several kinds of related problems
\cite{bh97,bks,cd01}, leading finally to the concept of fully
adaptive multiresolution schemes \cite{CKMP2002,dahmen,Muller,RSTB03}.
Overviews on multiresolution methods for conservation laws are
 given by Chiavassa, Donat, and M\"{u}ller \cite{chiavassaplewa}
and M\"{u}ller \cite{Muller}.
 The basic aim  of this approach is to accelerate a given finite
 volume scheme on a uniform grid at the cost of an at most
 controllable loss of accuracy, that is, the accelerated scheme should
  be of the same order than the original one. The
 principle of the multiresolution analysis is to represent a set of
 data given on a fine grid as values on a coarser grid plus a series
 of differences, called {\em details},
 at different levels of nested dyadic grids.  These
 differences contain    information on the local regularity
 of the solution. An appealing  feature of this data
 representation is that these details  are small in regions where
 the solution is smooth. By  thresholding  small
  details  (cells whose coefficients are smaller   than a
 prescribed tolerance are removed),  a locally refined adaptive grid is defined. This
 threshold is chosen such  that the discretization
 error of the reference scheme is balanced with the accumulated thresholding
 error  introduced in each time step.  Significant speed-up of the computation
 and data compression  is achieved for long-time evolution
problems, large systems, multidimensional domains, and
solutions with  sharp fronts.

The present paper serves two purposes.  On one hand,  the
adaptive multiresolution scheme for parabolic PDEs \cite{RSTB03} and
strongly degenerate parabolic PDEs in one space dimension
\cite{brss2,brss} is extended
 to two-dimensional systems of (possibly degenerate)
  parabolic PDEs. These  equations produce solutions that
 vary smoothly wherever the solution causes the PDE
to be parabolic, but produce sharp fronts, or even discontinuities,
  close to    solution values  at which the equation degenerates,
 so adaptive multiresolution methods are a proper device
to efficiently capture these fronts.  Similar features (that is,
 solutions with  steep gradients)   also  appear in a combustion model
 of reaction-diffusion type.
 The analysis made in \cite{ctm_rs04} is extended here
 to the study of two interacting flame balls. In this
case, an adaptive strategy is equally very useful, especially when the flame
front is well localized in space, since fine grids are only needed in
 small subregions of the computational domain.
 We will utilize then an adaptive multiresolution
scheme applied to a reference finite volume discretization with
explicit time integration.

On the other hand, chemical
 reactions are known to involve a large range of temporal scales,
 especially in long-time evolutions. Then  an adaptive time
stepping strategy is recommendable. Earlier efforts in this
direction, which include \cite{brss,chiu,DRS06,fehlberg} and the
references therein, were based on using the same time
 step to advance the solution on all parts of the computational
 domain, and controlling the time step  through  an
embedded pair of Runge-Kutta schemes (known
as Runge-Kutta-Fehlberg (RKF) schemes). In these procedures,
 one compares   the numerical solution after each  time
step with  an (approximate) reference solution, and adjusts the time step
 if the discrepancy is unacceptable. In contrast to  this approach,  we here adapt
the locally varying time stepping strategy recently introduced for
multiresolution schemes for conservation laws and multidimensional
systems by Lamby, M\"{u}ller, and Stiriba \cite{LMS}
 and M\"uller and Stiriba \cite{MS}. This  strategy
  is not precisely (time-)adaptive for scalar equations, since  the time step
for each level remains the same for all times. However, in the case of nonlinear systems,
coupling of components entering the CFL condition makes it necessary to compute the time
step after each iteration, according the evolving CFL condition, and therefore we have a
scheme adaptive in time. Our
 results in terms of CPU time savings are encouraging and the strategy  is
consistent  with  a CFL condition, in contrast to the
approach based one the RKF device. We mention that
   M\"uller and Stiriba \cite{MS} also combine  local time  stepping and
multiresolution for implicit schemes, and that
more details are also given in  the germinal papers of Berger
and Oliger \cite{BO}, Osher and Sanders \cite{OS} and the references
therein.

The remainder of this paper is organized as follows. In
Section~\ref{bbrs_sec:reac-dif}   the reaction-diffusion systems
 studied herein are briefly described. In
Section~\ref{bbrs_sec:fv} we present the reference finite volume methods to
which we apply the multiresolution device, and in Section~\ref{bbrs_sec:MR}
  the adaptive multiresolution strategy, as well
as  the required graded tree data structure, are outlined.
 In Section~\ref{bbrs_sec:err-analysis} we analyze the
 error   of the adaptive multiresolution scheme, and
deduce the optimal choice of the threshold. This choice
 ensures  that the discretization
 error of the reference scheme is balanced with the accumulated thresholding
 error which is introduced in each time step.
In Section~\ref{bbrs_sec:LTS} we
address the local time stepping strategy applied to the multiresolution
strategy, and in  Section~\ref{bbrs_sec:algimp}, we outline the
overall multiresolution procedure.
Finally, in Section~\ref{bbrs_sec:results}  the method is applied to
different scenarios. Example~1 corresponds to a single-species reaction
diffusion equation, Examples~2 and~3 deal with
 the thermo-diffusive model for the
interaction between flame balls, Examples~4 and~5  shows the results
of Turing-type  pattern formation produced by a reaction-diffusion
system, and  Example~6 arises from  a model
of chemotaxis with growth. Conclusions of our study are collected in
Section~\ref{bbrs_sec:conc}.
 All  numerical results clearly reveal high
resolution and improvement in terms of compression of memory
and  savings in computational effort.

\section{A class of reaction-diffusion systems} \label{bbrs_sec:reac-dif}
\subsection{A single-species reaction-diffusion model}
Model~1 is  the following initial-boundary value problem
for a scalar reaction-diffusion equation,
 where $\mathbf{x} =(x,y)$ and
 $(x,y,t)\in  Q_T := \Omega \times [0,T]$, $\Omega\subset \RR^2$:
\begin{subequations}\label{bbrs_eq:reac-dif}
\begin{align}
u_t& =f(u, \x)+\Delta A(u), \label{bbrs_eq2.1a} \\
 u(\x,0 ) & = u_0(\x) \quad\text{on $\Omega$}, \label{bbrs_model1init} \\
 \nabla A (u) \cdot \mathbf{n} & =0 \quad
\text{on $\Sigma_T:=\partial\Omega\times[0,T]$.}
 \label{bbrs_model1bound}
\end{align}
\end{subequations}
This  problem may serve as a scalar prototype degenerate
reaction-diffusion  model.  Here, the zero-flux boundary
condition \eqref{bbrs_model1bound} implies that the reaction-diffusion domain is
isolated from the external environment. For $f(u,\x) = f(u)$,
  \eqref{bbrs_eq2.1a}  appears in   \cite{Murray}
 in an  ecological setting,
 where $u$ denotes the population density of a species,
 and  $f(u)$ is its dynamics, where
it is assumed that   $f(0)=0$ and  $f'(0)\neq 0$.
 For example, $f(u) = u(1-u) -u^2/(1+u^2)$
 corresponds to the  population
 dynamics of the spruce band-worm  \cite{Murray},
 and   models the growth of the population  by
  a logistic expression and the rate
of mortality due to predation by other animals.
 We modify this expression by a radial spatial  factor,
and use
\begin{align} \label{bbrs_fxueq}
f(u,\x):=10\left(\exp(-5r)u(1-u)+\bigl(\exp(-5r)-1\bigr)\frac{u^2}{1+u^2}\right),
 \quad r:=\sqrt{(x-0.5)^2+(y-0.5)^2},
\end{align}
which means that the birth of individuals is concentrated near the
center $(0.5,0.5)$, and mortality increases with increasing
 distance   from the origin.
On the other hand, most  standard  spatial  models of population
dynamics simply assume  that  $A(u) = Du$, where the
 constant   diffusion coefficient
$D>0$   measures the dispersal efficiency of the species under
consideration. Motivated by  Witelski \cite{Witelski},
 who advanced degenerate diffusion in the context of population dynamics,
 we utilize herein the strongly degenerate diffusion coefficient
\begin{align}  \label{bbrs_model1audef}
 A(u)=\begin{cases}
0& \text{for $u\leqslant u_{\mathrm{c}}$,} \\
D(u-u_{\mathrm{c}})& \text{otherwise,}\end{cases}
\end{align}
where $u_{\mathrm{c}} >0$ is an assumed critical (threshold)
 value of~$u$ beyond which diffusion will take place.
 Model~1 gives rise to Example~1 of  Section~\ref{bbrs_sec:results}.

The difficulty in the well-posedness analysis
 of the problem \eqref{bbrs_eq:reac-dif}
 lies in  the boundary condition \eqref{bbrs_model1bound}
 when  $A$ is strongly degenerate. It is quite difficult to
 give a correct formulation of the zero flux boundary conditions.
 For  the case of non-homogeneous Dirichlet boundary conditions,
 however,   Mascia, Porretta, and Terracina
\cite{Porretta:01} demonstrated existence and uniqueness of
 $L^{\infty}$ entropy solutions.
In the special case where the function $A$ is strictly increasing, the
classical framework of variational solutions of parabolic equations is
 sufficient  to satisfy this wish.

\subsection{A two-species reaction-diffusion model}\label{bbrs_sec:2species}
Model~2 is given  by the following
initial-boundary value problem for a reaction-diffusion system
 on $Q_T$:
\begin{subequations}\label{bbrs_turing-sys}
\begin{align}
u_t&=\gamma f(u,v)+\Delta A(u) \quad \text{on $Q_T$,} \label{bbrs_turingu}  \\
v_t&=\gamma g(u,v)+d\Delta B(v) \quad \text{on $Q_T$,} \label{bbrs_turingv}
\\
u(\x,0)&= u_0(\x),\quad v(\x,0)=v_0(\x)\quad  \text{for $\x\in \Omega$}, \\
\nabla A(u)\cdot\n&=\nabla B(u)\cdot\n=0 \quad
  \text{on $\Sigma_T$}. \label{bbrs_neumcond}
\end{align}
\end{subequations}
We study this system in the context of two applications, namely as a
 model of combustion and as a two-species model of mathematical
 biology.

For  $A(u)=B(u)=u$,
$d=1/\mathrm{Le}$
 and  $\gamma=1$,  \eqref{bbrs_turing-sys} represents
  a reduced dimensionless thermo-diffusive model
 describing a combustion process, where $\mathrm{Le}$
is  the Lewis number.
We restrict ourselves to  a simple chemical
reaction with only two reactants and one product, the first reactant and the
product being highly diluted in the second reactant; and we neglect
gravity.
 Since the chemical reaction takes place in a lean premixed
 gas, we focus on the
limiting reactant, and denote by~$v$ its normalized partial mass, while $u$
represents  normalized temperature.
 The reaction rates  are given by an Arrhenius law:
\begin{align} \label{bbrs_arrh}
f(u,v):=\frac{\beta^2}{2}v\exp\left(\frac{\beta(1-u)}{\alpha(1-u)-1}
\right), \quad  g(u,v):=-f(u,v),
\end{align}
where   $\alpha$ and~$\beta$ are the temperature rate and the
 dimensionless activation energy, called Zeldovich number, respectively.
 In Example~2
 of  Section~\ref{bbrs_sec:results},  this model is employed
 to simulate  the interaction between two flame balls,
as an extension of the applications of the same model that were
considered in  \cite{RS05,ctm_rs04}. Here, a   \emph{flame ball}
denotes a slowly
propagating spherical flame structure in a premixed gaseous mixture.

If radiation effects are taken into account, \eqref{bbrs_turingu} is replaced by
\begin{align} \label{bbrs_eq24anew}
u_t=\gamma f(u,v)+S(u)+\Delta A(u) \quad \text{on $Q_T$,}
\end{align}
where the dimensionless heat loss due to radiation $S$ follows the
Stefan-Boltzmann law
\begin{align} \label{bbrs_stefanboltz}
S(u)=\rho\bigl[ (u+\alpha^{-1} -1)^4-
(\alpha^{-1} -1)^4\bigr],
\end{align}
and the  dimensionless coefficient $\rho$ controls the radiation level.
Conditions  \eqref{bbrs_neumcond} imply that the process takes
place inside a box with adiabatic walls. See \cite{ctm_rs04} for
 details and a discussion of the case with one flame ball.
The interaction  of two flame balls including radiation  is simulated
 in Example~3 of  Section~\ref{bbrs_sec:results}.

On the other hand,   \eqref{bbrs_turing-sys} also arises
 in mathematical biology as a well-known reaction-diffusion system
modelling
the interaction between two chemical species with respective
concentrations~$u$ and~$v$. Under certain conditions,
 it produces stationary  solutions with  Turing-type spatial patterns
\cite{Murray,turing}. To simulate the formation of such a pattern,
we here select the  kinetics between each
species due to  Schnakenberg \cite{schnakenberg}:
\begin{align} \label{bbrs_fgkinetics}
f(u,v)=a-u+u^2v,\quad g(u,v)=b-u^2v.
\end{align}
Alternative choices of~$f$ and~$g$ that lead to Turing-type patterns
 are discussed in \cite{murray82,Murray}.
 For
\begin{align} \label{bbrs_aubu}
A(u)=B(u)=u,
\end{align}
this system has a uniform positive steady state $(u^0,v^0)$
 given by
$u^0=a+b$ and  $v^0=b/(a+b)^2$, where $b>0$ and $a+b>0$,
and under certain conditions, \eqref{bbrs_turing-sys}  has a unique solution.
See for instance  \cite{britton} for the proof of existence
and uniqueness.

We recall from \cite[Sect.~2.3]{Murray} some  results
 on  the conditions under which \eqref{bbrs_turing-sys} produces
a diffusion-driven instability giving rise to Turing-type pattern
  in the non-degenerate case.
 A necessary condition
 is satisfaction  of the inequalities
$f_u + g_v < 0$,  $f_u g_v - f_v g_u >0$,
 $d f_u + g_v >0$ and $( d f_u + g_v)^2 - 4d
 (f_u g_v - f_v g_u) >0$.
Evaluating these inequalities for the system \eqref{bbrs_turingu},
\eqref{bbrs_turingv} and the particular kinetics \eqref{bbrs_fgkinetics} yields the
inequalities
\begin{align} \label{bbrs_abineq}
 0 < b-a < (a+b)^3, \quad (a+b)^2 >0, \quad
 d(b-a) > (a+b)^3, \quad
\bigl(d(b-a)-(a+b)^3\bigr)^2 > 4d(a+b)^4.
\end{align}
To characterize the stationary pattern that arises from a choice of
$(a,b)$ that satisfies \eqref{bbrs_abineq}, we define
\begin{align*}
 L^{\pm} (a,b,d):= \frac{d(b-a) - (a+b)^3 \pm  \sqrt{[d(b-a)-(a+b)^3]^2 -
4d (a+b)^4} }{2d(a+b)}.
\end{align*}
The analysis of general rectangular domains
 \cite{Murray}  implies  that in the non-degenerate case,
 the unstable patterned solution of
 the initial-boundary value problem \eqref{bbrs_turing-sys} is given by
\begin{align*}
 \mathbf{w} (x,y,t) = \sum_{m,n} \mathbf{c}_{nm} \exp \bigl(
 \lambda(k^2)
t \bigr) \cos (n \pi x) \cos (m \pi y),
\end{align*}
where the constants $ \mathbf{c}_{nm}$ depend on a Fourier series
 of the initial conditions for $\mathbf{w}$, and the summation takes
 place over all numbers $n$ and $m$ that satisfy
\begin{align*}
 \gamma L^- (a,b,d) =: k_1^2 < k^2= \pi^2 (n^2 + m^2) < k_2^2
 := \gamma L^+(a,b,d),
\end{align*}
and $\lambda(k^2)$ is the positive solution of the following
 equation,  where $f_u$, $f_v$, $g_u$ and $g_v$ are evaluated at
$(u^0,v^0)$:
\begin{align*}
 \lambda^2 + \lambda \bigl( k^2 (1+d) - \gamma (f_u + g_v) \bigr) +
 dk^4 - \gamma ( d f_u + g_v) k^2 + \gamma^2 (f_u g_v-f_v g_u) =0.
\end{align*}
Example~4 of Section~\ref{bbrs_sec:results} presents a  numerical
solution of \eqref{bbrs_turing-sys} with the kinetics \eqref{bbrs_fgkinetics}
 and the diffusion coefficients \eqref{bbrs_aubu}, where parameters
are chosen according to the preceding discussion
 such that indeed a Turing-type  pattern is produced.
On the other hand, in Example~5, we  present a simulation where
\eqref{bbrs_aubu} is replaced by the degenerate diffusion functions
 \begin{align}\label{bbrs_degen_param}
A(u)=\begin{cases}
0& \text{for }u\leqslant u_{\mathrm{c}},\\
u-u_{\mathrm{c}}& \text{otherwise}\end{cases},\ B(u)=\begin{cases}
0& \text{for }u\leqslant v_{\mathrm{c}},\\
u-v_{\mathrm{c}}& \text{otherwise,}\end{cases} \quad
 u_{\mathrm{c}}, v_{\mathrm{c}} \geq 0.
\end{align}
It turns out that even if the stability analysis does {\em not}
 apply to the degenerate case,  our numerical experiments (Example 5)
lead to  the formation of a pattern.

We mention that the mathematical analysis of the system \eqref{bbrs_turing-sys}
is still an open problem because of the presence of the boundary condition \eqref{bbrs_neumcond}.
A successful technique for proving uniqueness of (entropy weak) solutions to
degenerate parabolic equations with Dirichlet boundary condition
is based on Kru\v{z}kov's method \cite{Kruzkov}.
Related to this approach we mention that
 Holden, Karlsen, and Risebro  \cite{holden2003}   prove
existence and uniqueness of entropy solutions of weakly coupled systems of
degenerate parabolic equations in an unbounded domain. Our system is an example of the degenerate
reaction-diffusion system analyzed in \cite{holden2003},
 but is equipped here with the boundary condition \eqref{bbrs_neumcond}.

\subsection{A chemotaxis-growth system}
We assume that $\Omega\subset \RR^2$ is convex, bounded and open.
 Model~3 is  the following generalization of
 the Keller-Segel model \cite{horstmann2,kellersegel1} for chemotactical movement:
\begin{subequations}\label{bbrs_eq:Kel-Seg}
\begin{align}
u_t& = \nabla \cdot  \bigl(\sigma\nabla u- u \nabla \chi(v) \bigr)+g(u)
 \quad  \text{on $Q_T$}, \\
v_t & =h(u,v) + d \Delta v  \quad  \text{on  $Q_T$,} \\
u(\x,0)&=u_0(\x),\quad  v(\x,0)=v_0(\x)  \quad \text{on $\Omega$,} \\
 \nabla u\cdot\n &=\nabla v\cdot\n=0 \quad \text{on $\Sigma_T$}.
\end{align}
\end{subequations}
The system \eqref{bbrs_eq:Kel-Seg}
describes the aggregation of slime molds caused by their
chemotactical features.
Cell migration appears in numerous
  biological phenomena. In the case of chemotaxis, cells
(or an organism) move in response to a chemical gradient. Specifically,
\eqref{bbrs_eq:Kel-Seg} corresponds to the model proposed by Mimura and
 Tsujikawa   \cite{mimura} for
 the spatio-temporal aggregation patterns shown by
the bacteria {\it Escherichia coli}. This model incorporates four
effects: diffusion, chemotaxis, production of chemical substance,
and population growth.  In the absence of growth, $u=u(\x,t)$
usually represents the density of the cell population of the amoebae
{\it Dictyostelium discoideum}, $v=v(\x,t)$ is the
 concentration of the chemo-attractant
({\it cAMP: cyclic Adenosine Monophosphate}), and $\chi$
denotes the chemotactical sensitivity function, which may  be given by
\begin{align} \label{bbrs_chinuv}
 \chi (v) = \nu v, \quad \nu >0,
\end{align}
where $\nu$ is a chemotactical parameter.
The function~$g$  takes into account
the growth rate of the population, and can be given by
\begin{align} \label{bbrs_gudef}
 g(u) = u^2 (1-u).
\end{align}
 Moreover,
  $\sigma>0$ and and $d>0$~are constant diffusion rates
for both components. The function~$h$ describes the rates of
production and degradation of the chemo-attractant; here, we choose
\begin{align} \label{bbrs_huvdef}
h(u,v)=\alpha u-\beta v, \quad  \alpha,\beta\geqslant 0.
\end{align}
For this case it is known that if $0\leqslant u_0\in L^2(\Omega)$,
$0\leqslant v_0\in H^{1+r}(\Omega)$, and  $\partial \Omega$ is smooth
enough, \eqref{bbrs_eq:Kel-Seg} possesses a unique global solution (see, e.g.,
 \cite{biler}); and if $u_0$ and  $v_0$ are radial and
$\smash{\|u_0\|_{L^1}}$
is sufficiently large, then  $\smash{\|u(t)\|_{L^2}}$ blows up in finite time.
On the other hand, Efendiev, Kl\" {a}re, and Lasser  \cite{efendiev}
 analyzed the fractal dimension  of the exponential attractor
in dependence of~$\nu$. Our Example~6 of
Section~\ref{bbrs_sec:results} is based on examples presented
 in  \cite{efendiev}, and presents numerical
 solutions of \eqref{bbrs_eq:Kel-Seg}  for various values of~$\nu$.

\section{Finite Volume Schemes}\label{bbrs_sec:fv}
We employ a standard finite volume
 scheme to  discretize a reaction-diffusion equation, which is
described here for   a uniform
 grid. The   rectangular  spatial
 domain $\Omega \subset \mathbb{R}^2$ is partitioned into  control
volumes $\smash{ (\Omega_{ij})_{(i,j)\in\Lambda}}$, where $\Lambda$
 is an index set, defining
$\smash{\Omega_{ij}:=[x_{i-1/2},x_{i+1/2}]\times[y_{j-1/2},y_{j+1/2}]}$,
$\smash{\Delta x:=x_{i+1/2}-x_{i-1/2}}$, $\smash{\Delta y:=y_{j+1/2}-y_{j-1/2}}$, for all $(i,j)\in\Lambda$, and
$\smash{\widetilde{\Delta x}}:=\min \{\Delta x,\Delta y\}$.
 The  cell
average of  a quantity~$q$  at time~$t$ is defined by
\begin{equation}
\bar{q}_{ij}(t)=\frac{1}{|\Omega_{ij}|}\iint_{\Omega_{ij}}q(\x,t) \, d\x.
\end{equation}
\subsection{Discretization of Models~1 and~2} The  finite
 volume  scheme is described here for \eqref{bbrs_eq:reac-dif} and as it applies
 to the first equation of \eqref{bbrs_turing-sys}; for the second
equation of \eqref{bbrs_turing-sys}, we replace~$u$ by~$v$,
$f(u,v)$ by~$g(u,v)$,  and~$A(u)$ by~$dB(v)$. Integrating the respective
equation and averaging over $\Omega_{ij}$ yields
\begin{align} \label{bbrs_eq2.1}
\frac{1}{|\Omega_{ij}|}\iint_{\Omega_{ij}} u_t (\x,t)\, d\x=
\frac{1}{|\Omega_{ij}|}\iint_{\Omega_{ij}}\mathcal{D}\Bigl(u(\x,t),\nabla
A\bigl(u(\x,t)\bigr)\Bigr)\, d\x\ +\frac{1}{|\Omega_{ij}|}
\iint_{\Omega_{ij}}f(u(\x,t))\, d\x,
\end{align}
where $\mathcal{D}$ denotes the right-hand side of the PDE under
 consideration except for the reaction term. For the two-dimensional
  case and on a cartesian grid, $\mathcal{D}$
is discretized via
\begin{gather*}
\bar{\mathcal{D}}_{ij}:=-\frac{1}{\Delta x}\left(\bar{F}_{i+1/2,j}-
\bar{F}_{i-1/2,j}\right)-\frac{1}{\Delta y}\bigl(\bar{F}_{i,j+1/2}-
\bar{F}_{i,j-1/2}\bigr), \\
\bar{F}_{i+1/2,j}:=-\frac{1}{\Delta x}
\bigl( A(\bar{u}_{i+1,j})-A(\bar{u}_{ij}) \bigr), \quad\bar{F}_{i,j+1/2}:=-
 \frac{1}{\Delta y} \bigl( A(\bar{u}_{i,j+1})-A(\bar{u}_{ij})\bigr).
\end{gather*}
The reaction term is approximated by
$\smash{\bar{f}_{ij}\approx f(\bar{u}_{ij},\bar{v}_{ij})}$.
If we incorporate a  first-order Euler time integration for both components,
then the corresponding interior marching formula for Model~2 is
\begin{align}\label{bbrs_marchingf_model2}
\bar{u}_{ij}^{n+1}= \bar{u}_{ij}^n+\Delta t \gamma\bar{f}_{ij}+
\Delta t\bar{\mathcal{D}}_{ij}\bigl(\mathcal{S}(\bar{u}_{ij}^n),\widetilde{\Delta x}\bigr), \quad
\bar{v}_{ij}^{n+1} = \bar{v}_{ij}^n+\Delta t \gamma\bar{g}_{ij}+
d\Delta t\bar{\mathcal{D}}_{ij}\bigl(\mathcal{S}(\bar{v}_{ij}^n),\widetilde{\Delta x}\bigr),
\end{align}
where $\mathcal{S}(\cdot)$ denotes the  stencil utilized for computing
$\bar{\mathcal{D}}_{ij}$. According to \cite{ck,holden2003},
this scheme is stable under the CFL condition
\begin{equation}\label{bbrs_cfl-turing}
\lambda\gamma\big(\|f_u\|_\infty+\|f_v\|_\infty+\|g_u\|_\infty+
\|g_v\|_\infty\big)+4\mu d \bigl(\|A'\|_\infty+\|B' \|_\infty
\bigr)\leqslant 1.
\end{equation}
Here
$\smash{\lambda:=\Delta t/ \widetilde{\Delta x}}$,
 $\smash{\mu:=\Delta t / \smash{\widetilde{\Delta x}}^2}$.

\subsection{Discretization of Model~3}
We  define the difference
 operators $\delta_x^{\pm} V_{ij} := \pm (V_{i\pm 1,j}- V_{ij})$ and
    $\delta_y^{\pm} V_{ij} := \pm (V_{i,j\pm 1}- V_{ij})$.
Then  a suitable second
order difference operator for a general term $\nabla\cdot(Q\nabla u)$ is
$$
\nabla \cdot (Q\nabla u)\approx \frac{1}{\Delta
  x^2}\delta_x^+ \bigl(Q_{i+1/2,j}\delta_x^-u_{ij}\bigr)
+\frac{1}{\Delta y^2}\delta_y^+\bigl(Q_{i,j+1/2} \delta_y^-u_{ij}\bigr).
$$
Integrating the corresponding equations, averaging over $\Omega_{ij}$
and discretizing yields  the  following interior marching formula:
\begin{align}
\label{bbrs_marching_chemo}
\begin{split}
\bar{u}_{ij}^{n+1}=&  \bar{u}_{ij}^n+
\frac{\sigma\Delta t}{\Delta x^2} \delta_x^+\delta_x^- \bar{u}_{ij}^n
+ \frac{\sigma\Delta t}{\Delta y^2} \delta_y^+\delta_y^-
\bar{u}_{ij}^n
 +\frac{\Delta t}{\Delta x^2}\big(\delta_x^+\bigl(Q_{i-1/2,j}^n  \delta_x^-
  \bar{v}_{ij}^n \bigr)
 \bigr) \\ &   +\frac{\Delta t}{\Delta y^2}\bigl(\delta_y^+ \bigl(
 Q_{i,j-1/2}^n
 \delta_y^- \bar{v}_{ij}^n\bigr) \bigr)+g\bigl(\bar{u}_{ij}^n\bigr),  \\
\bar{v}_{ij}^{n+1}=& \bar{v}_{ij}^n+\Delta t h\bigl(\bar{u}_{ij}^n,\bar{v}_{ij}^n\bigr)+
\frac{d\Delta t}{\Delta x^2} \delta_x^+\delta_x^- \bar{v}_{ij}^n
+ \frac{d\Delta t}{\Delta y^2} \delta_y^+\delta_y^- \bar{v}_{ij}^n, \\
Q_{i,j+1/2}^n & : = \frac{1}{2}\big(\chi'\bigl(\bar{v}_{ij}^n\bigr)\bar{u}_{ij}^n
+\chi'\bigl(\bar{v}_{i,j+1}^n \bigr)\bar{u}_{i,j+1}^n \big), \quad
Q_{i+1/2,j}^n  := \frac{1}{2}\big(\chi'\bigl(\bar{v}_{ij}^n \bigr)
\bar{u}_{ij}^n
+\chi'\bigl(\bar{v}_{i+1,j}^n \bigr)\bar{u}_{i+1,j}^n \bigr).
\end{split}
\end{align}
Analogously to \eqref{bbrs_cfl-turing}, the scheme
\eqref{bbrs_marching_chemo} is stable under the corresponding CFL condition
\begin{align}\label{bbrs_cfl-chemo}
\lambda\bigl(\|h_u\|_\infty+\|h_v\|_\infty+
\|g' \|_\infty\bigr)+4\mu d \bigl(\sigma+\|\chi' \|_\infty \bigr)\leqslant 1.
\end{align}
 The left-hand sides of
 \eqref{bbrs_cfl-turing} and \eqref{bbrs_cfl-chemo}  obviously
evolve in time, so in practice, at each time step we obtain $\Delta t$ from
these conditions, and~$\lambda$ and~$\mu$ are  not constants; rather,
 they are adjusted in each time step.

 \section{Conservative adaptive multiresolution discretization}\label{bbrs_sec:MR}
 In this section we recall some basic properties of the
 multiresolution discretization and the data structure.
 For a more detailed description, we refer  to \cite{brss2,RSTB03}.
  We will organize the numerical solution and corresponding differences at
 different levels, in a \emph{dynamic graded tree structure}: whenever
 a \emph{node} is included in the tree,  all other nodes corresponding to
 the same spatial region in coarser resolutions  are also included.
  The tree structure is mainly needed for ease of navigation,
  which contributes to accelerating the scheme; data compression would
  also be possible by other techniques.  The adaptive grid
 corresponds to a set of nested dyadic grids generated by refining
 recursively a given cell depending on the local regularity of the solution.

  \begin{figure}[t]
     \begin{center}
     \includegraphics[width=.5\textwidth]{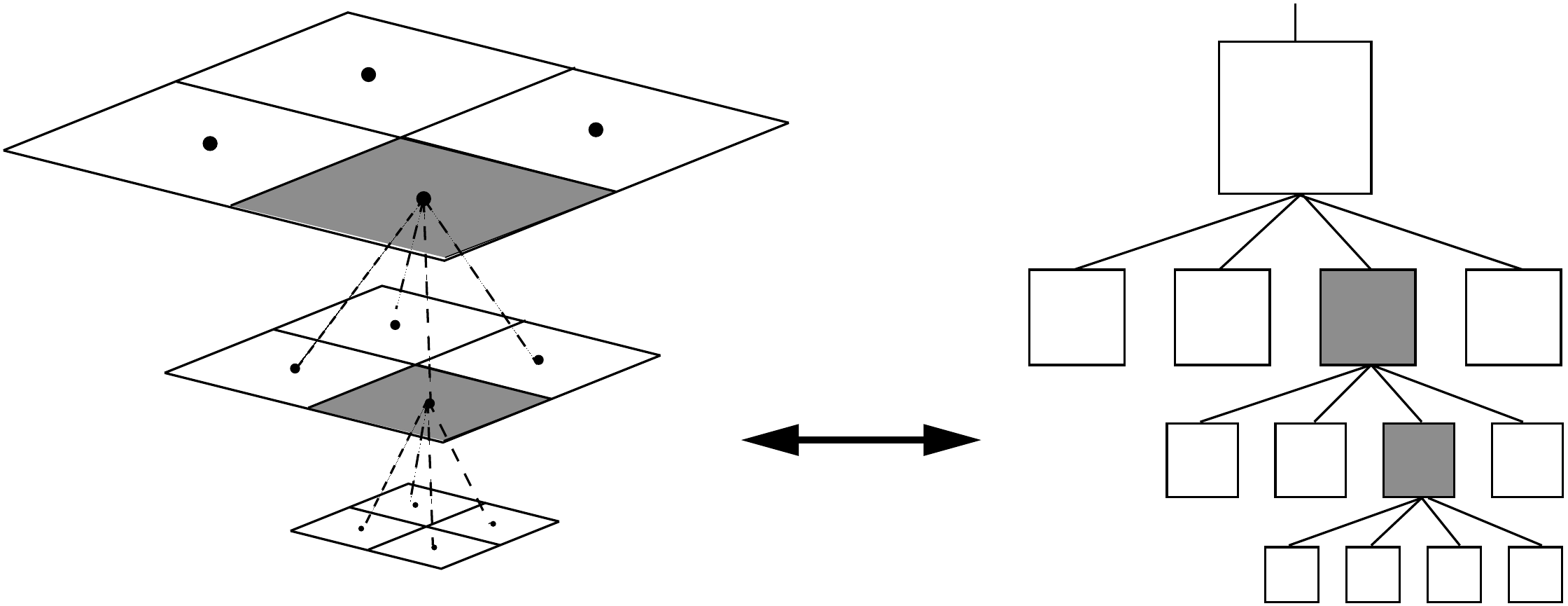}
     \end{center}
     \caption{\it Graded tree data structure (``quadtree''), after \cite{Muller}.}
     \label{bbrs_fig:arbol}
 \end{figure}

We  denote by  {\em root}  the basis   of the tree.
In two space dimensions, a parent node has four
  sons, and the sons of the same parent are called
 \emph{brothers}. A node without sons is called a \emph{leaf}.
 A given node has $s'=2$ nearest neighbors in each direction, called
 \emph{nearest cousins}, needed for the computation of the fluxes of leaves;
 if these nearest cousins do not
  exist,  we create them as \emph{virtual leaves}. Brothers are also considered
  nearest cousins.  Figure~\ref{bbrs_fig:arbol} illustrates the graded tree structure.
  The leaves of the tree are the control volumes from which we form the
adaptive mesh. Here, the property of the tree being graded means that grid
refinement and
coarsening is governed by that two neighboring control volumes cannot differ
by more than one level in the tree. This is equivalent with the notion of
"one-irregular rule" (see e.g., \cite{moore00,moore07}).
  We denote by $\Lambda$ the set
 of indices of existing nodes, by $\mathcal{L}(\Lambda)$ the
 restriction of $\Lambda$ to the leaves, and by $\Lambda_l$ the
 restriction of $\Lambda$ to a multiresolution level $l$, $0\leq
 l\leq L$. We  denote by
 $\bar{\ww}_{i,j,l}=(\bar{u}_{i,j,l},\bar{v}_{i,j,l})$  the vector of cell
 averages for both components of the solution (and the obvious simplification
 $\bar{\ww}_{i,j,l}=\bar{u}_{i,j,l}$ for a single-species problem) located at
 spatial position  $(i,j)$ at level $l$, and by $\smash{\bar{\mathbf{W}}_l}$ the set
 of cell averages for all nodes at level $l$.
 To estimate the cell averages of a level $l$ from those of the
 next finer level $l+1$, we use the projection operator
 $P_{l+1\to l}$, which is exact,  unique, and in our case is defined by
 \begin{align}\label{bbrs_multisol2d}
 \bar{\ww}_{i,j,l}=\frac{1}{4}\sum_{e_1,e_2\in \{ 0,1\}
 }\bar{\ww}_{2i+e_1,2j+e_2,l+1}.
 \end{align}
 To estimate the cell averages of a level $l+1$ from those of level $l$, we
 employ the prediction operator  $P_{l\to l+1}$,  which provides an approximation
 $\hat{\mathbf{w}}$ by  interpolation of $\smash{\bar{\mathbf{W}}_l}$ at level $l+1$.
 This operator is  local
 in the sense that the interpolation for a son is made from the cell averages
 of its parent and the $s$~nearest cousins of its parent; and it is
   consistent  with the projection in the sense that it is conservative
 with respect to the coarse grid cell averages or equivalently,
 $P_{l+1\to l}\circ P_{l\to l+1}= \mathrm{Id}$.
 For a regular grid structure in two dimensions,
  we use a polynomial interpolation
 introduced in  \cite{bh97}:
 \begin{equation}
 \hat{\ww}_{2i+e_1,2j+e_2,l+1}=\bar{\ww}_{i,j,l}-(-1)^{e_1}Q_x-(-1)^{e_2}Q_y
 +(-1)^{e_1e_2}Q_{xy},\quad e_1,e_2\in \{ 0,1\},
 \end{equation}
 where
 \begin{align} \label{bbrs_papa2D}
 \begin{split}
 Q_x&:=\sum_{n=1}^s\tilde{\gamma}_n(\bar{\ww}_{i+n,j,l}-\bar{\ww}_{i-n,j,l}),\quad
 Q_y:=\sum_{p=1}^s\tilde{\gamma}_p(\bar{\ww}_{i,j+p,l}-\bar{\ww}_{i,j-p,l}), \\
 Q_{xy}&:=
 \sum_{n=1}^s\tilde{\gamma}_n\sum_{p=1}^s\tilde{\gamma}_p(\bar{\ww}_{i+n,j+p,l}-
      \bar{\ww}_{i+n,j-p,l}-\bar{\ww}_{i-n,j+p,l}+\bar{\ww}_{i-n,j-p,l}).
 \end{split}
 \end{align}
 The chosen accuracy order of the multiresolution method for our cases
 is $r=s+1=3$, where $s$ is the number of required nearest uncles for each
 spatial direction. The
 corresponding coefficients are $\tilde{\gamma}_1=\smash{-\frac{22}{128}}$ and
 $\tilde{\gamma}_2=\smash{\frac{3}{128}}$.
 Nevertheless, one may select here an  arbitrarily higher  order of  accuracy.

As stated before, the adaptive grid consists in the set of leaves
$\mathcal{L}(\Lambda)$, which forms a partition of the computational domain~$\Omega$.

 A \emph{detail} is the difference between the exact and the predicted value
 $$\bar{d}^u_{i,j,l}:=\bar{u}_{i,j,l}-\hat{u}_{i,j,l},\quad
 \bar{d}^v_{i,j,l}:=\bar{v}_{i,j,l}-\hat{v}_{i,j,l}.$$

 For \emph{multicomponent solutions}, there are  many
 possible definitions of a scalar detail $\smash{\bar{d}_{i,j,l}}$
that   is calculated from the details of the components
  (in our case, $\smash{ \bar{d}^u_{i,j,l}}$ and
 $\smash{ \bar{d}^v_{i,j,l}}$),
  depending mainly on the nature of the
  problem.  Roussel and Schneider \cite{RS05}  utilize the Euclidian norm
  of the details,
 $\smash{\bar{d}_{i,j,l}=((\bar{d}^u_{i,j,l})^2+(\bar{d}^v_{i,j,l})^2)^{1/2}}$.
 In our
 case, given the nature of our problems, we could simply select one
  component
 $\smash{\bar{d}_{i,j,l}=\bar{d}^u_{i,j,l}}$, as was done
 by Sj\"{o}green \cite{sjoeg} for the compressible Euler equations,
 but in order to guarantee that the computations
 of the refinement and coarsening procedures are always on the safe side, in the
 sense
 that we always prefer to keep a node with a detail pair containing at least
one value over the threshold
\begin{align} \label{bbrs_eq4.4}
 \eps_l=2^{2(l-L)}\eps,
\end{align}
  we will use
 $\smash{\bar{d}_{i,j,l}=\min\{\bar{d}^u_{i,j,l},\bar{d}^v_{i,j,l}\}}$
   and
 $\smash{\bar{d}_{i,j,l}=\max\{\bar{d}^u_{i,j,l},\bar{d}^v_{i,j,l}\}}$
 for the refinement and coarsening procedures (see details in Algorithm~7.1),
 respectively,   similar to Harten's choice  \cite{harten:1995} for the Euler
  equations of gas dynamics.

 Since  a parent has four sons, the consistency property of
 the prediction operator implies that  knowledge of the cell average
 values of the four sons is equivalent to that of the cell average
 value of the parent node and three independent details. Repeating this
 operation recursively on $L$ levels, we get the
 \emph{multiresolution transform} on the cell average values
 $\smash{\bar{\mathbf{M}}:
 \bar{\mathbf{W}}_L\mapsto(\bar{D}_L,\ldots,\bar{D}_1,\bar{\mathbf{W}}_0)}$,
 where $\smash{\bar{D}_l=(\bar{d}_{i,j,l})_{(i,j)}}$.
 This means that  knowledge of
 the cell averages of all leaves is equivalent to that  of the
 cell average  value of the root and of the details of all nodes of the tree.

 After thresholding, i.e., removing nodes where the detail is below the prescribed
 tolerance $|\bar{d}_{i,j,l}|<\eps_l$, where $\eps_l$ is given by \eqref{bbrs_eq4.4};
 a \emph{safety zone} is added  to the tree, which means that
 one finer level is added to the tree in all possible nodes without violating the
 graded tree data structure. This is done by splitting each leaf into four new leaves
 in such a way that the new tree remains graded. This device, which was
  proposed e.g. in \cite{harten:1995,RSTB03},   ensures that the graded tree
 will represent adequately the solution in the next time step, and depends strongly
on the assumption of finite propagation speed of sharp fronts.

\begin{figure}[t]
     \begin{center}
      \includegraphics[width=.3\textwidth]{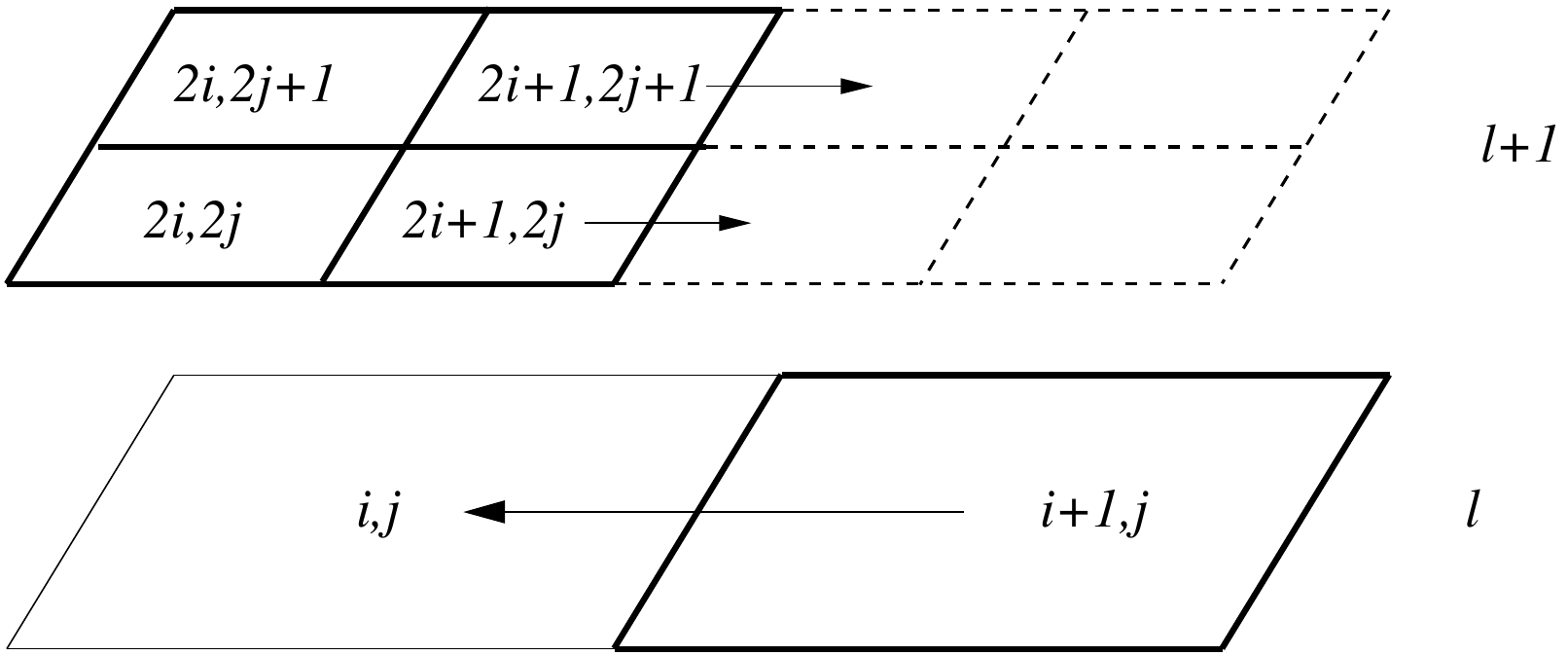}
     \end{center}
     \caption{\it Conservative flux computation for coarser levels.}
     \label{bbrs_fig:fluxes_2D}
 \end{figure}

 Now, assume that the tree has only two levels $l$ and $l+1$ (straightforwardly
 extensible to an arbitrarily larger tree).
 To ensure conservativity of the scheme, we compute only the
  fluxes at level $l+1$ and we set the ingoing flux on the leaf
 at level $l$ equal to the sum of the outgoing fluxes on the leaves of
 level $l+1$ (see Figure~\ref{bbrs_fig:fluxes_2D})
 \begin{equation}\label{bbrs_flux_conserv}
 F_{(i+1,j,l)\to(i,j,l)}=F_{(2i+1,2j,l+1)\to(2i+2,2j,l+1)}+
 F_{(2i+1,2j+1,l+1)\to(2i+2,2j+1,l+1)}.
 \end{equation}
 This choice decreases the number of costly flux evaluations without
 loosing the conservativity in the flux computation, and this presents a real
 advantage when using a graded tree structure, see
 e.g. \cite{RSTB03}. This advantage is lost for a non-graded
tree structure, in which case these data
 (fluxes for leaves on an immediately finer level) are not always available.

 It is important to remark that in the case of systems, since we are
 dealing with reaction-diffusion systems where species tend  to attract
 each other, we  manage the multiresolution framework and the data
 structure as  one unified mesh with two components per control volume.
  This means that we
   construct only one graded tree and apply only one thresholding strategy
 for both species; however, there are other cases where it is preferable to
 organize  the different species in separate adaptive meshes, for
 example  when the species segregate spatially, as
   in systems of conservation laws modelling traffic
flow and polydisperse sedimentation \cite{bk06}.

\section{Error analysis of the adaptive
multiresolution scheme}\label{bbrs_sec:err-analysis}
Using the main properties of the reference finite volume scheme
on a uniform grid at the finest level $L$,
such as the contraction property in $L^1$ norm, CFL stability condition and
order of approximation in space, in \cite{CKMP2002,RSTB03},
the authors decompose the global error between the cell average
values of the exact solution vector at the level $L$, denoted by
$\bar{\ww}^L_{\mathrm{ex}}=(\bar{u}^L_{\mathrm{ex}},\bar{v}^L_{\mathrm{ex}})$,
and those of the multiresolution computation with
a maximal level $L$, denoted by $\bar{\ww}^L_{\mathrm{MR}}$, into two errors
\begin{align}
\bigl\|\bar{\ww}^L_{\mathrm{ex}} - \bar{\ww}^L_{\mathrm{MR}} \bigr\| \leq
\bigl\|\bar{\ww}^L_{\mathrm{ex}} - \bar{\ww}^L_{\mathrm{FV}} \bigr\| +
\bigl\|\bar{\ww}^L_{\mathrm{FV}} - \bar{\ww}^L_{\mathrm{MR}} \bigr\|.
\end{align}
The first error on the right-hand side, called \emph{discretization
  error},
is the one of the reference finite volume scheme on a uniform grid at the
finest level $L$. In many circumstances (see e.g. \cite{brss,CKMP2002}),  this error can
be bounded by
\begin{align}
\bigl\|\bar{\ww}^L_{\mathrm{ex}} - \bar{\ww}^L_{\mathrm{FV}}
\bigr\| \leq C_1  2^{-\tilde{\alpha} L},   \quad  C_1 > 0,
\end{align}
where $\tilde{\alpha}$ is the convergence order of the finite volume scheme. Based on
preliminary numerical experiments (obtained in a similar fashion as in \cite{brss}),
for our examples we obtain the approximate value $\tilde{\alpha} = 2.18$.

For the second error, called \emph{perturbation error},
in \cite{CKMP2002} the authors assume that  the
details on a level~$l$ are
deleted when they are smaller than a prescribed tolerance $\epsilon_l$.
Under this assumption, they show that  if the discrete time evolution
operator is contractive in the chosen norm, and if the tolerance
$\epsilon_l$ at the level $l$ is given by
 \eqref{bbrs_eq4.4},
then the perturbation error accumulates in time and  satisfies
\begin{align} \label{bbrs_equ:accu}
\bigl\|\bar{\ww}^L_{\mathrm{FV}} - \bar{\ww}^L_{\mathrm{MR}} \bigr\| \leq
C_2  n \epsilon, \quad  C_2 > 0,
\end{align}
where $n$ denotes the number of time steps. At a fixed time $T=n  \Delta t$,
this gives
\[
\bigl\|\bar{\ww}^L_{\mathrm{FV}} - \bar{\ww}^L_{\mathrm{MR}} \bigr\| \leq
C_2
\frac{T}{\Delta t} \epsilon,   \quad  C_2 > 0.
\]
Motivated by their analysis, and according to the global CFL condition 
\eqref{bbrs_cfl-turing} of the reference finite volume scheme defined for the
discretization of Model~2, we have
\[
\Delta t \leqslant \frac{ \widetilde{\Delta x}{}^2}{\|f_u\|_\infty +\|f_v\|_\infty+
\|g_u\|_\infty +\|g_v\|_\infty+ \widetilde{\Delta x}  4d(\|A'\|_\infty+\|B'\|_\infty)}.
\]
If we write $\smash{\widetilde{\Delta x}=\sqrt{|\Omega|}2^{-L}}$,  this yields
\[
\Delta t = C_3 \frac{|\Omega|2^{-2L}}{\|f_u\|_\infty +\|f_v\|_\infty+
\|g_u\|_\infty +\|g_v\|_\infty+
\sqrt{|\Omega|}2^{-L}4d(\|A'\|_\infty+\|B' \|_\infty)},
\quad 0<C_3\leqslant 1.
\]

The main idea of the adaptive multiresolution
scheme is to perturb the solution given by a finite volume scheme  on a
uniform discretization (reference mesh) in such a way that the total error, i.e., the error
between the exact solution and the adaptive solution that is projected to the reference
fine mesh, is of the same order as the discretization error. For this purpose,
one has to balance the discretization error and the perturbation error. 
With this idea in mind, it is possible to derive in a similar fashion the optimal choice
for the threshold parameter $\epsilon$ for the adaptive multiresolution scheme. 

As in \cite{brss}, we need that $\epsilon\propto 2^{-\tilde{\alpha}L}$ or
$$\epsilon |\Omega|2^{2L}\Bigl(\|f_u\|_\infty +\|f_v\|_\infty+
\|g_u\|_\infty +\|g_v\|_\infty+ \sqrt{|\Omega|}2^{-L}4d\bigl(\|A'\|_\infty+\|B'\|_\infty\bigr)\Bigr)
 \propto 2^{-\tilde{\alpha}L}.$$
 In this way, for the numerical computations of Model~2, we may set the so-called
\emph{reference tolerance} $\varepsilon_{\mathrm{R}}$   to
\begin{align} \label{bbrs_equ:epsref1}
\epsilon_{\mathrm{R}} = C \frac{2^{-(\tilde{\alpha}+2)L}}{|\Omega|(\|f_u\|_\infty +
\|f_v\|_\infty+\|g_u\|_\infty +\|g_v\|_\infty) +|\Omega|^{3/2}2^{L}\,
4d(\|A'\|_\infty+\|B' \|_\infty) }.
\end{align}
Analogously,  for Model~3, the reference tolerance may be set to
 \begin{align}\label{bbrs_equ:epsref2}
 \epsilon_{\mathrm{R}} &= C \frac{2^{-(\tilde{\alpha}_1+2)L}}{|\Omega|(\|h_u\|_\infty +
 \|h_v\|_\infty+\|g' \|_\infty) +|\Omega|^{3/2}2^{L}\,
 4d(\sigma+\|\chi'\|_\infty) },
\end{align}
where $\tilde{\alpha}_1$ is a value of the convergence rate for Model~3.

Note that all the $L^\infty$ norms in \eqref{bbrs_equ:epsref1} and \eqref{bbrs_equ:epsref2} are
computed numerically. To determine an acceptable value for the factor $C$ (which, of
course, depends on $T$, $C_1$, $C_2$ and $C_3$), a series of computations with different
tolerances are needed in each case, prior to final computations. Essentially, we  select the
largest available candidate value for $C$ such that  the same order of accuracy (same slopes for
the error computation) as  that of the reference finite volume scheme is maintained.
 This procedure basically generalizes the treatment in \cite{brss2} of spatially one-dimensional
strongly degenerate parabolic  equations. In \cite{DRS06} the authors prove for
 scalar, one-dimensional, nonlinear conservation laws, that the threshold error is stable in
 the sense that the constant $C$ is uniformly bounded and, in particular, does not depend on
 the threshold value $\epsilon_{\mathrm{R}}$, the number of refinement levels $L$ and the
 number of time steps $n$. In our case, even when a rigorous proof is still missing for the systems
 considered in the present work, from the previous deduction we see a similar behavior for $C$.

We also mention that as in previous works \cite{brss2,brss,CKMP2002}, here the reference 
tolerance remains fixed for all times, though it is certainly possible
 to recompute the reference tolerance at each time step.

\section{Local time stepping} \label{bbrs_sec:LTS}
We utilize a version of the locally varying time
stepping strategy advanced by M\"{u}ller and Stiriba
\cite{MS}, and summarize here its principles. The basic
idea is to enforce a local CFL condition by using the same CFL number for
all levels, and the strategy consists in evolving all leaves on level
$l$ using the local time step
\begin{equation}\label{bbrs_deltatlocal}
\Delta t_l=2^{L-l}\Delta t,\quad l=L-1,\ldots,0,
\end{equation}
where $\Delta t=\Delta t_L$ corresponds to the time step on the finest
 level $L$. This strategy allows to increase the time step for the major
part of the adaptive mesh without violating  the CFL stability  condition.

Clearly, portions of the solution lying on
different resolution levels need to be synchronized to obtain a conservative scheme.
This will be achieved after $2^l$ time steps using $\Delta t_l$: all
leaves forming the adaptive mesh are synchronized in time, so one time
step with $\Delta t_0$ is equivalent to $2^L$ intermediate time steps with
$\Delta t_L$. In order to additionally save  computational effort, we
  update  the tree structure (refinement and coarsening) only each odd
intermediate time step $1,3,\ldots,2^L-1$ (as suggested in
 \cite{BO}), and furthermore,  the projection and prediction
operators are performed only on the levels occupied by the leaves of
the current tree, i.e., we do not update the tree structure by prediction from
the root cell, but from the coarsest leaves (we  refer to this as
\emph{partial grid adaptation}). For the rest of the intermediate
time steps, we  use the current tree structure. Notice that the updating
of the tree is still done in each \emph{global} time step. For the sake of
synchronization and conservativity of the flux computation, for coarse
levels (levels without
leaves), we employ \emph{the same} fluxes ($\smash{\bar{\mathcal{D}}_{ij,l}}$ and
$\smash{\bar{f}_{ij,l}}$) computed in the previous intermediate time step,
because the cell averages on these levels are the same that in the previous
intermediate time step. Only for finer levels (levels containing leaves),
we compute fluxes, and do so in the following way: if there is a leaf at the
corresponding cell edge and at the same resolution level~$l$, we simply perform
a flux computation using the brother leaves, and the virtual leaves at the same
level if necessary; and if there is a leaf at the corresponding cell edge but
on a finer resolution level $l+1$ (in this case we refer to this edge as an
\emph{interface edge}), the flux will be determined as in \eqref{bbrs_flux_conserv},
that is, we compute the fluxes  at a level $l+1$ on the same edge, and we set the
ingoing flux on the corresponding edge at level~$l$ equal to the sum of the
outgoing fluxes of the son cells of level $l+1$ (for the same edge).
We recall that the graded tree structure ensures that two neighboring control volumes
of the adaptive mesh do not differ by more than one resolution level, which is
equivalent to the satisfaction of the one-irregular rule.
In order to always have at hand the computed fluxes as in \eqref{bbrs_flux_conserv},
we need to perform the locally varying time stepping recursively from
fine to coarse levels. If at any instance of the procedure there is a missing value,
we can project the value from the sons nodes or we can predict this value from
the parent nodes.  For illustrative purposes, we give an example of an interior
first-order flux calculation for Model~2, to complete a full macro time step,
by the following algorithm (we show only the flux calculation for 
$(i,j,l)\to (i+1,j,l)$ and $(i,j,l)\to (i,j+1,l)$ since 
the other fluxes are obtained analogously):

\begin{Alg}[Locally varying intermediate time stepping] \hfill
\begin{enumerate}
   \item Grid adaptation (provided the former sets of leaves and
          virtual leaves).
   \item {\bf do} $k=1,\ldots,2^L$ (and therefore the local time steps are
      $n+2^{-L},n+2\cdot 2^{-L}, n+ 3\cdot 2^{-L},\ldots,  n+1$)
   \begin{enumerate}
      \item Synchronization:
      \item[]{\bf do} $l=L,\ldots,1$
      \begin{itemize}
         \item[]{\bf do} $i=1, \dots, |\tilde{\Lambda}|_x(l)$, $j=1, \dots,
           |\tilde{\Lambda}|_y(l)$
         \begin{itemize}
            \item[] {\bf if} $1\leqslant l\leqslant \tilde{l}_{k-1}$
              {\bf then}
            \begin{itemize}
                \item[]  {\bf if} $(i,j,l)$ is a virtual leaf {\bf then}
		  \begin{itemize}
		    \item[] ${\displaystyle
                  \bar{F}^{n+k2^{-L}}_{(i,j,l)\to(i+1,j,l)} \leftarrow
               \bar{F}^{n+(k-1)2^{-L}}_{(i,j,l)\to(i+1,j,l)}\vphantom{\int} }$
                 \item[]  ${\displaystyle
                   \bar{f}_{ij,l}^{n+k2^{-L}} \leftarrow
                   \bar{f}_{ij,l}^{n+(k-1)2^{-L}}, \quad
                   \bar{g}_{ij,l}^{n+k2^{-L}} \leftarrow
                   \bar{g}_{ij,l}^{n+(k-1)2^{-L}}\vphantom{\int}}$
            \end{itemize}
		  \item[] {\bf endif}
            \end{itemize}
            \item[]{\bf else}
            \begin{itemize}
	        \item[] {\bf if} $(i,j,l)$ is a leaf {\bf then}
		  \begin{itemize}
                \item[] ${\displaystyle
                   \bar{f}_{ij,l}^{n+k2^{-L}} \leftarrow f
                   \bigl(\bar{u}_{ij,l}^{n+k2^{-L}},\bar{v}_{ij,l}^{n+k2^{-L}}
                   \bigr), \quad \bar{g}_{ij,l}^{n+k2^{-L}} \leftarrow g
                   \bigl(\bar{u}_{ij,l}^{n+k2^{-L}},\bar{v}_{ij,l}^{n+k2^{-L}}
                   \bigr)\vphantom{\int}}$
                \item[] {\bf if}
                     $(i+1,j,l)$ is a leaf {\bf or} $(i,j+1,l)$ is a leaf {\bf then} compute fluxes by
                    \item[] $\quad {\displaystyle  \bar{F}_{(i,j,l)\to(i+1,j,l)}
                        \leftarrow -\frac{1}{h(l)}\bigl(A(\bar{u}_{i+1,j,l})
                        -A(\bar{u}_{i,j,l})\bigr)\vphantom{\int\limits_x^x}}$
                    \item[] $\quad {\displaystyle  \bar{F}_{(i,j,l)\to(i,j+1,l)}
                        \leftarrow -\frac{1}{h(l)}\bigl(A(\bar{u}_{i,j+1,l})
                        -A(\bar{u}_{i,j,l})\bigr)\vphantom{\int\limits_x^x}}$
                \item[]{\bf endif}
		\item[] {\bf if}
                    $(2i+2,2j,l+1)$, $(2i+2,2j+1,l+1)$ are leaves (interface
                   edges) {\bf then} compute fluxes by
                     \item[] $\quad {\displaystyle\vphantom{\int}
                          \bar{F}_{(i,j,l)\to(i+1,j,l)}\leftarrow
                          \bar{F}_{(2i+2,2j,l+1)\to(2i+1,2j,l+1)}
                         +\bar{F}_{(2i+2,2j+1,l+1)\to(2i+1,2j+1,l+1)}
                          \vphantom{\int}}$
                     \item[] $\quad {\displaystyle\vphantom{\int}
                          \bar{F}_{(i,j,l)\to(i,j+1,l)}\leftarrow
                          \bar{F}_{(2i,2j+2,l+1)\to(2i,2j+1,l+1)}
                         +\bar{F}_{(2i+1,2j+2,l+1)\to(2i+1,2j+1,l+1)}
                          \vphantom{\int}}$
                \item[] {\bf endif}
            \end{itemize}
  \item[] {\bf endif}
            \end{itemize}
            \item[] {\bf endif}
         \end{itemize}
         \item[]{\bf enddo}
      \end{itemize}
      \item[]{\bf enddo}
   \item Time evolution:
      \item[]{\bf do} $l=1,\ldots,L$, $i=1, \dots, |\tilde{\Lambda}|_x(l)$, $j=1, \dots, |\tilde{\Lambda}|_y(l)$
      \begin{itemize}
         \item[] {\bf if} $1\leqslant l\leqslant \tilde{l}_{k-1}$ {\bf then}
             there is no evolution:
         \begin{itemize}
             \item[] ${\displaystyle \bar{u}_{ij,l}^{n+(k+1)2^{-L}}
                \leftarrow \bar{u}_{ij,l}^{n+k2^{-L}}, \quad
              \bar{v}_{ij,l}^{n+(k+1)2^{-L}}
                \leftarrow \bar{v}_{ij,l}^{n+k2^{-L}}}$
         \end{itemize}
         \item[]{\bf else}
         \begin{itemize}
            \item[] Interior marching formula only for the leaves $(i,j,l)$:
            \item[] ${\displaystyle \bar{u}_{ij,l}^{n+(k+1)2^{-L}}
                  \leftarrow \bar{u}_{ij,l}^{n+k2^{-L}}
                   +\gamma \Delta t_l \bar{f}_{ij,l}^{n+k2^{-L}}+
                  \Delta t_l\bar{\mathcal{D}}_{ij,l} \bigl( \mathcal{S}
                   \bigl(\bar{u}_{ij,l}^{n+k2^{-L}}\bigr),h(l) \bigr)
                   \vphantom{\int}}$
            \item[]  ${\displaystyle   \bar{v}_{ij,l}^{n+(k+1)2^{-L}}
                      \leftarrow  \bar{v}_{ij,l}^{n+k2^{-L}}
                     +\gamma \Delta t_l  \bar{g}_{ij,l}^{n+k2^{-L}}+
               d\Delta  t_l\bar{\mathcal{D}}_{ij,l} \bigl( \mathcal{S}
              \bigl(\bar{v}_{ij,l}^{n+k2^{-L}}\bigr),h(l) \bigr) }$
            \end{itemize}
         \item[] {\bf endif}
      \end{itemize}
      \item[]{\bf enddo}
      \item Partial grid adaptation each odd intermediate time step:
      \item[] {\bf do} $l=L,\ldots,\tilde{l}_k+1$
	\begin{itemize}
     \item[] Projection from the leaves.
       \end{itemize}
              {\bf enddo}
      \item[] {\bf do} $l=\tilde{l}_k,\ldots,L$ 	\begin{itemize}
     \item[] Thresholding, prediction, and
        addition of the safety zone.
	  \end{itemize}
	  {\bf enddo}
   \end{enumerate}
   \item[]{\bf enddo}
\end{enumerate}
\end{Alg}

Here, $\tilde{l}_k$ denotes the coarsest level containing leaves in the intermediate step $k$ (as
introduced in \cite{MS}),
$h(l)$ is the mesh size on level $l$, and $|\tilde{\Lambda}|_z(l)$ is the size of the set formed
by leaves and virtual leaves per
resolution level $l$ in the direction $z$. The interior marching formula  is the modified
marching formula \eqref{bbrs_marchingf_model2}
for Model~2, for the intermediate time steps $k=1,\ldots,2^L$, for the leaf in the
position $(i,j)$ at level $l$.

Our scheme is formulated for a first-order,
 explicit Euler time discretization. Generalizations for
 higher order schemes are given in \cite{MS} for
 Crank-Nicolson schemes and in \cite{DRS062} for Runge-Kutta schemes,
 respectively.

\section{Algorithm implementation} \label{bbrs_sec:algimp}
Now we give a brief description of the multiresolution procedure used to solve
the test problems.
\begin{Alg}[Multiresolution procedure] \hfill
\begin{enumerate}
\item Initialization of parameters.
\item Creation of the initial graded tree structure:
  \begin{enumerate}
      \item Create the root of the tree and compute its cell average value.
      \item Split the cell, compute the cell average values in the sons and
compute the corresponding details.
      \item Apply the thresholding strategy for the splitting of the
	new sons: If $d_{i,j,l}>\epsilon_l$ then split the son
(here we use $\smash{d_l=\min \{d_l^u,d_l^v\}}$).
      \item Repeat this until all sons have details below the required
tolerance $\varepsilon_l$.
  \end{enumerate}
\item {\bf do}  $n=1, \dots,  total\_ time\_ steps$ \label{bbrs_doloop}
\begin{enumerate}
\item Determination of the leaves and virtual leaves sets.
\item Time evolution with global time step: \label{bbrs_xxx} Compute the
discretized divergence operator for all the leaves.
\item \label{bbrs_updating_tree} Updating the tree structure:
  \begin{itemize}
      \item Recalculate the values on the nodes and the virtual nodes by
projection from the leaves. Compute the details for all positions $(\cdot,\cdot,l)$
for $l\geqslant \tilde{l}_k$.

{\bf if} $|\bar{d}_{i,j,l}|<\eps_l$
(here we use $d_l=\max(d_l^u,d_l^v)$) in a node {\bf and} in its brothers {\bf then}
the cell and its brothers are \emph{deletable}.

{\bf endif}
      \item {\bf if} some node and all its sons are deletable {\bf and} the sons
are leaves without virtual sons, {\bf then} delete sons (coarsening).
\begin{itemize}
\item[] {\bf if} this node has no
sons {\bf and} it is not deletable {\bf and} it is not at level $l=L$, {\bf then} create sons (refinement).
\item[] {\bf endif}
\end{itemize}
{\bf endif}
      \item Update the values in the new sons by prediction operator from the
former leaves.
\end{itemize}
\end{enumerate}
\item[] {\bf enddo}
\item Output: Save meshes, leaves and cell averages.
\end{enumerate}
\end{Alg}
Here  $total\_ time\_ steps$ stands for the total time steps needed to reach
$T_{\text{final}}$ using $\Delta t$ as the maximum time step allowed by the
CFL condition using the finest space step.

When using a locally varying time stepping, replace step (\ref{bbrs_xxx}) by
the new step
\begin{enumerate}
\item[(3)]  {\bf do}  $n=1, \dots,  total\_ time\_ steps$
\begin{enumerate}
\item \emph{Determination of the leaves and virtual leaves sets}.
\item \emph{Time evolution with local time stepping:
 Compute the discretized divergence operator for all the leaves and
virtual leaves}
\item {\bf do} $k=1, \dots , 2^L$ \emph{($k$ counts intermediate time
  steps)}
\begin{itemize}
\item {\em Compute the
intermediate time steps depending on the position of the
corresponding leaf as explained in Section~\ref{bbrs_sec:LTS}.}
\item {\bf if} $k$ is odd {\bf then}  \emph{update the tree
  structure:}
\begin{itemize}
      \item \emph{Recalculate the values on the nodes and the virtual nodes by
projection from the leaves. Compute the details in the whole tree.}

{\bf if} \emph{$|\bar{d}_{i,j,l}|<\eps_l$
(here we use $d_l=\max(d_l^u,d_l^v)$) in a node {\bf and} in its brothers {\bf then}
the cell and its brothers are \emph{deletable}.}

{\bf endif}
      \item {\bf if} \emph{some node and all its sons are deletable {\bf and} the sons
are leaves without virtual sons, {\bf then} delete sons (coarsening).}
\begin{itemize}
\item[] {\bf if} \emph{this node has no
sons {\bf and} it is not deletable {\bf and} it is not at level $l=L$, {\bf then}
create sons (refinement).}
\item[] {\bf endif}
\end{itemize}
{\bf endif}
      \item \emph{Update the
values in the new sons by prediction operator from the
former leaves.}
\end{itemize}
\item[] {\bf endif}
\end{itemize}
\item[] {\bf enddo}
\item[] \emph{(Now, after $2^L$ intermediate
  steps,  all nodes are synchronized.)}
\end{enumerate}
\item[] {\bf enddo}
\end{enumerate}
Here  $total\_ time\_ steps$ for the new step stands for the total time steps needed to reach
$T_{\text{final}}$, with $\Delta t_0$ as the maximum time step allowed by the
CFL condition using the coarsest space step.

Notice that with such a process,  we obtain high-order approximation in
the smooth regions and mesh refinement near discontinuities as a consequence
of the polynomial accuracy in the multiresolution prediction operator,
even if the reference finite volume scheme is low-order accurate.

Bihari and Harten \cite{bh97} use the following quantity,
 which we call
 \emph{data compression rate},
\begin{equation}\label{bbrs_mu}
\eta:=\frac{N_{Lx}N_{Ly}}{N_{Lx}N_{Ly}/2^{2L}+|\mathcal{L}(\Lambda)|},
\end{equation}
to  measure  the possible improvement in  data compression,
 whose feasibility  in turn  strongly  depends  on
  a smart implementation to navigate inside the tree,
see for example   \cite{RSTB03}. Here,  $N_{Lx}N_{Ly}$
is the number of cells in the finest grid, and $|\mathcal{L}(\Lambda)|$ is the
size of the set of leaves. We measure the \emph{speed-up} between the CPU time of
the numerical solution obtained by the FV method and the CPU time of the numerical
solution obtained by the multiresolution method:
$V=(\mathrm{CPU\, time})_{\mathrm{FV}}/(\mathrm{CPU\, time})_{\mathrm{MR}}$.

To measure errors between a reference solution $u$ and an approximate solution $u_{MR}$
obtained using multiresolution, we will use $L^p$-errors:
$e_p=\smash{ \|u^n-u^n_{\mathrm{MR}}\|_p}$, $p=1,2,\infty$, where
\begin{gather*}
e_\infty=\max_{1\leq i \leq N_{Lx},1\leq j\leq N_{Ly}}\bigl|u^n_{i,j}-
{u^n_{\mathrm{MR}}}_{i,j,L} \bigr|, \\
e_p=\left(\frac{1}{N_{Lx}N_{Ly}}\sum_{i=1}^{N_{Lx}}\sum_{j=1}^{N_{Lx}}
\bigl|u^n_{i,j}-{u^n_{\mathrm{MR}}}_{i,j,L} \bigr|^p\right)^{1/p},\quad p=1,2.
\end{gather*}
Here $\smash{{u^n_{\mathrm{MR}}}_{i,j,L}}$ is the value on the finest level $L$ obtained
by prediction from the corresponding leaf.

\section{Numerical results}\label{bbrs_sec:results}

\begin{figure}[t]
\begin{center}
\begin{tabular}{cc}
\includegraphics[width=0.44\textwidth,height=0.352\textwidth]{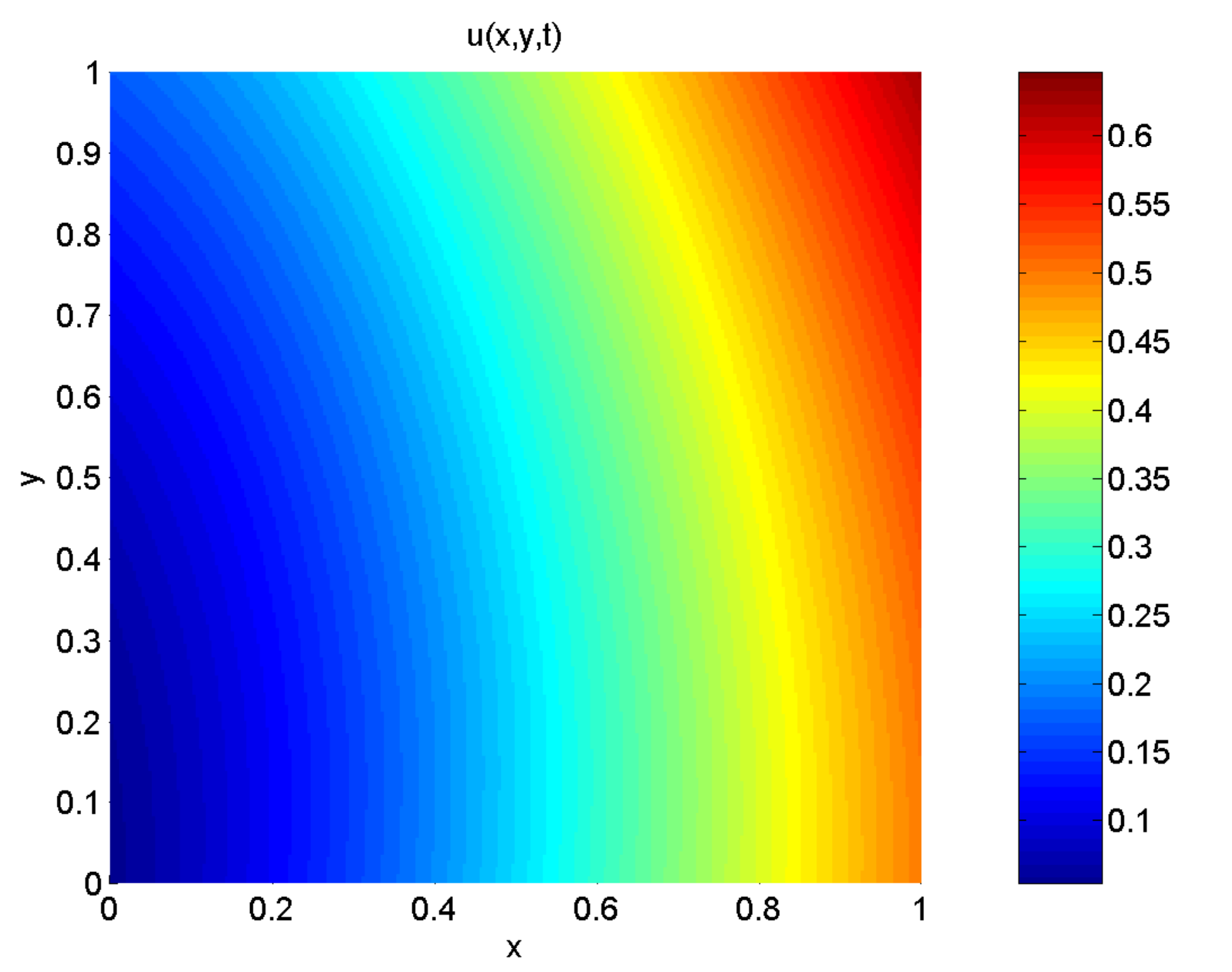}&
\includegraphics[width=0.352\textwidth,height=0.352\textwidth]{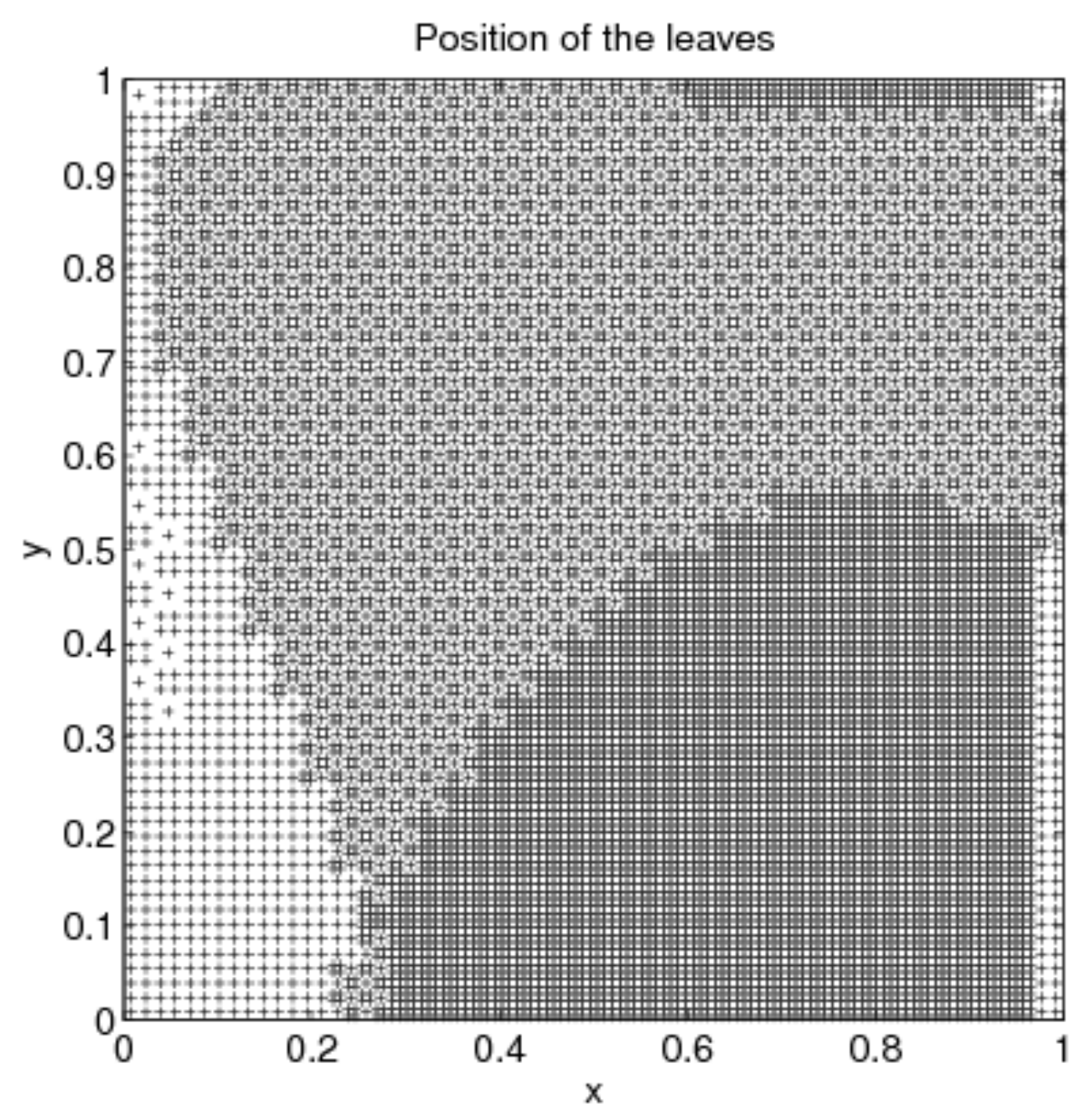}
\\
\includegraphics[width=0.44\textwidth,height=0.352\textwidth]{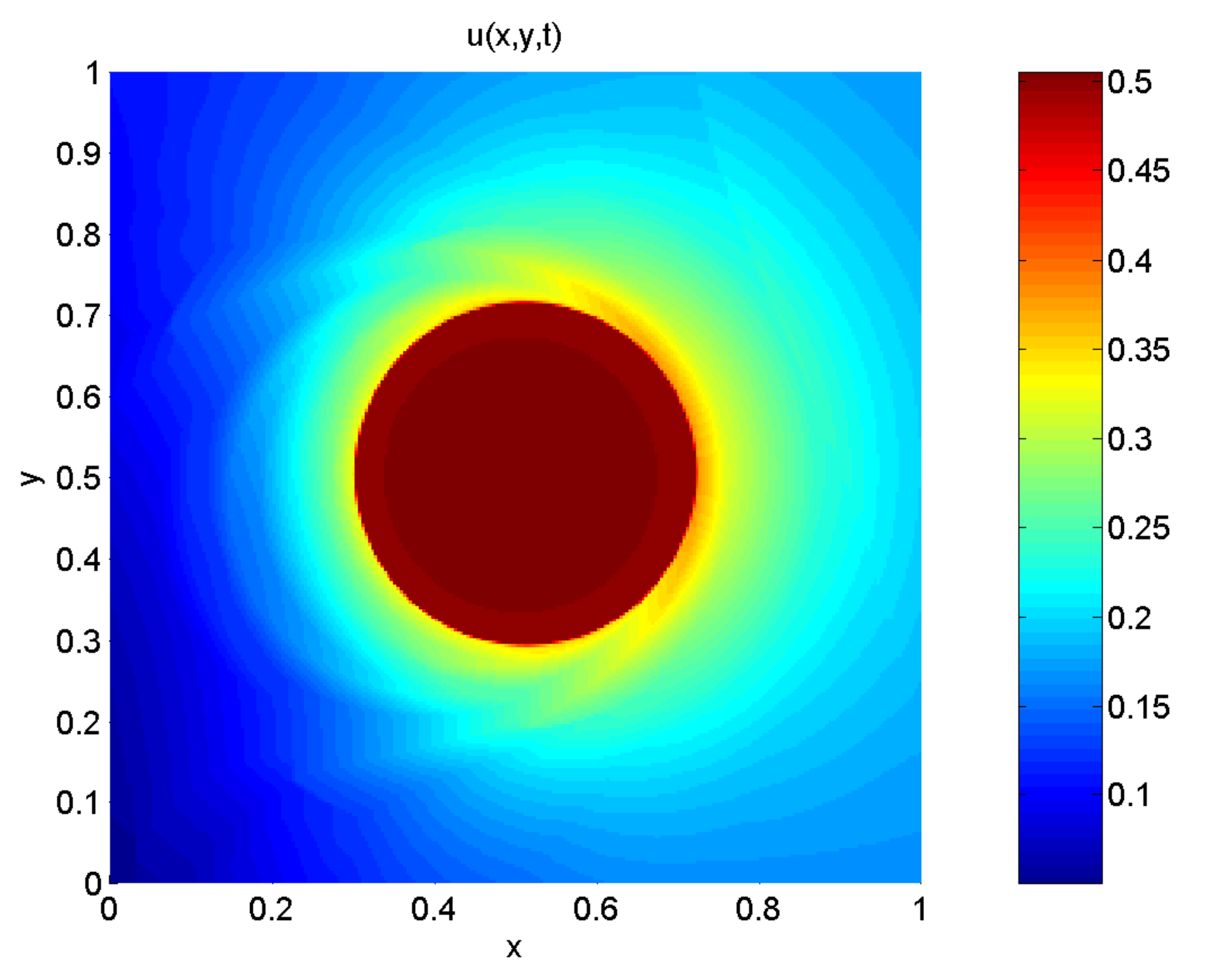}&
\includegraphics[width=0.352\textwidth,height=0.352\textwidth]{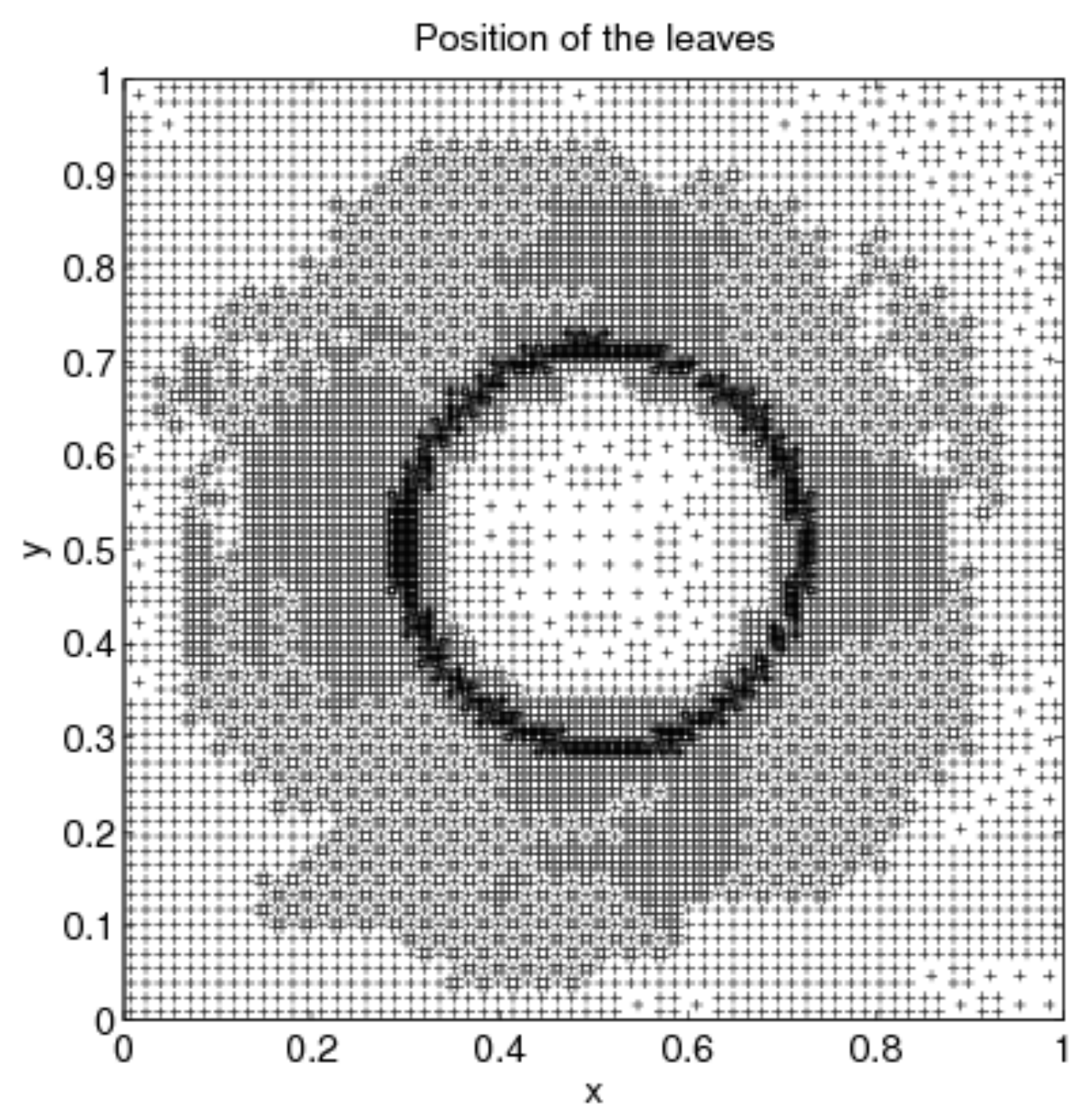}
\\
\includegraphics[width=0.44\textwidth,height=0.352\textwidth]{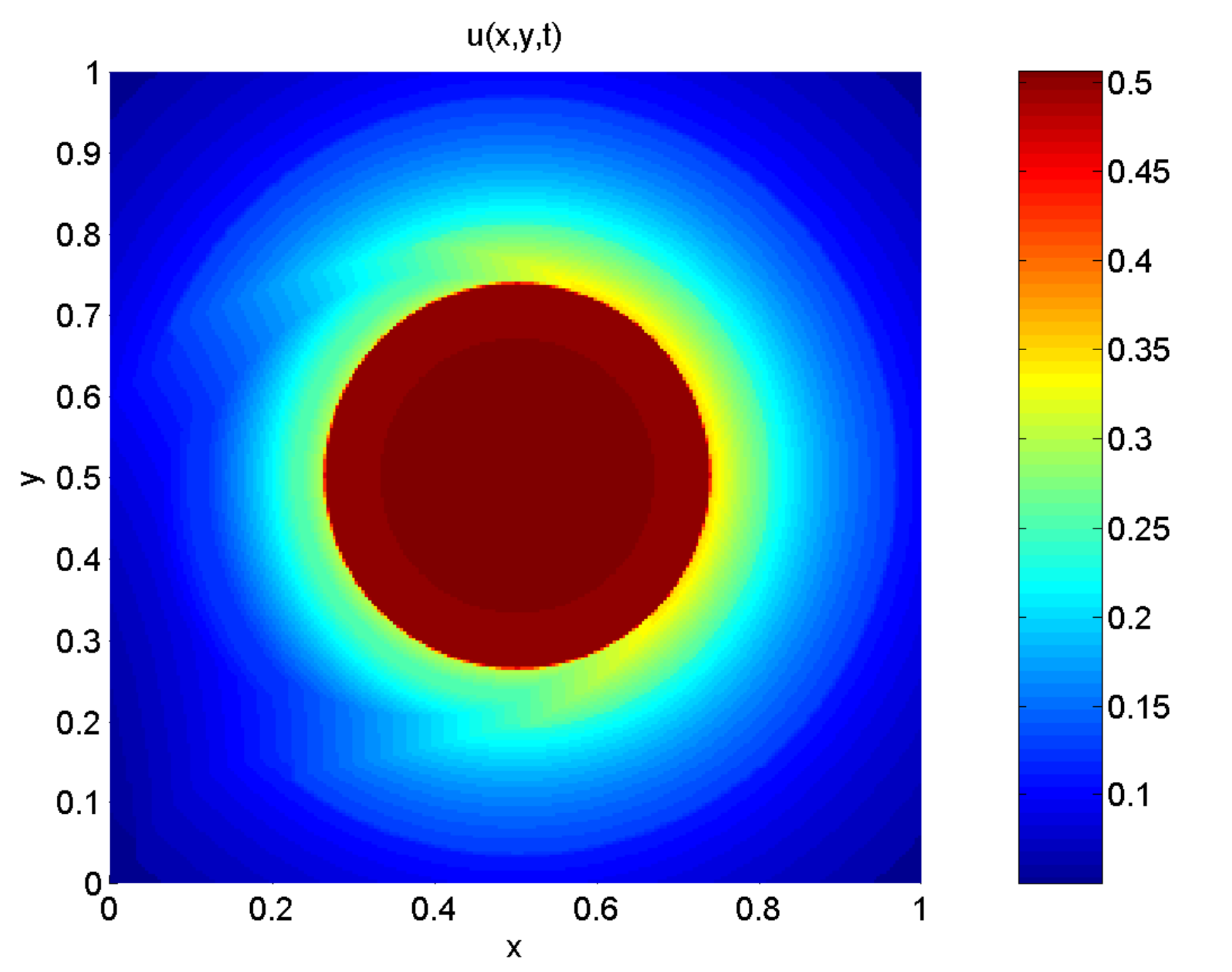}&
\includegraphics[width=0.352\textwidth,height=0.352\textwidth]{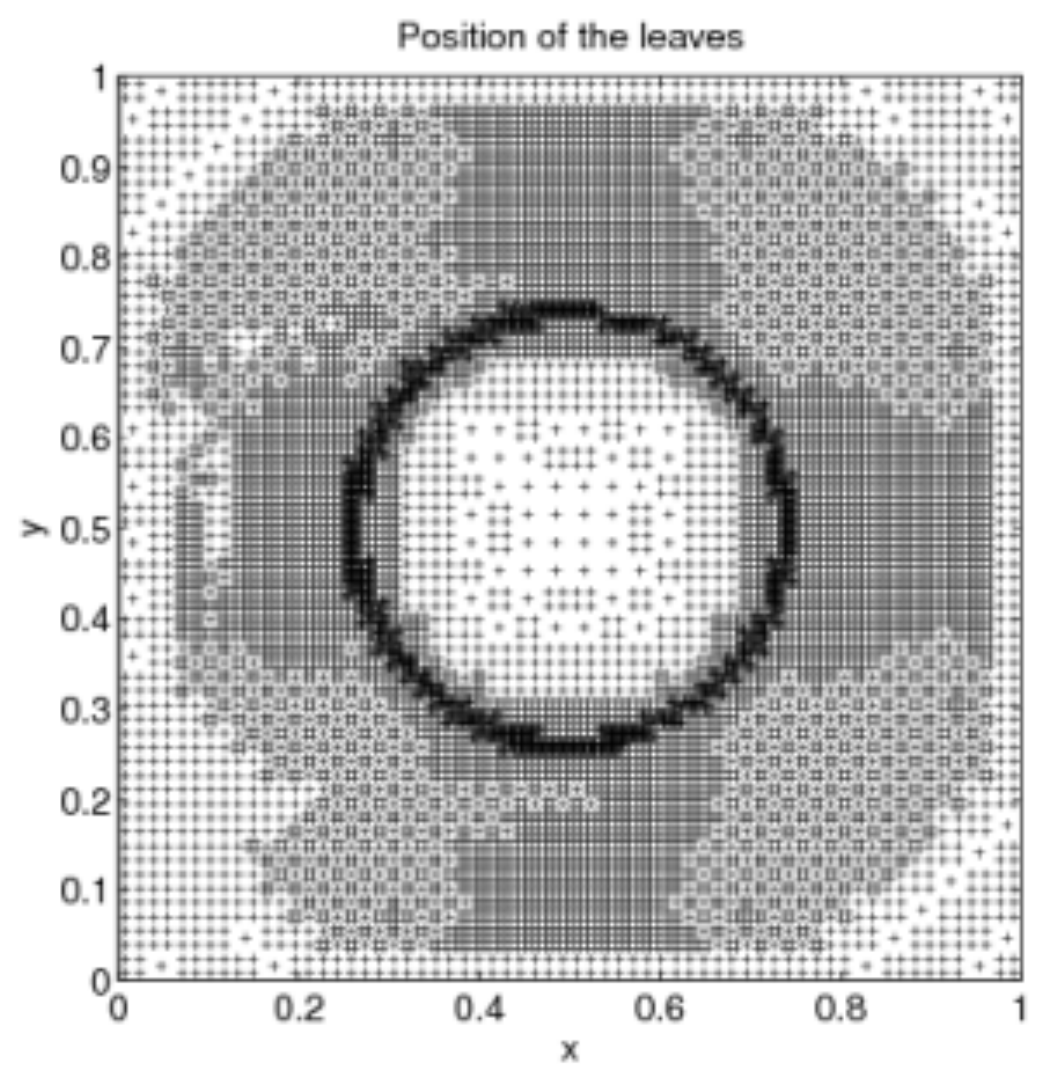}
\end{tabular}
\caption{\it Example~1 (single-species model): Numerical solution
(left) and leaves (right)
 of the corresponding tree data structure at times $t=0$ (top),
   $t=0.5$ (middle)  and $t=3$ (bottom).}
\label{bbrs_fig:one-a}
\end{center}
\end{figure}

 \begin{figure}[t]
 \begin{center}
 \begin{tabular}{ccc}
 (a) & (b) & (c)\\ 
\includegraphics[width=0.305\textwidth]{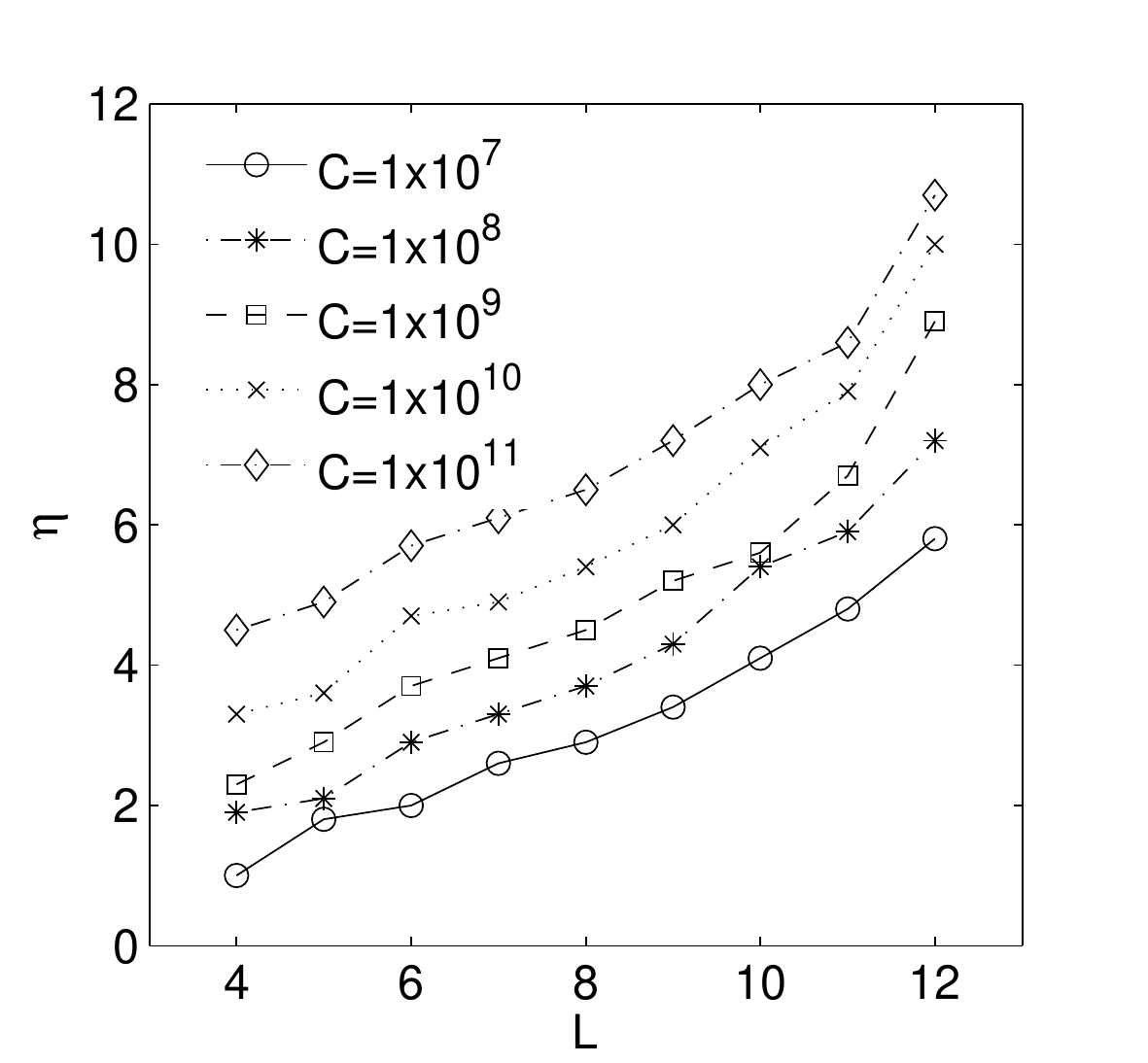}&
\includegraphics[width=0.305\textwidth]{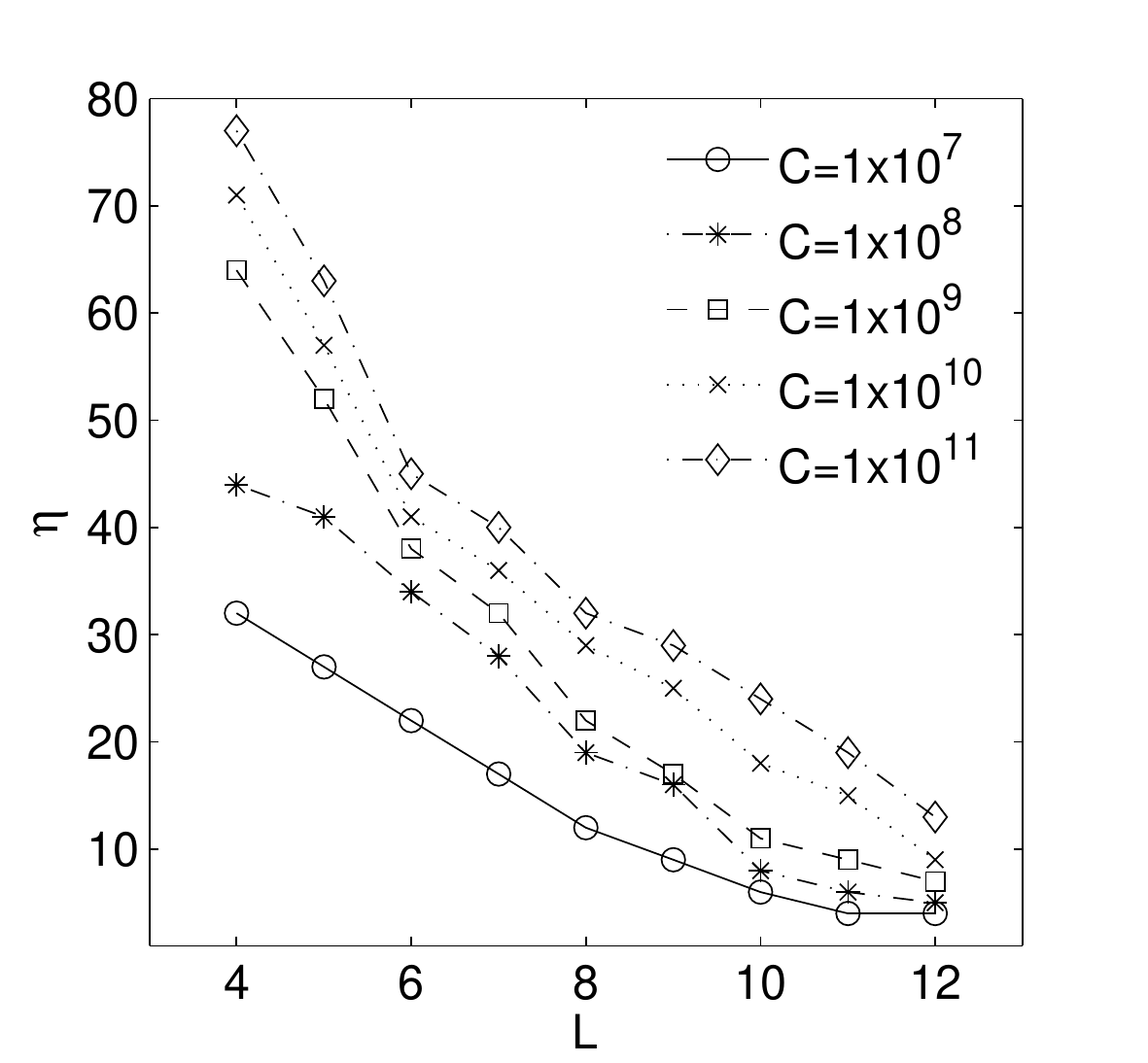}&
\includegraphics[width=0.305\textwidth,height=0.28\textwidth]{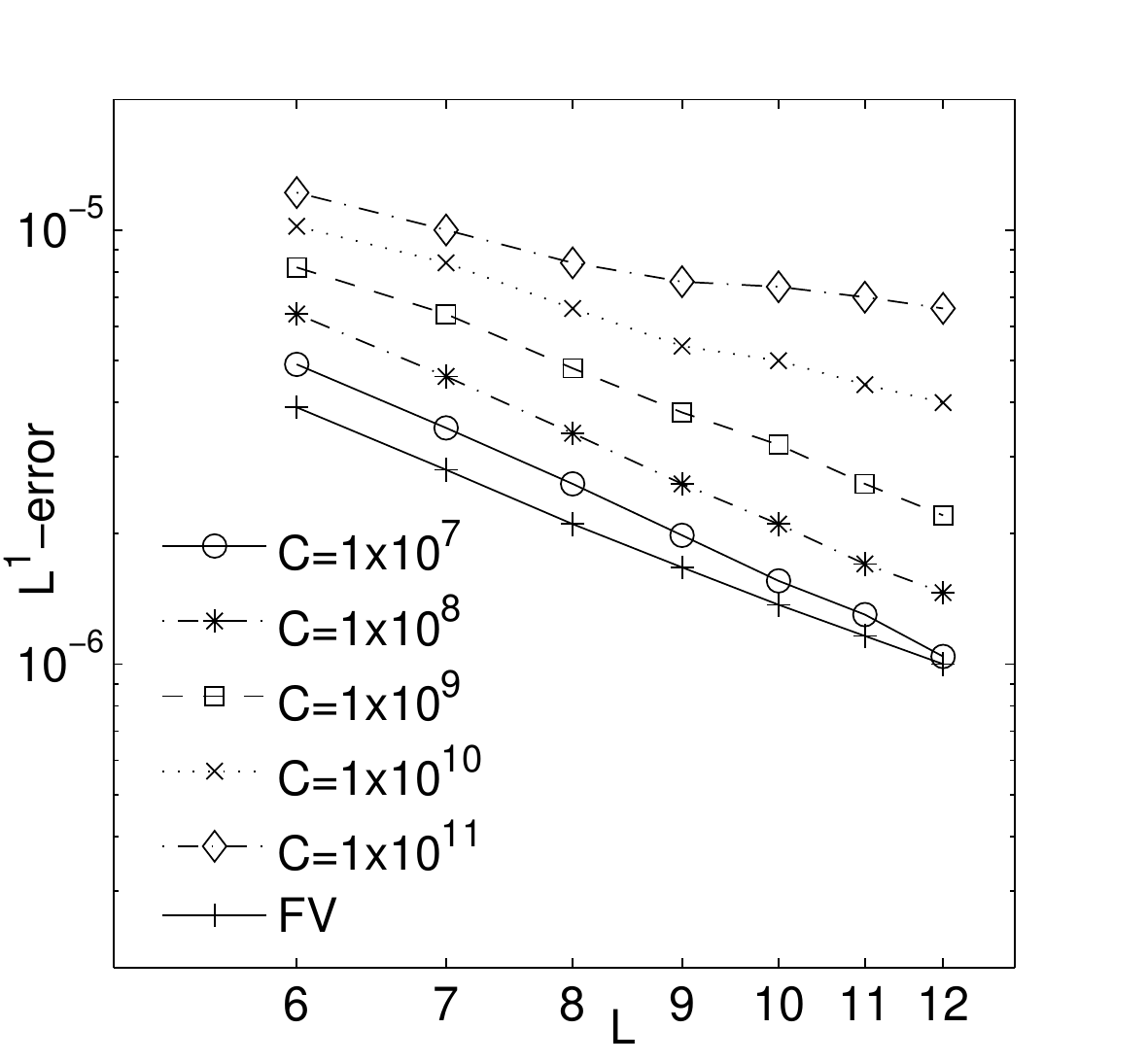}
 \end{tabular}
 \caption{\it Example~1 (single-species model): (a) data compression
rate~$\eta$, (b) speed-up factor~$V$ and (c) $L^1$-errors for different
levels~$L$ and values of~$C$.
The simulated time is $t=0$.}
 \label{bbrs_fig:ex1_factorC}
 \end{center}
 \end{figure}

\subsection{Example~1: Single-species model}
For this example,  we consider \eqref{bbrs_eq:reac-dif} with a strongly degenerate diffusion
term \eqref{bbrs_model1audef}, where we choose $D:=1$ and
$u_{\mathrm{c}}:=0.5$,   a square domain $\Omega=[0,1]^2$, and
 the function  $f(u,\x)$ given by \eqref{bbrs_fxueq}.
Figure~\ref{bbrs_fig:one-a} displays
   the numerical solution starting from the
 smoothly varying   initial function
\begin{align*}
u_0(x,y)=0.5\big(1+\sin(1.1(x-\cos(0.7y))\big)
\cos\bigl( 0.5(y-\sin(1.3x))\bigr).
\end{align*}
We choose a maximal resolution level of $N_L=256^2=65536$ control volumes
on the finest grid.  Figure~\ref{bbrs_fig:ex1_factorC}
 illustrates how the factor $C$ in \eqref{bbrs_equ:epsref1} is selected
 in our case as the optimal value from a finite selection  of test values (each
 value giving a different value for the reference tolerance \eqref{bbrs_equ:epsref1}
 $\varepsilon_{\mathrm{R}}$).
Figures~\ref{bbrs_fig:ex1_factorC} (a) and~(b) indicate  that for all displayed levels,
the multiresolution procedure is in every case (for different
values of $C$) cheaper (in terms of both acceleration and memory savings)
 than the corresponding reference FV computation  on the finest
grid. Figure~\ref{bbrs_fig:ex1_factorC}~(c) indicates that the
computations obtained using $C=1.0\times10^{9}$ (and hence
$\varepsilon_{\mathrm{R}} =9.43\times10^{-4}$) are sufficiently accurate, in the sense
that with these choices,  we keep the same slope for the $L^1$-error as the FV calculations while
increasing~$V$ and~$\eta$. We remark that here, as in  previous
  works  that use  similar methods
(see  e.g.  \cite{CKMP2002}), there actually exists a range of threshold parameters that preserve
the same slope for errors with respect to reference solution, for  which $C=1.0\times10^{9}$
is an average value. Here,  we compute errors using a reference FV solution on a
fine grid with $N_L=2048^2=4194304$ control volumes. (Here and in
 all other examples, we calculate errors in the approximate sense
 with respect to a reference solution.)

\subsection{Examples~2 and~3: Interaction between two flame balls}
\begin{figure}[t]
\begin{center}
\begin{tabular}{cc}
\includegraphics[width=0.44\textwidth,height=0.352\textwidth]{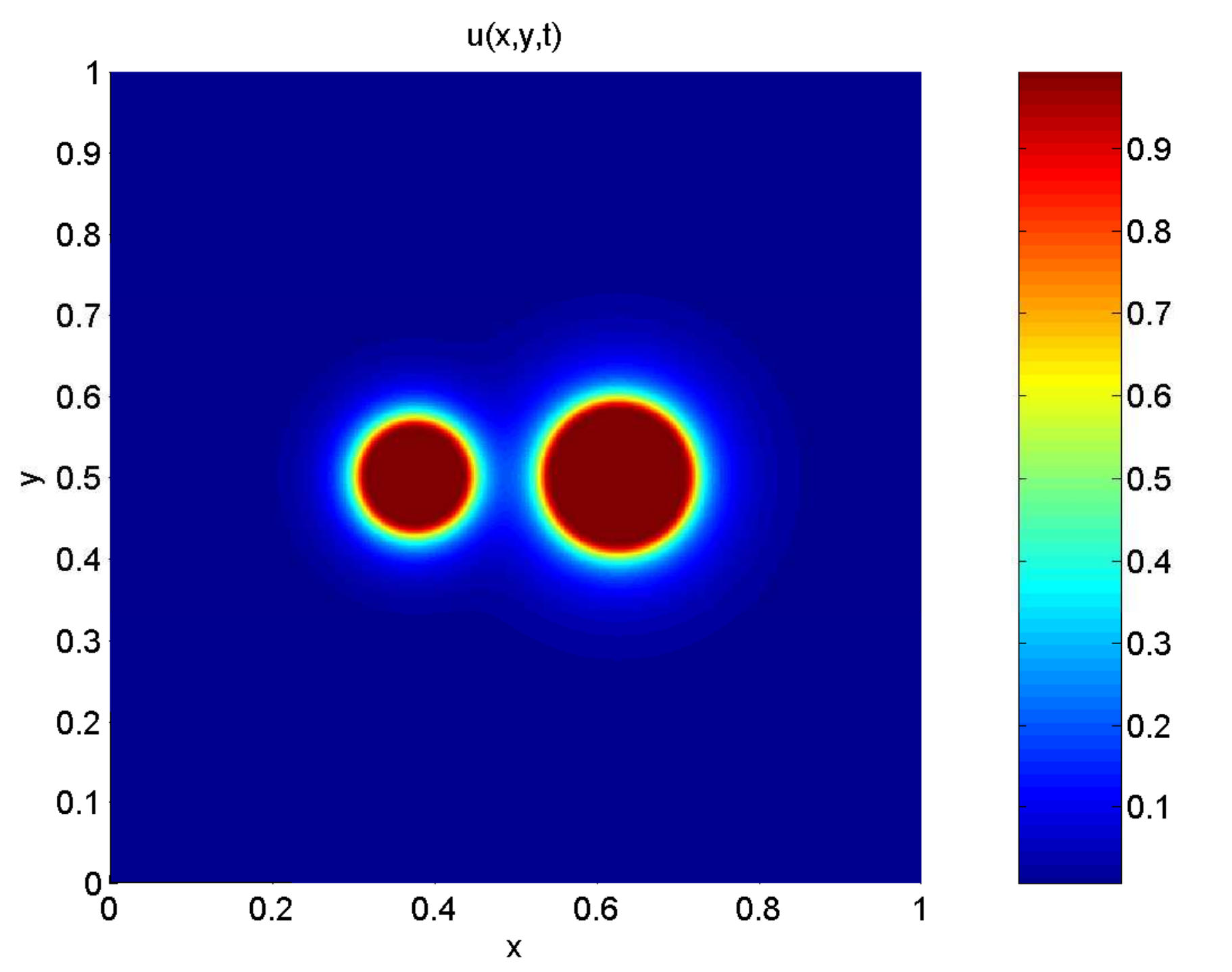}&
\includegraphics[width=0.352\textwidth,height=0.352\textwidth]{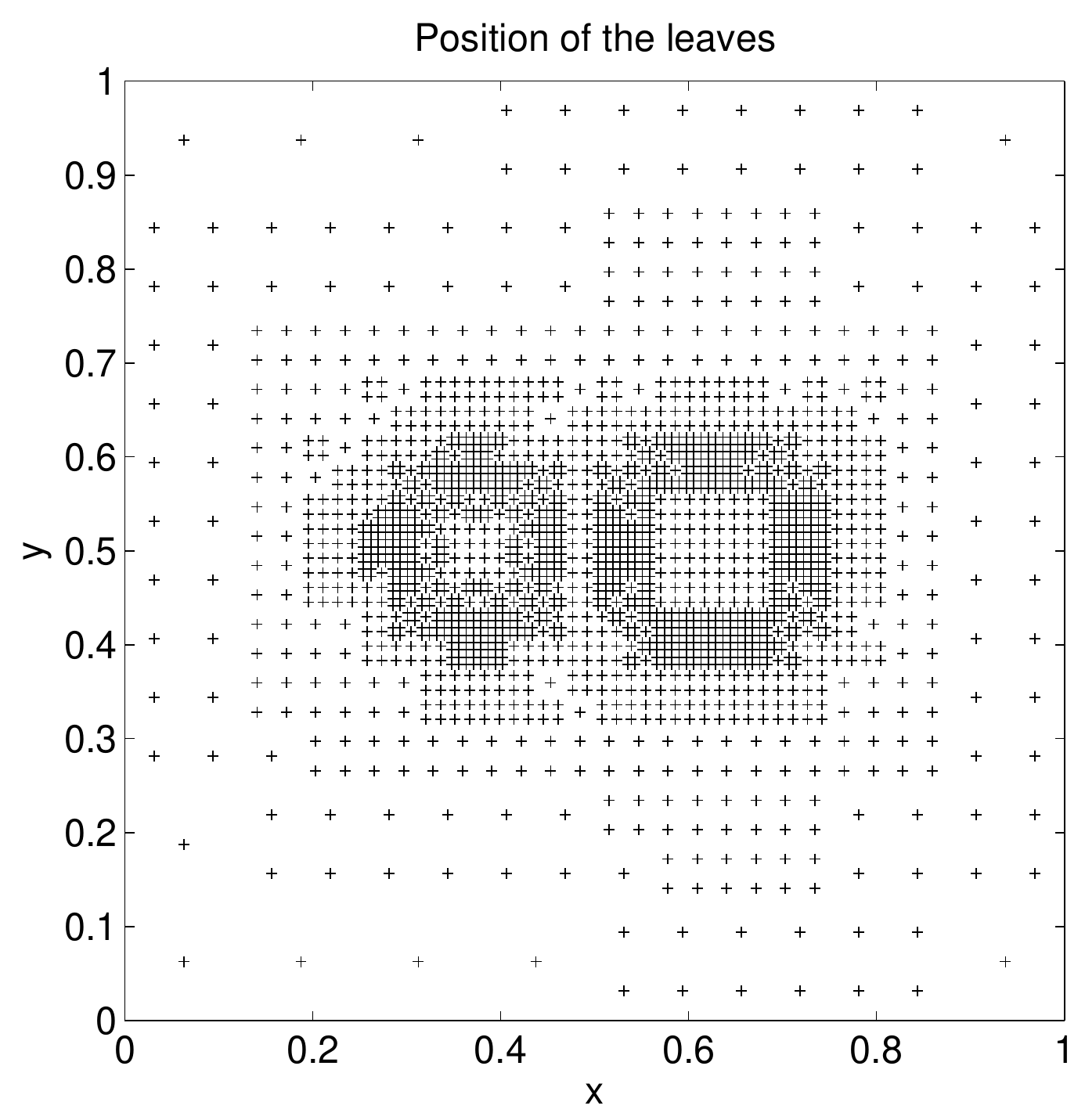}\\
\includegraphics[width=0.44\textwidth,height=0.352\textwidth]{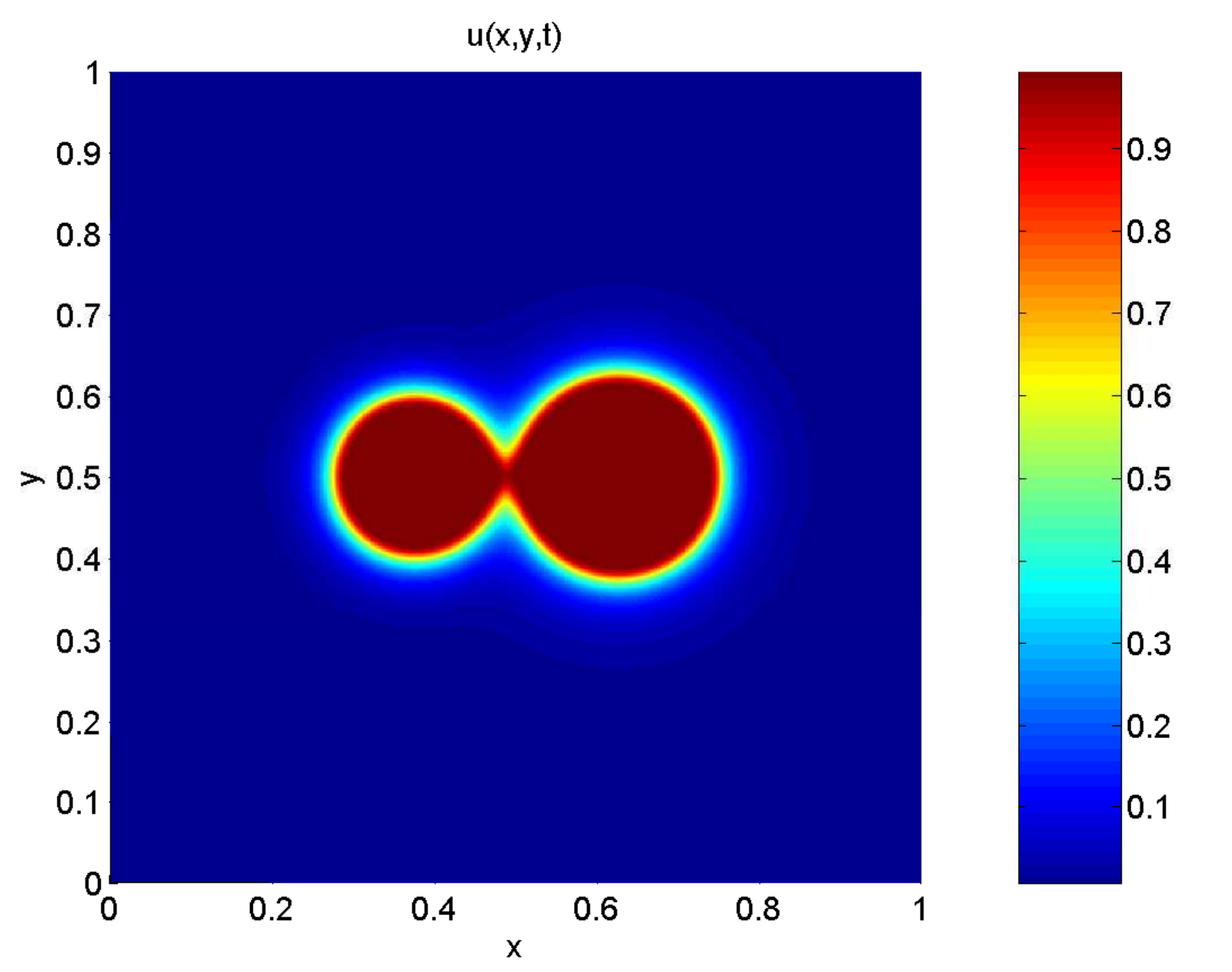}&
\includegraphics[width=0.352\textwidth,height=0.352\textwidth]{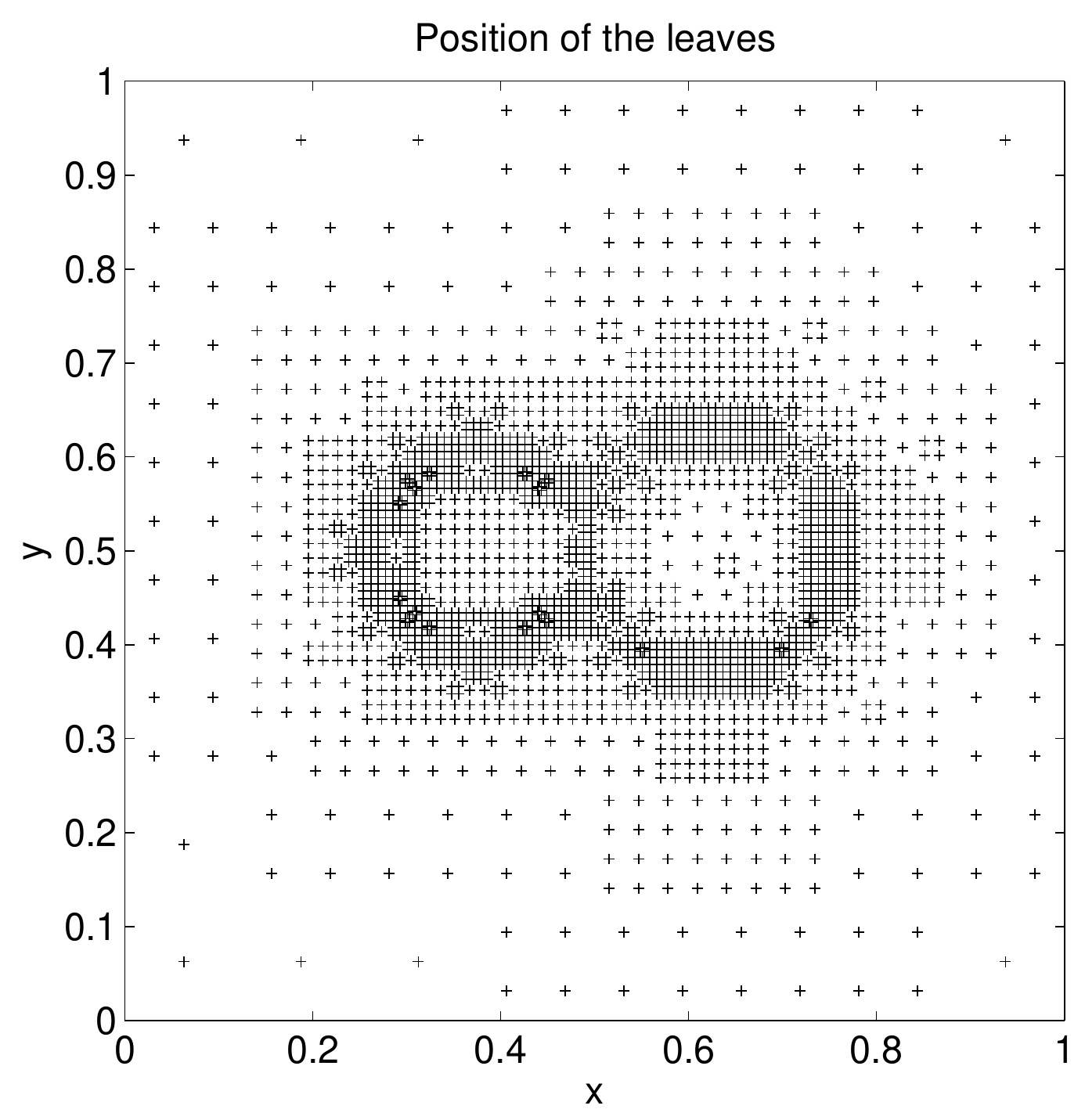}\\
\includegraphics[width=0.44\textwidth,height=0.352\textwidth]{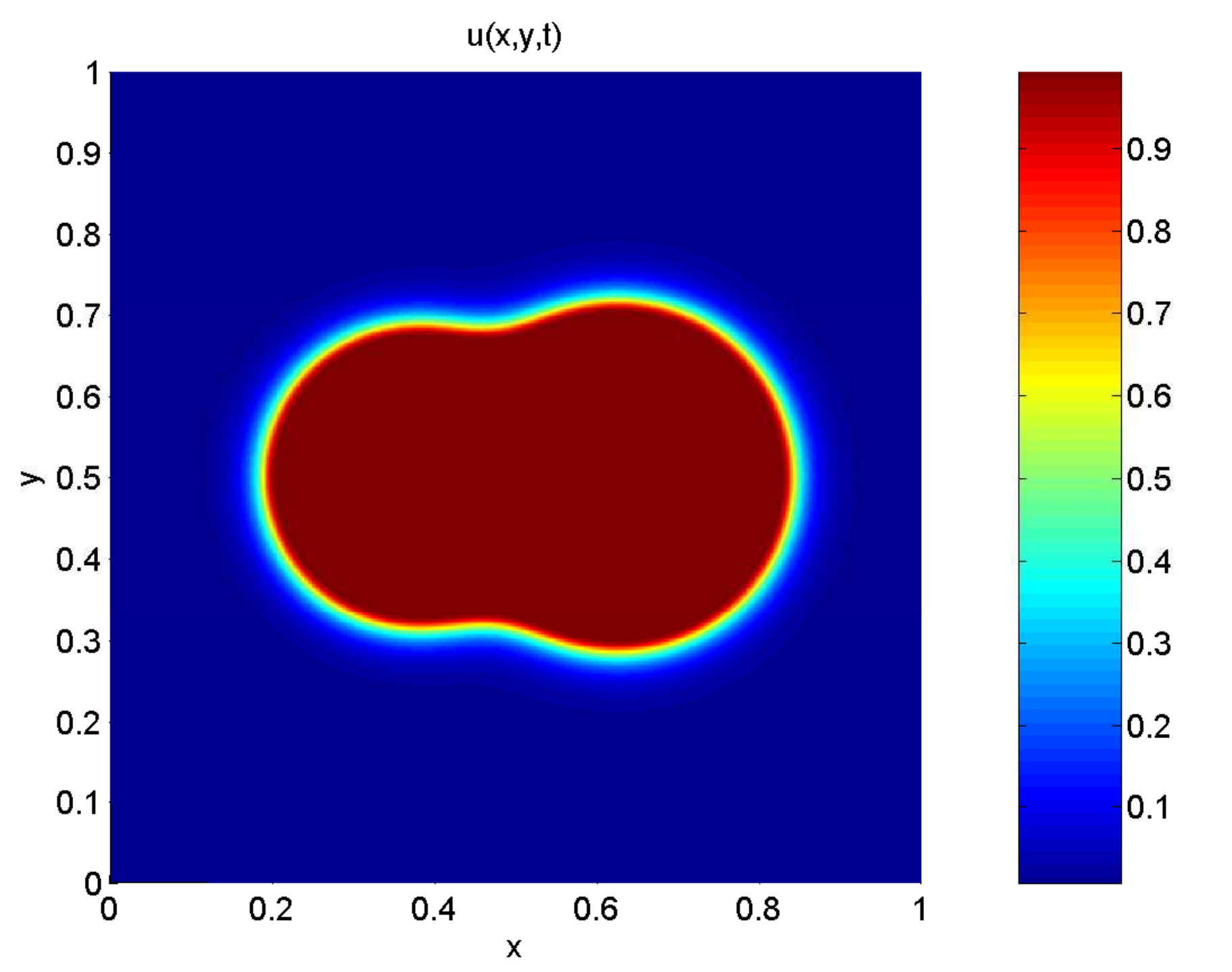}&
\includegraphics[width=0.352\textwidth,height=0.352\textwidth]{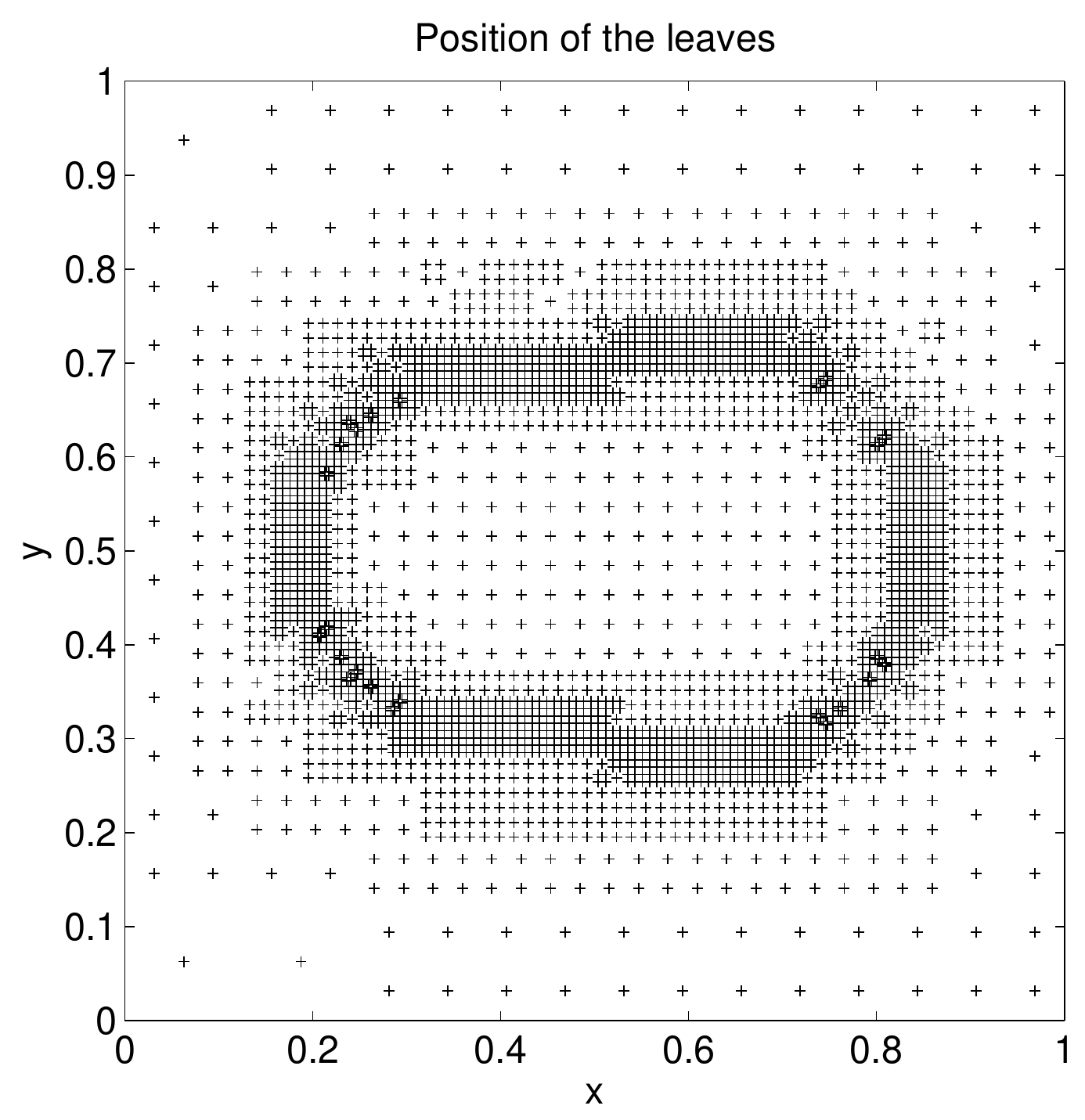}
\end{tabular}
\caption{\it Example~2 (interaction of two flame balls without radiation):
   Numerical solution for species $u$ (left) and  leaves of the corresponding
tree data structure (right) at times $t=2$ (top), $t=4$ (middle)  and
  $t=10$ (bottom).} \label{bbrs_fig:two_flame}
\end{center}
\end{figure}

 We study \eqref{bbrs_turing-sys} as a dimensionless model for the interaction between
two flame balls of different sizes. We consider a square domain
$\Omega=[-30,30]^2$ and  that the walls are sufficiently far from the
flame balls so that their influence is negligible. Physical
parameters characterizing the gaseous mixture and the combustion
 process are chosen as in \cite{RS05,ctm_rs04}.
We use the  parameters
$\alpha=0.64$ and $\beta=10$. The initial data is given by
$u(x,y,0)=u_0(r_1,r_2),\ v(x,y,0)=v_0(r_1,r_2)$ with  $r_1^2=(x-x_1)^2+y^2,\
r_2^2=(x-x_2)^2+y^2$, where
\begin{align} \label{bbrs_u0}
u_0(r_1,r_2)&:=\begin{cases}
1 &\text{if  $r_1< a$ or $r_2<b$,} \\
\max\bigl\{\exp(1-r_1/a),\exp(1-r_2/b)\bigr\}  & \text{otherwise,}
\end{cases} \quad
 v_0 (r_1, r_2) := 1- u_0(r_1,r_2).
\end{align}
In Example~2, we   simulate the process without  radiation, i.e.,
$\rho=0$ and hence $S(u)=0$. We set the Lewis number to $\mathrm{Le}=1$. Here $x_1=-7.5$,
$x_2=7.5$ and $a=1.8$, $b=2.5$ are the respective $x-$position and initial radii
of the two flame balls. This choice ensures that there is no interaction between
the two flame balls at $t=0$ and that there is no extinction of the flame balls.
We simulate the process until $t=10$, and Figure~\ref{bbrs_fig:two_flame} shows from
left to right the temperature and reaction rate configuration obtained using the
fully adaptive multiresolution scheme, and the position of the dynamic graded tree
leaves, which form  the corresponding adaptive  mesh. The different
times correspond from the top to the bottom to: before ($t=2$), during ($t=4$)
and after ($t=10$) direct interaction between the two flame balls, when the balls
tend to create a new circular flame structure. We choose the following
 multiresolution parameters: the maximal resolution level
$L=10$ corresponding to $N_L=512^2=262144$ control volumes in the finest grid,
 and the  reference tolerance $\varepsilon_{\mathrm{R}}=4.94\times10^{-3}$.

 For comparison purposes, we introduce the global chemical reaction rate
\begin{align*}
 R(t):=\iint_\Omega f(u,v)\, dx\,dy.
\end{align*}
Errors in different norms, reaction rates, information on  data compression and
speed-up rate for different methods at three different times are depicted in
Table \ref{bbrs_table:example_2}. Due  to the particular shape of
solutions, which is nearly constant away from the combustion front,
by using multiresolution, one can obtain very high rates of data compression, speed-up
and low errors.

\begin{figure}[t]
\begin{center}
\begin{tabular}{cc}
\includegraphics[width=0.48\textwidth,height=0.384\textwidth]{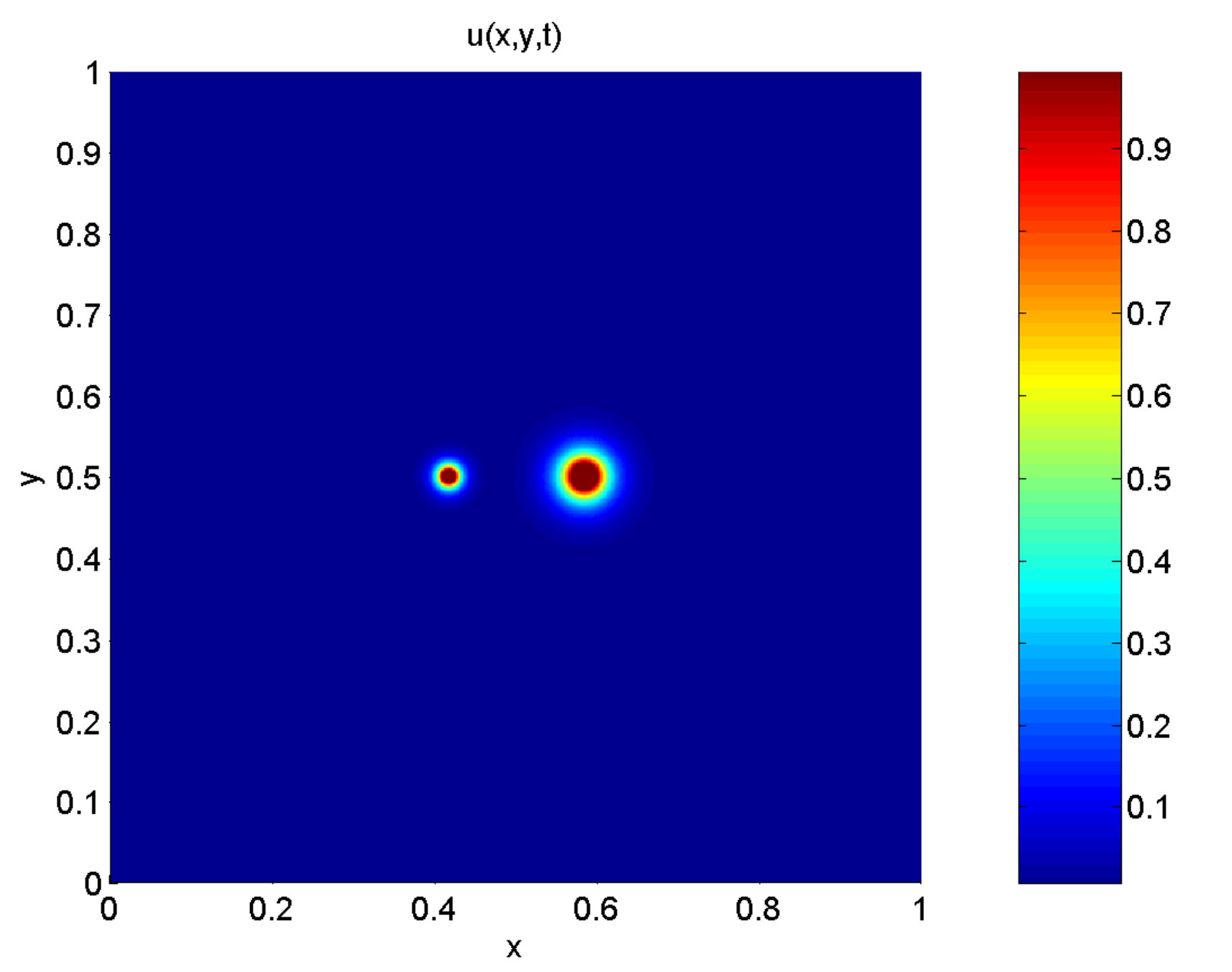}&
\includegraphics[width=0.384\textwidth,height=0.384\textwidth]{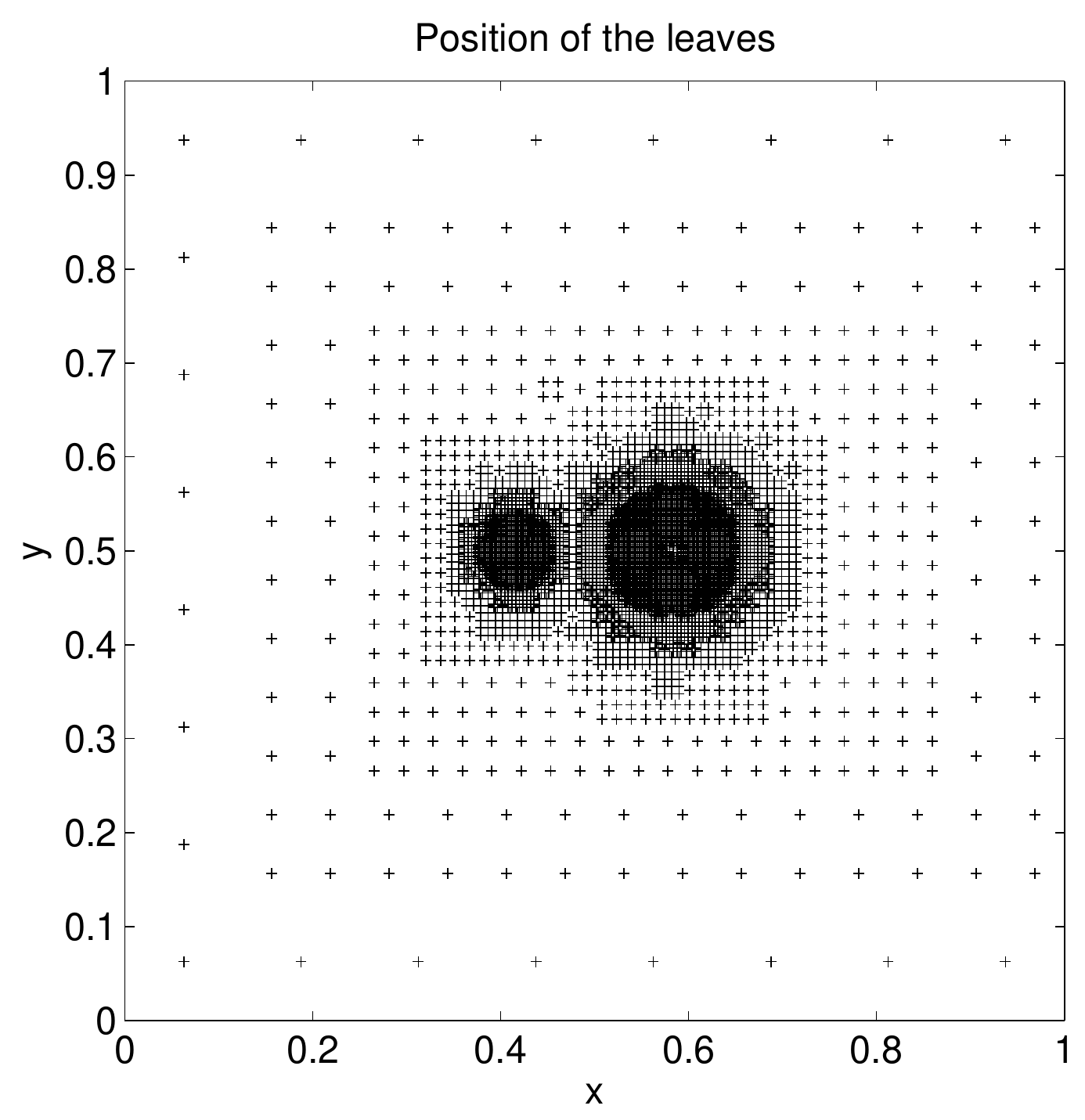}\\
\includegraphics[width=0.48\textwidth,height=0.384\textwidth]{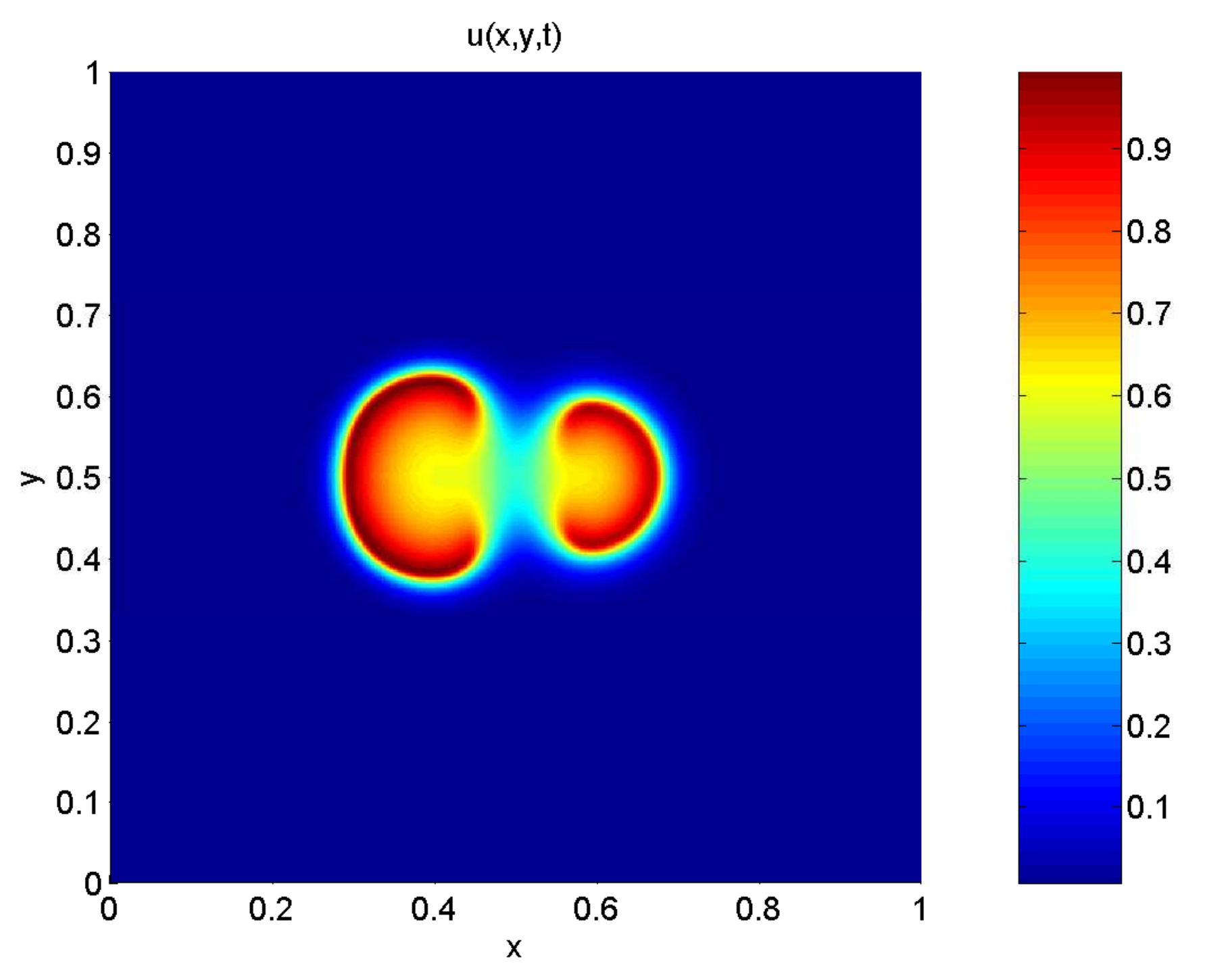}&
\includegraphics[width=0.384\textwidth,height=0.384\textwidth]{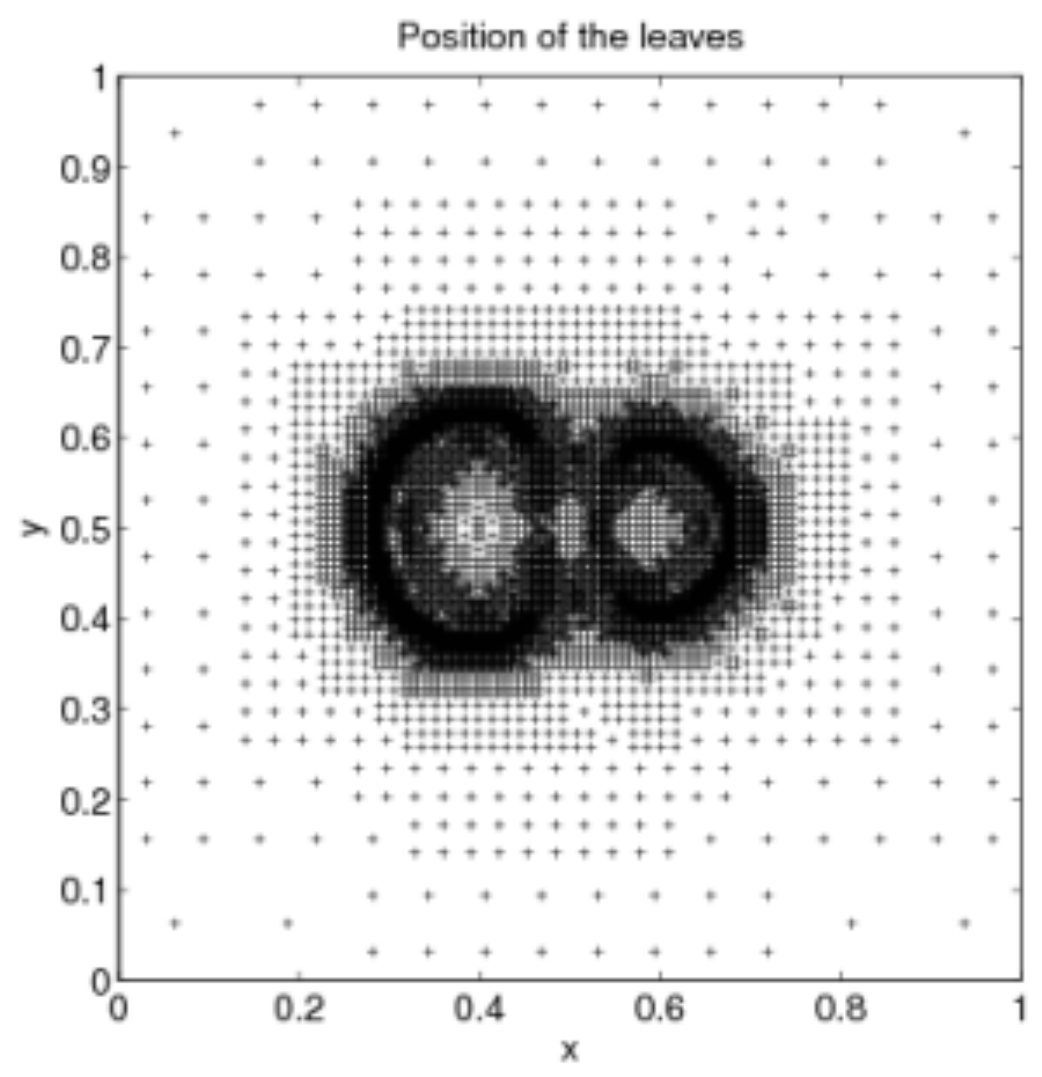}
\end{tabular}
\caption{\it Example~3  (interaction of two flame balls with radiation):
   Numerical solution for species $u$ (left) and  and leaves of the corresponding
tree (right) at times $t=0$ (top) and $t=5$ (bottom).} \label{bbrs_fig:two_flame_b}
\end{center}
\end{figure}

\begin{figure}[t]
\begin{center}
\begin{tabular}{cc}
\includegraphics[width=0.48\textwidth,height=0.384\textwidth]{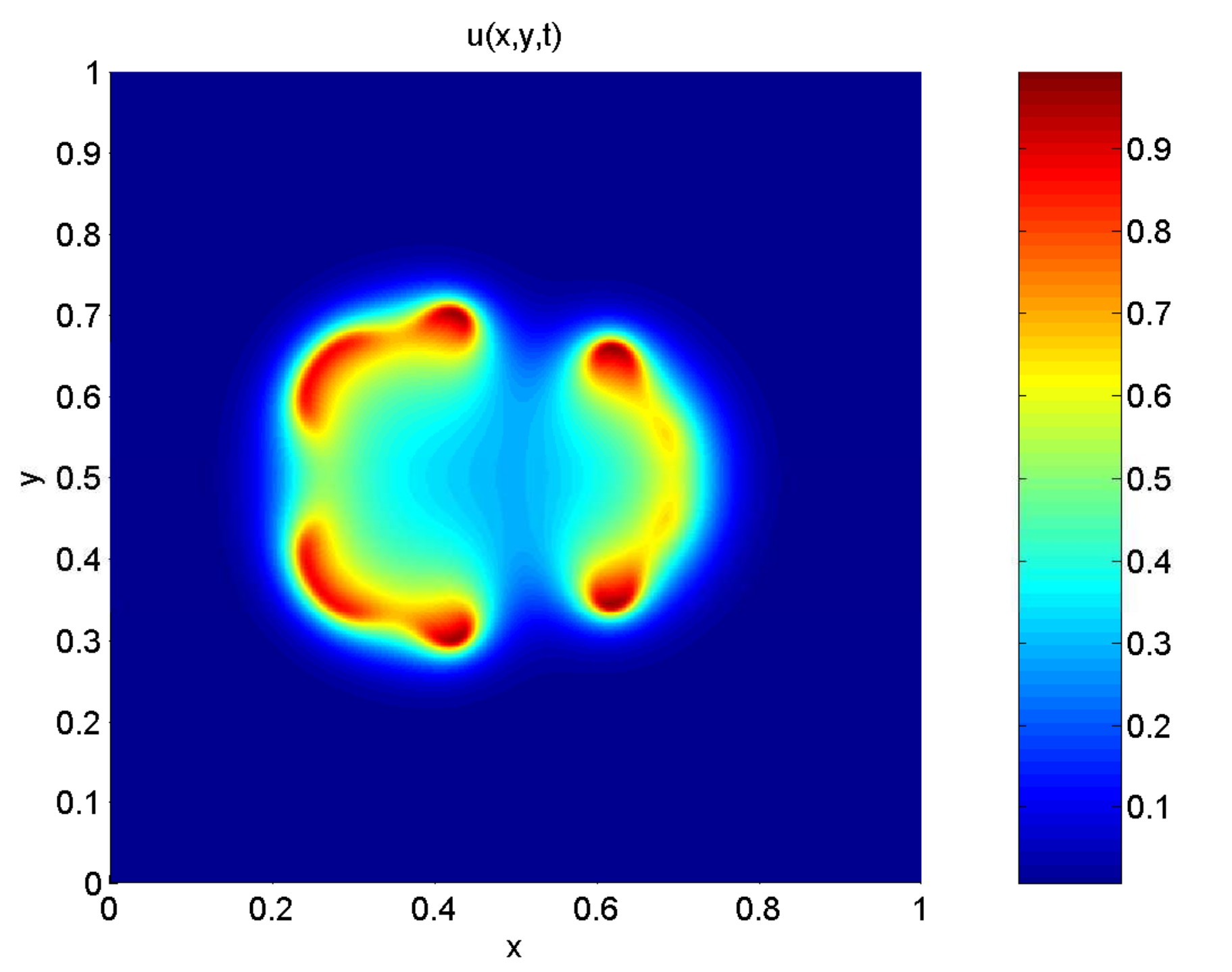}&
\includegraphics[width=0.384\textwidth,height=0.384\textwidth]{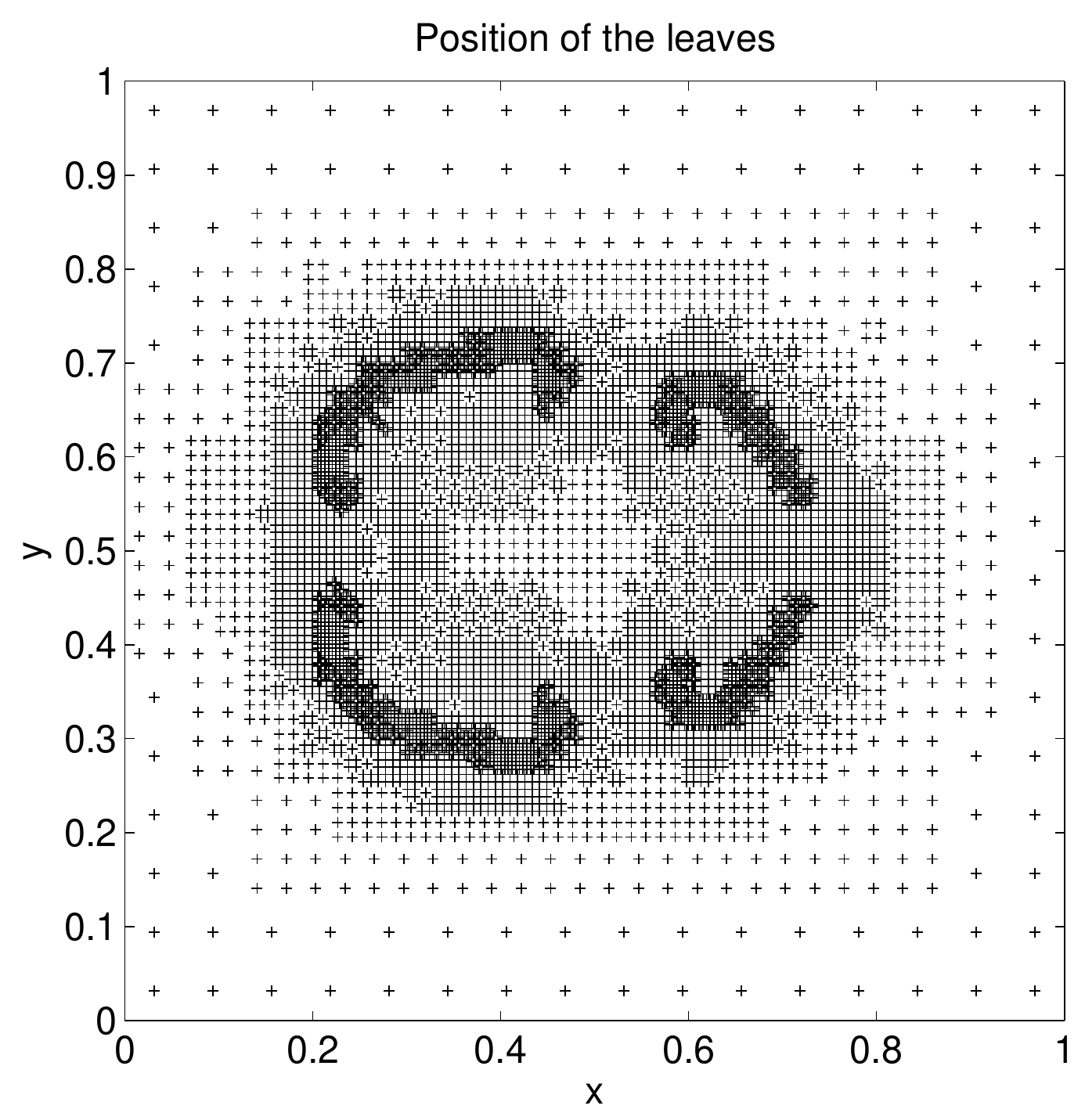}\\
\includegraphics[width=0.48\textwidth,height=0.384\textwidth]{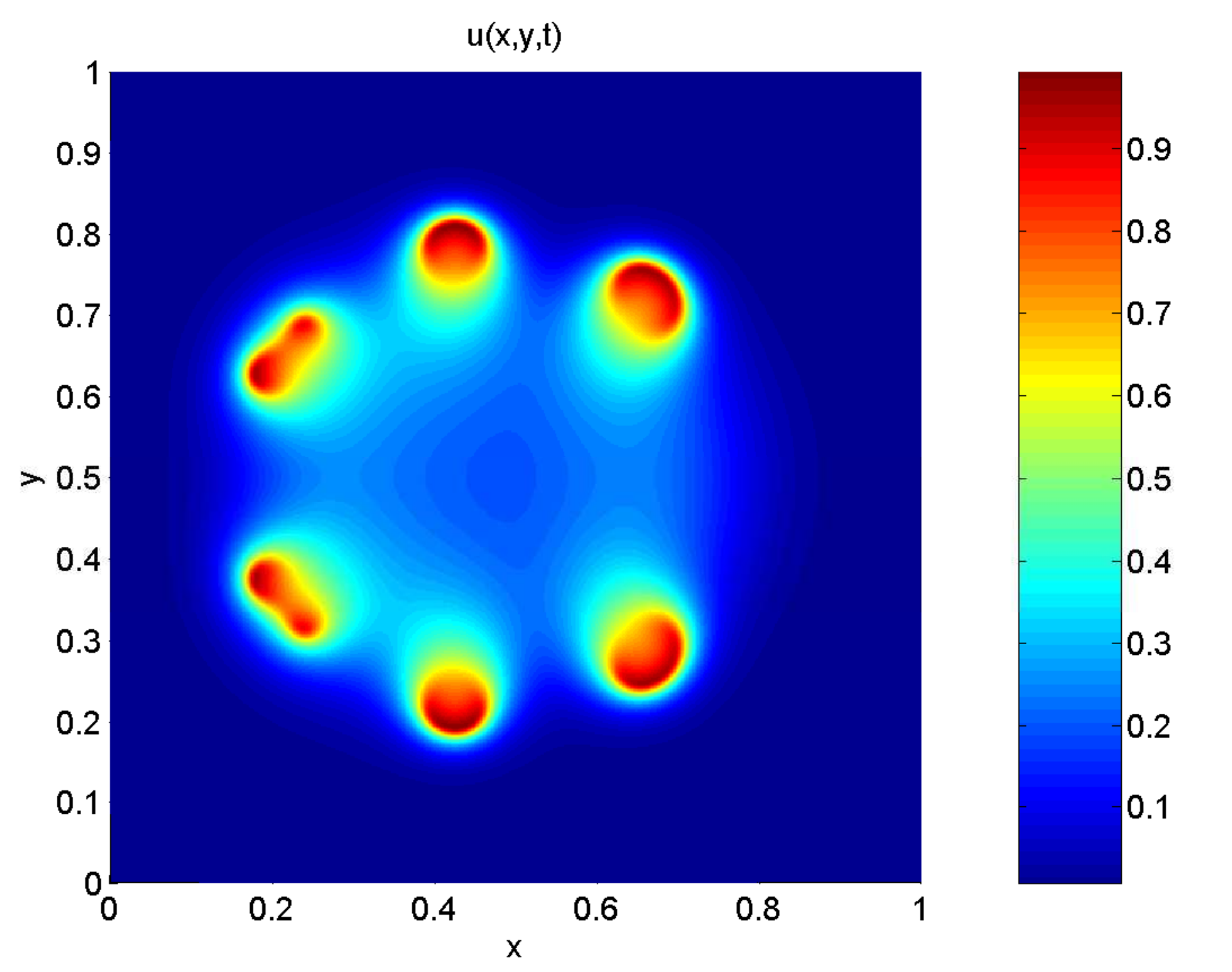}&
\includegraphics[width=0.384\textwidth,height=0.384\textwidth]{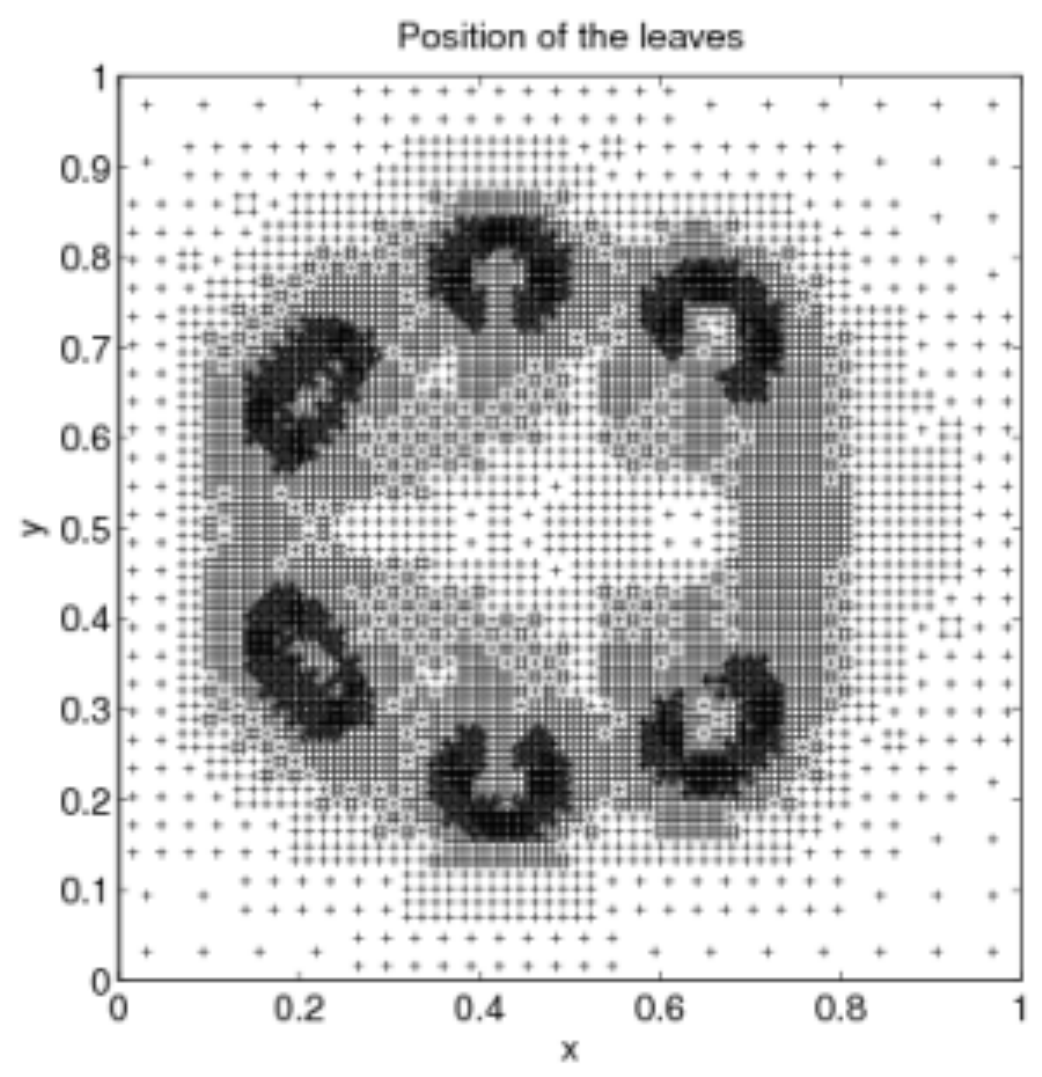}
\end{tabular}
\caption{\it Example~3 (interaction of two flame balls with radiation):
  Numerical solution for species $u$ and leaves of the corresponding
tree (right) at times $t=12$ (top)  and
  $t=20$ (bottom).} \label{bbrs_fig:two_flame_c}
\end{center}
\end{figure}

In Example~3, we simulate the case with radiation,  i.e. we use \eqref{bbrs_eq24anew} and
\eqref{bbrs_stefanboltz},  where $\rho=0.05$ and $\mathrm{Le}=0.3$. Now $x_1=-5$, $x_2=5$
and $a=0.5$, $b=1$ are the respective $x-$position and initial radii of the two flame
balls. We simulate the process until $t=10$ and Figures~\ref{bbrs_fig:two_flame_b}
and~\ref{bbrs_fig:two_flame_c} show the scenario for this case. First,  the balls grow spherically,
tend to create a new flame structure, and then their fronts tend to extinguish when they
touch each other, while the radiation effect causes
the entire flame front to split. Here the maximal resolution level is set to
$L=10$ corresponding to $N_L=512^2=262144$ control volumes in the finest grid,
 and
the  reference tolerance is $\varepsilon_{\mathrm{R}}=7.43\times10^{-3}$.

Notice that the
multiresolution procedure automatically detects the higher gradient regions and
uses this information to adaptively represent the solution by the refinement and
coarsening of the mesh, i.e., by  the  adaptive addition and removal of  control
volumes on these areas.

\begin{figure}[t]
\begin{center}
\begin{tabular}{cc}
(a) & (b) \\
 \includegraphics[width=0.4\textwidth,height=0.35\textwidth]{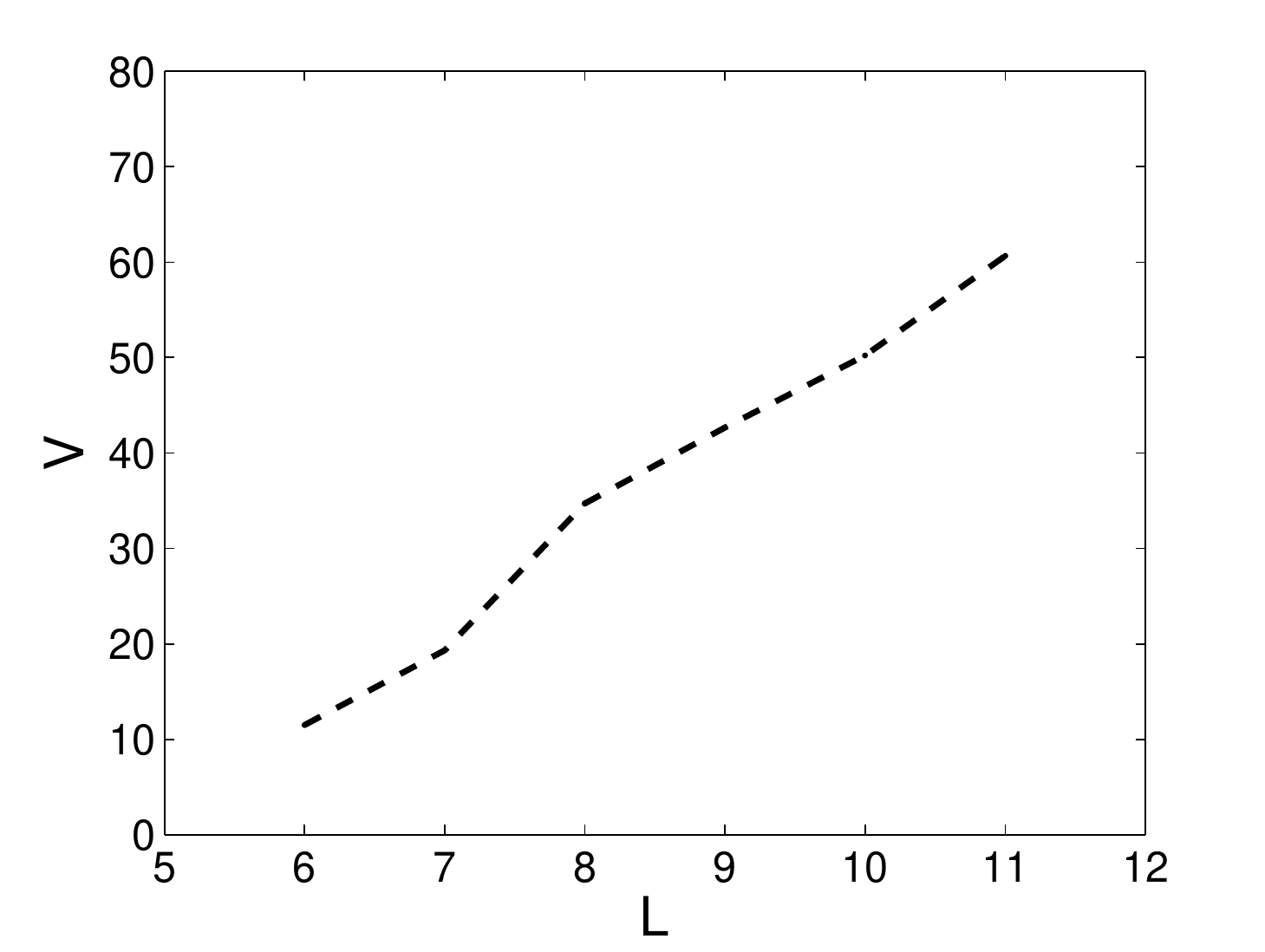}&
 \includegraphics[width=0.4\textwidth,height=0.35\textwidth]{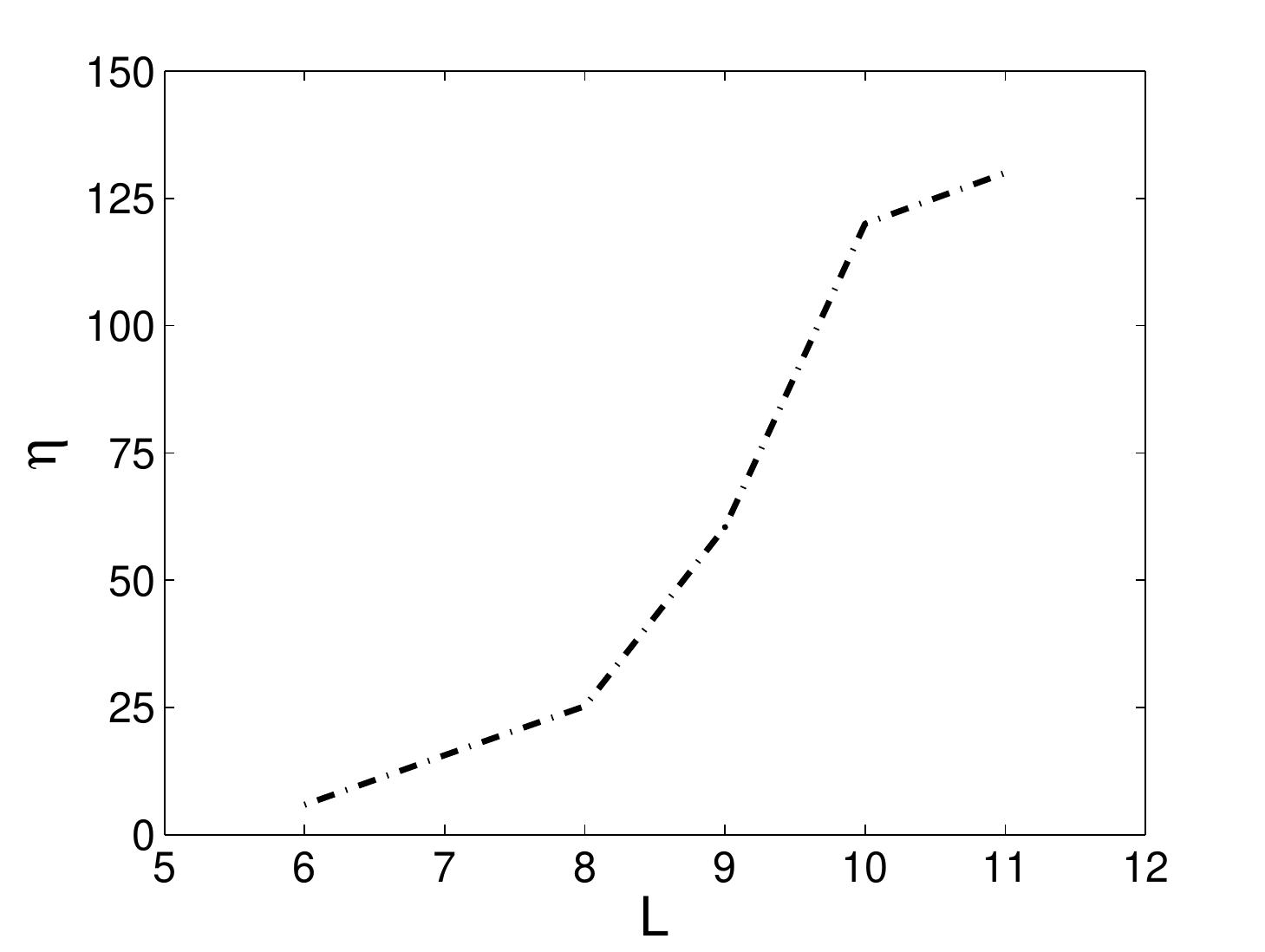}\\
(c) & (d) \\
 \includegraphics[width=0.4\textwidth,height=0.35\textwidth]{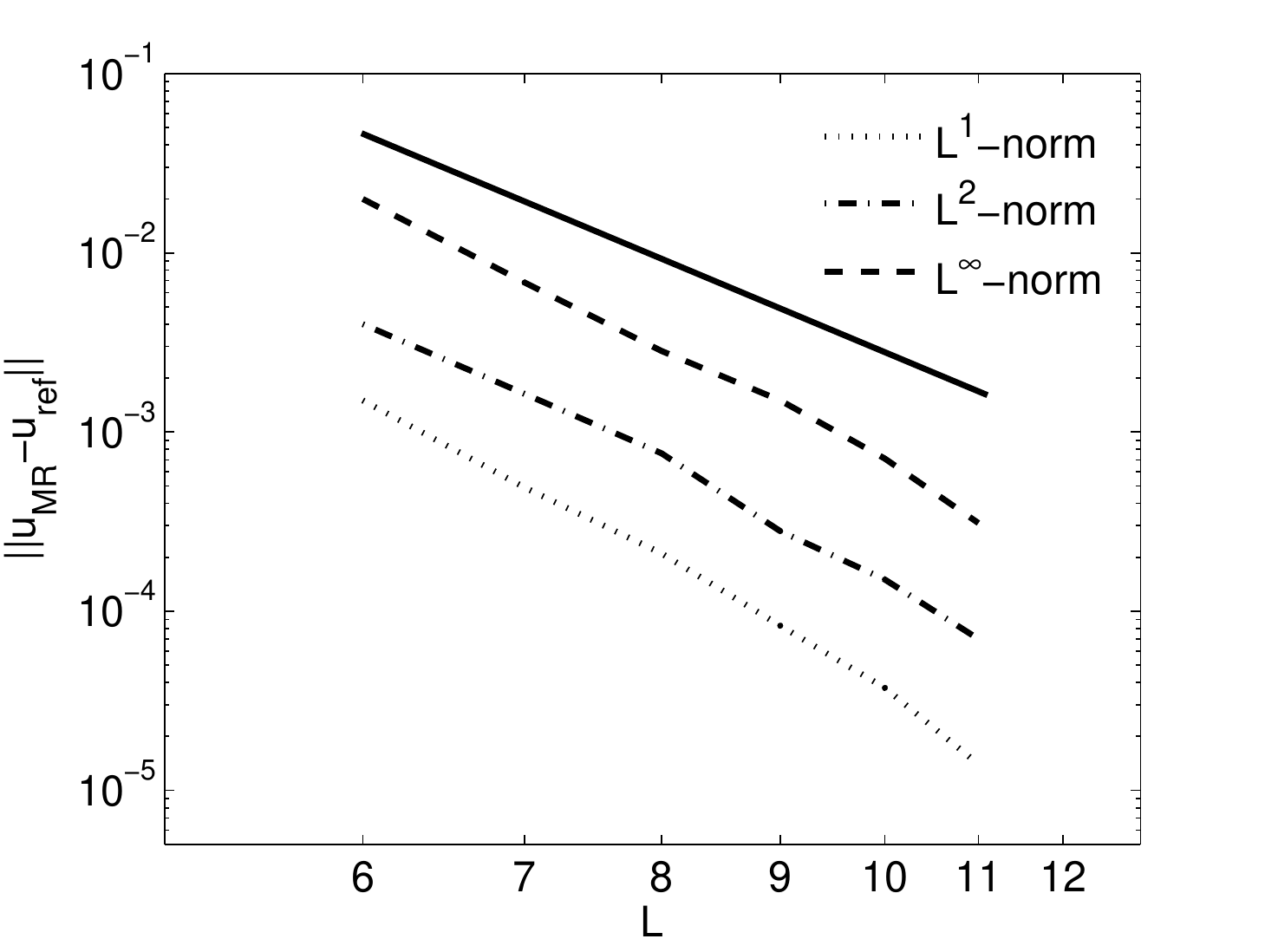}&
 \includegraphics[width=0.4\textwidth,height=0.35\textwidth]{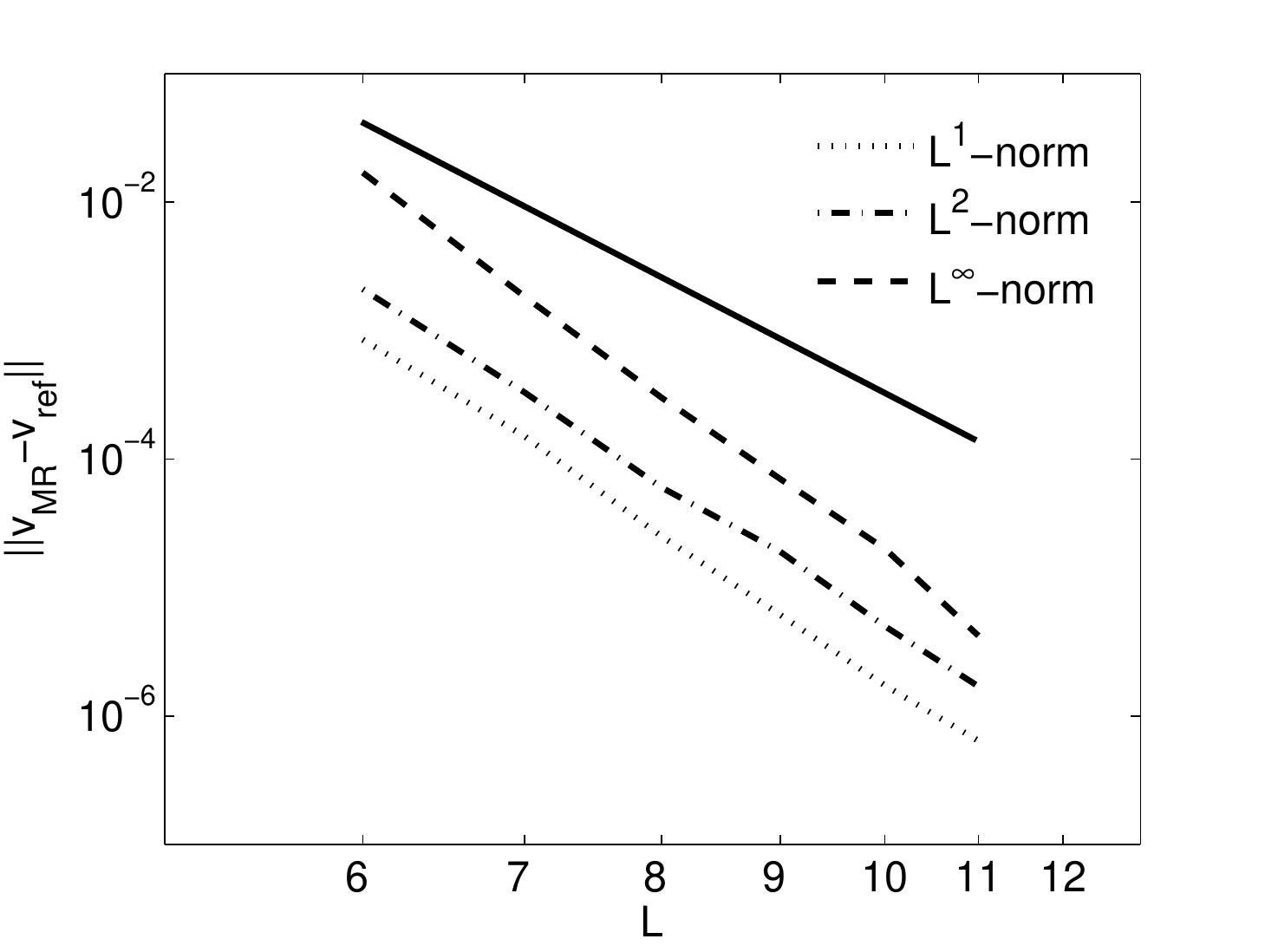}
\end{tabular}
\caption{\it Example~2 (interaction of two flame balls without radiation):
(a)  speed-up rate~$V$, (b)  data compression rate $\eta$,
for different levels, at time $t=4.0$; (c)
errors $\|\bar{u}_{\mathrm{MR}}-\bar{u}_{\mathrm{ref}}\|_1$,
$\|\bar{u}_{\mathrm{MR}}-\bar{u}_{\mathrm{ref}}\|_2$,
$\|\bar{u}_{\mathrm{MR}}-\bar{u}_{\mathrm{ref}}\|_\infty$
 and  (d)  $\|\bar{v}_{\mathrm{MR}}-\bar{v}_{\mathrm{ref}}\|_1$,
$\|\bar{v}_{\mathrm{MR}}-\bar{v}_{\mathrm{ref}}\|_2$
and $\|\bar{v}_{\mathrm{MR}}-\bar{v}_{\mathrm{ref}}\|_\infty$
respectively for different levels $L$, at time $t=4$.}
\label{bbrs_fig:ex2_cpu_eta}
\end{center}
\end{figure}

The $L^1$, $L^2$ and $L^\infty$ errors between the numerical solution
obtained by our multiresolution scheme for different multiresolution
levels $L$
and the reference solution (obtained by finite volume approximation in a uniform
fine grid with $2^{2\cdot 14}$ control volumes) for Example~2  are depicted
in Figure~\ref{bbrs_fig:ex2_cpu_eta} (c) and (d). The slopes indicate a
 rate of convergence slightly better than two.

For Example~3, we  apply the locally varying time stepping (LTS)
strategy  detailed in Section~\ref{bbrs_sec:LTS}.  We choose the maximum CFL
number allowed by \eqref{bbrs_cfl-turing}, which is $\mathsf{CFL}_0=1$ for the coarsest level.
For the remaining levels we use $\mathsf{CFL}_l=2^l\mathsf{CFL}_0$,
 which   means that we
perform each macro time step with $\Delta t=\Delta t_0=2^L\Delta t_L$ as given by
\eqref{bbrs_deltatlocal}.

\begin{figure}[t]
\begin{center}
\begin{tabular}{cc}
(a) & (b) \\
 \includegraphics[width=0.4\textwidth,height=0.35\textwidth]{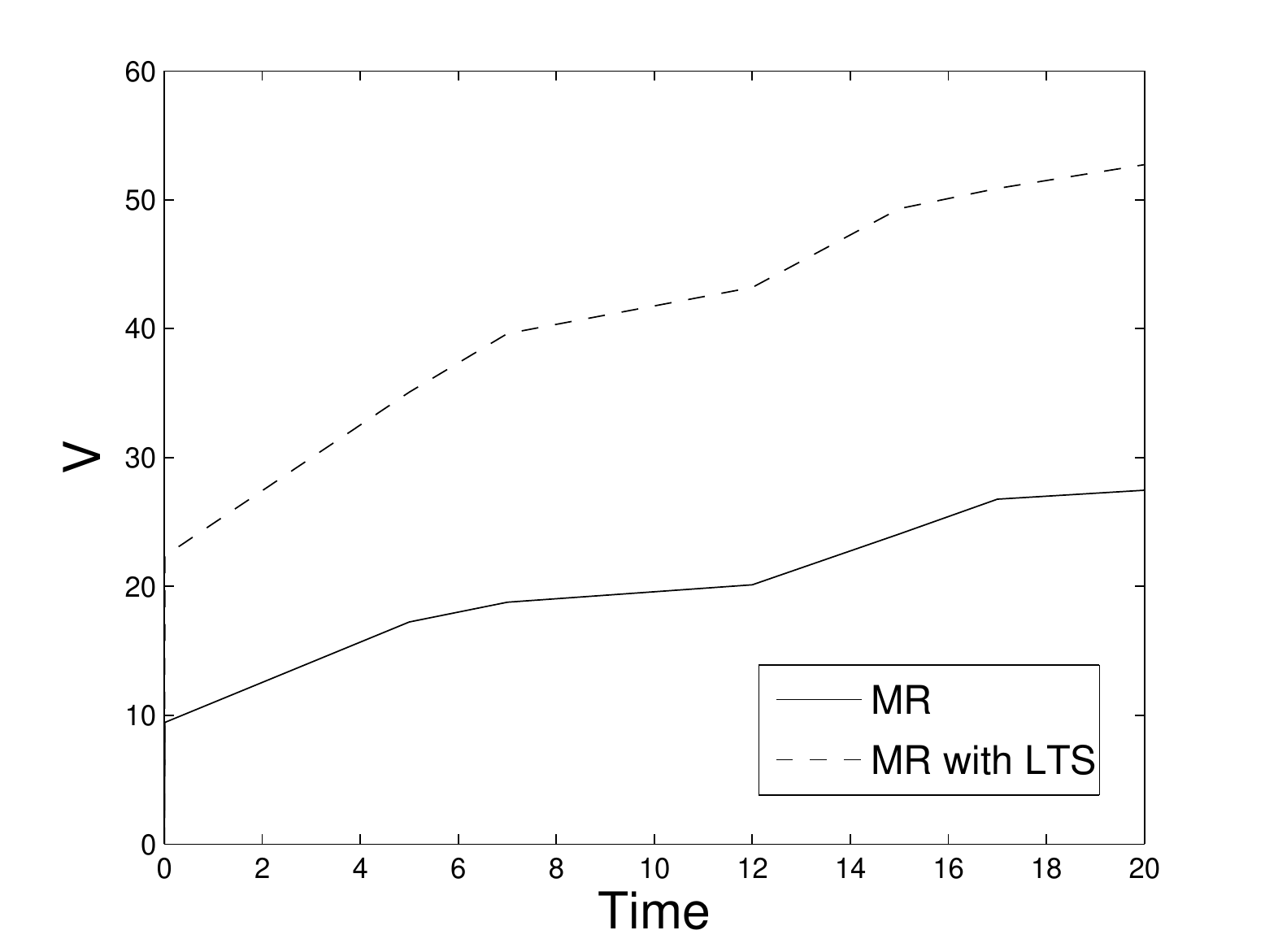}&
 \includegraphics[width=0.4\textwidth,height=0.35\textwidth]{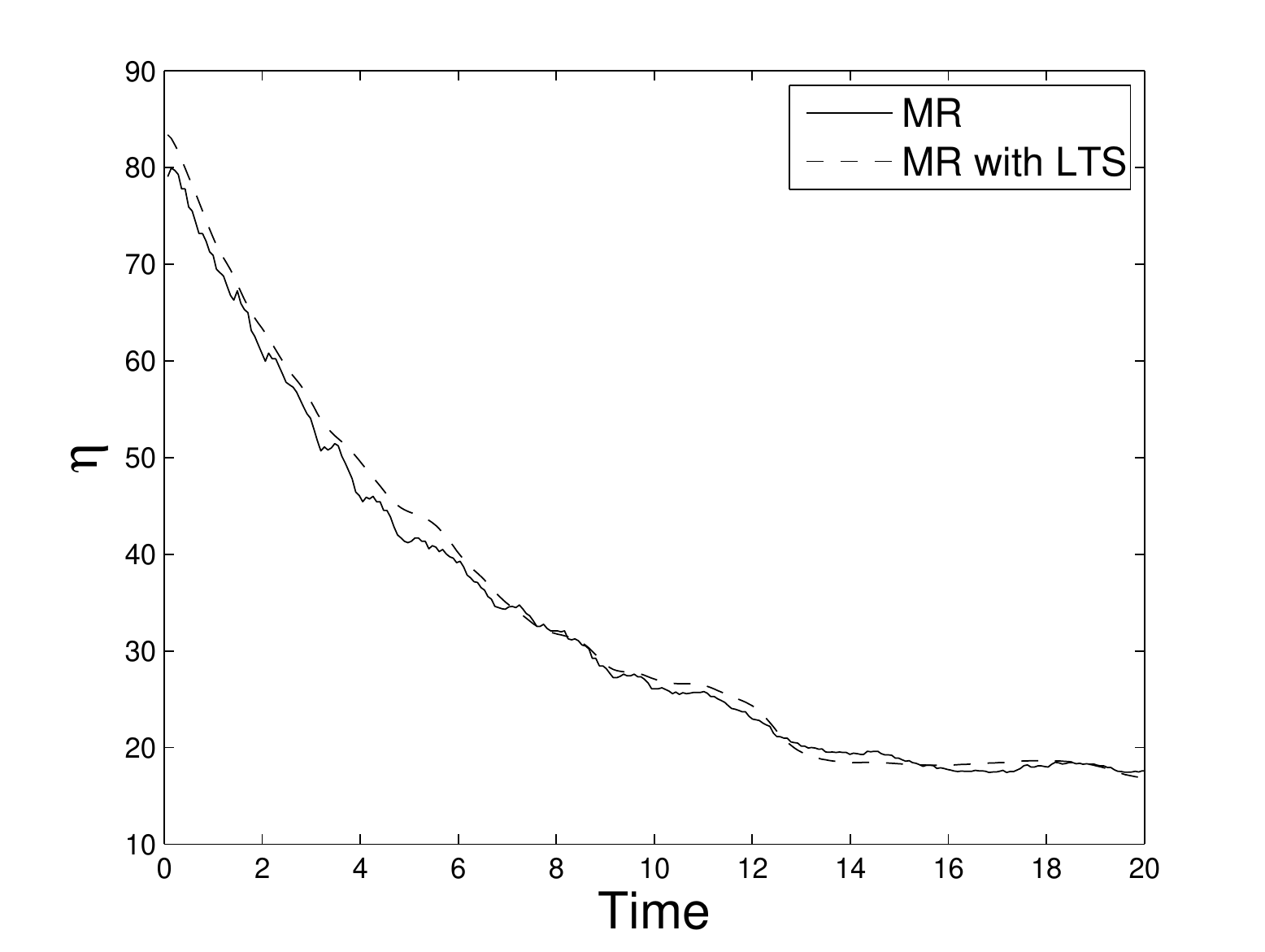}\\
(c) & (d) \\
 \includegraphics[width=0.4\textwidth,height=0.35\textwidth]{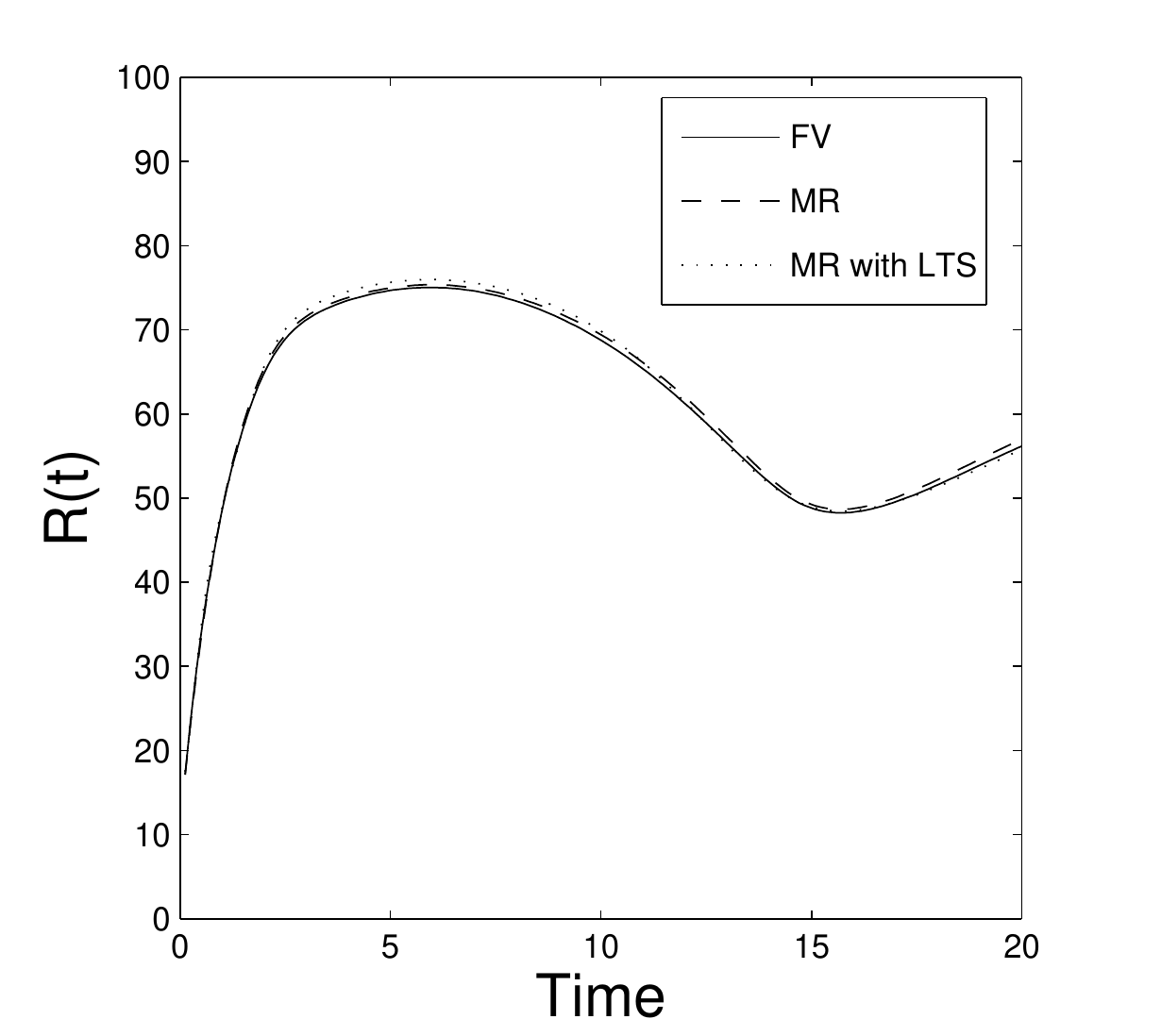}&
 \includegraphics[width=0.4\textwidth,height=0.35\textwidth]{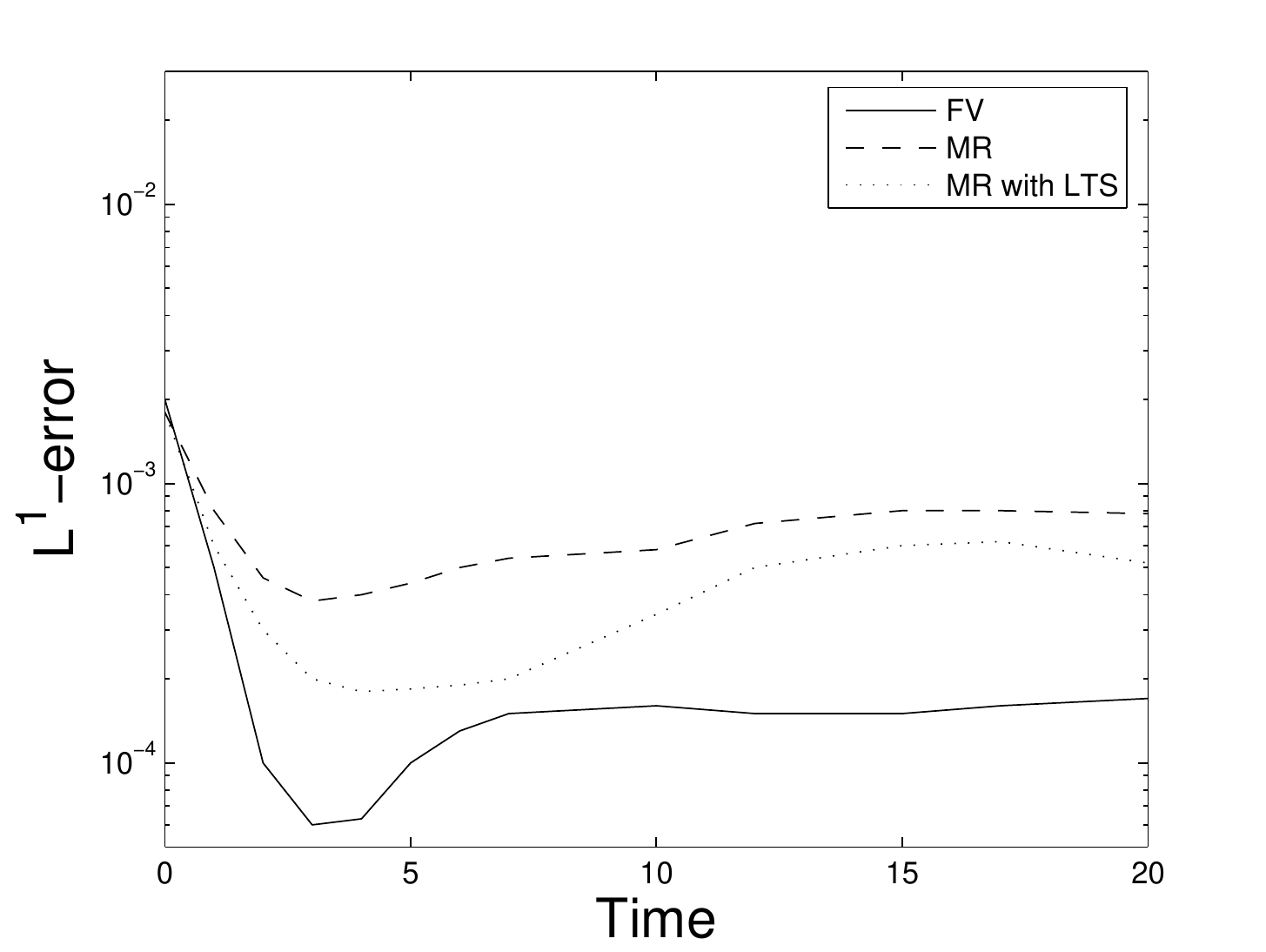}
\end{tabular}
\caption{\it Example~3 (interaction of two flame balls with radiation): Time evolution of speed-up,
data compression, and total reaction rates;
 and $L^1$-errors for different methods. $L=10$ multiresolution levels and
reference tolerance $\varepsilon_{\mathrm{R}}=7.43\times10^{-3}$.} \label{bbrs_fig:ex2b-varios}
\end{center}
\end{figure}

 In Figure~\ref{bbrs_fig:ex2b-varios} we
compare speed-up, data compression rate and total reaction rate for the finite
volume reference scheme, the multiresolution scheme with global time step, and the
multiresolution method with level-dependent time stepping. Notice that with
LTS, the  speed-up rate is approximately doubled for all times,  while the
compression rate and the total reaction rate remain  of the same order as the
multiresolution computation with global time step.

\begin{table}[t]
\begin{center}
\begin{tabular}{|l|c|c|c|c|c|c|c|r|}
\hline
Time&  $V$& $\eta$& $L^1-$error & $L^2-$error &$L^\infty-$error& Method & $R(t)$  $\vphantom{\int_X^X}$ \\
\hline
\hline
$t=2.0$
& 12.47& 138.2613& $5.41\times10^{-3}$&$5.77\times10^{-3}$&$2.46\times10^{-2}$& MR & 56.7230  $\vphantom{\int_X^X}$\\
&      &         & &&& FV & 56.0078 $\vphantom{\int_X^X}$\\
\hline
$t=4.0$
& 20.56& 113.4331&$6.39\times10^{-3}$&$8.42\times10^{-4}$&$3.02\times10^{-2}$& MR   & 80.0374 $\vphantom{\int_X^X}$\\
&      &         & &&& FV     & 79.5247$\vphantom{\int_X^X}$\\
\hline
$t=10.0$
 & 34.42& 83.9129 & $5.20\times10^{-3}$&$4.90\times10^{-3}$&$5.49\times10^{-2}$& MR  & 98.9210 $\vphantom{\int_X^X}$\\
&      &         &&&& FV     & 98.7942 $\vphantom{\int_X^X}$\\
\hline
\end{tabular}
\end{center}

\vspace*{2mm}

\caption{\it Example~2 (interaction of two flame balls without radiation):  Corresponding simulated
time, speed-up rate~$V$, compression rate~$\eta$, errors, and total reaction rate~$R(t)$ for the $u$ species. }
\label{bbrs_table:example_2}
\end{table}

\begin{figure}[t]
\begin{center}
\begin{tabular}{cc}
\includegraphics[width=0.37\textwidth]{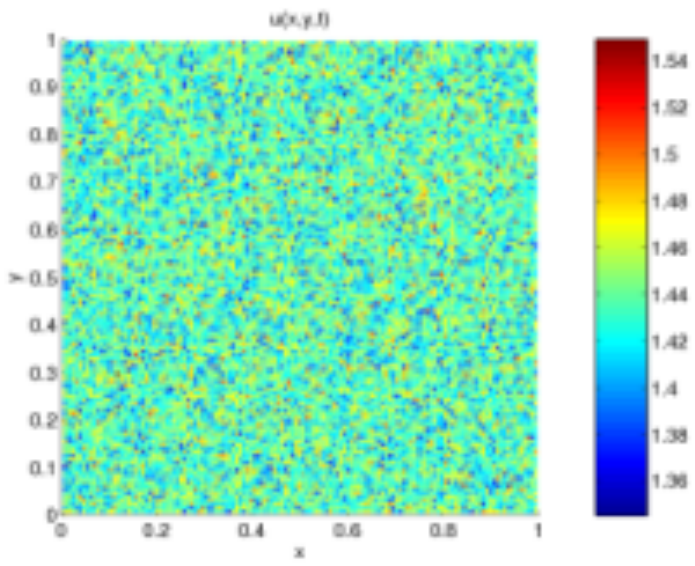}&
\includegraphics[width=0.37\textwidth]{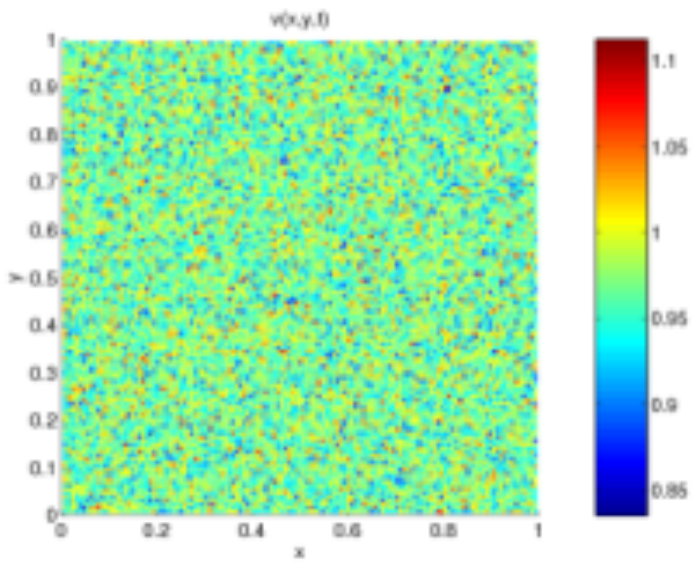}
\end{tabular}
\caption{\it Examples~4 and~5 (Turing pattern-formation):  Initial data
$u_0(x,y)$ (left) and $v_0(x,y)$ (right).}
\label{bbrs_fig:initial_data}
\end{center}
\end{figure}

\subsection{Examples~4 and~5: A Turing model of pattern formation}
We select the parameters
$a=-0.5$, $b=1.9$, $d=4.8$ and $\gamma=210$.
According to the discussion of  Section~\ref{bbrs_sec:2species}, these
parameters allow diffusion-driven instabilities to evolve.
 The  initial  concentration distribution is a  normally distributed random
 perturbation around the
stationary state $(u^0,v^0)$ for the non-degenerate case, with a variance
lower than the amplitude of the final patterns, see  Figure~\ref{bbrs_fig:initial_data}.
For the case of non-degenerate diffusion (Example~4),
we use $A(u)$ and $B(u)$ as given by  \eqref{bbrs_aubu}. For these
parameters, the steady state is  $(u^0=1.4,  v^0=0.96939)$.

\begin{figure}[t]
\begin{center}
\begin{tabular}{cc}
(a) & (b) \\
\includegraphics[width=0.4\textwidth,height=0.35\textwidth]{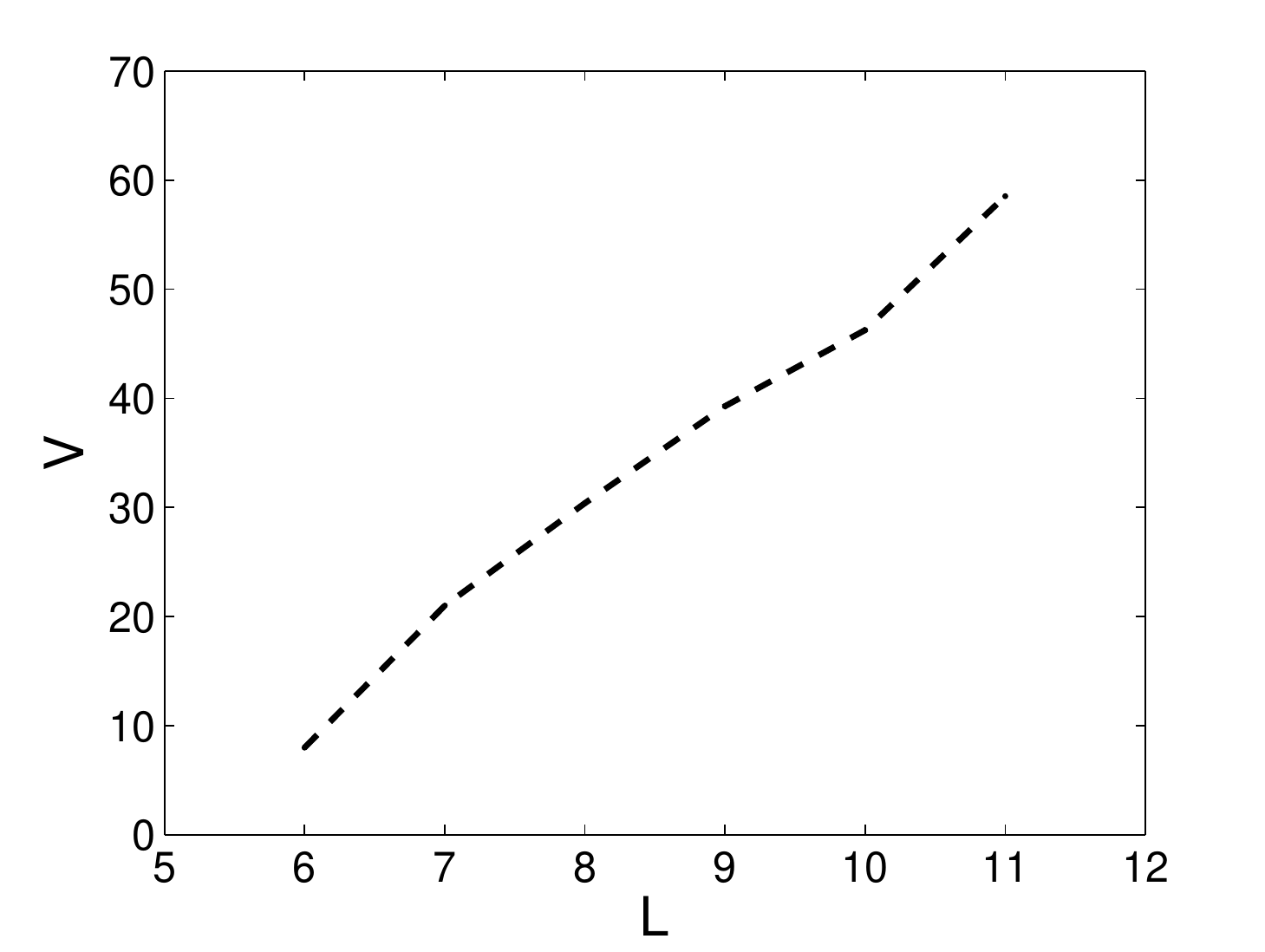}&
\includegraphics[width=0.4\textwidth,height=0.35\textwidth]{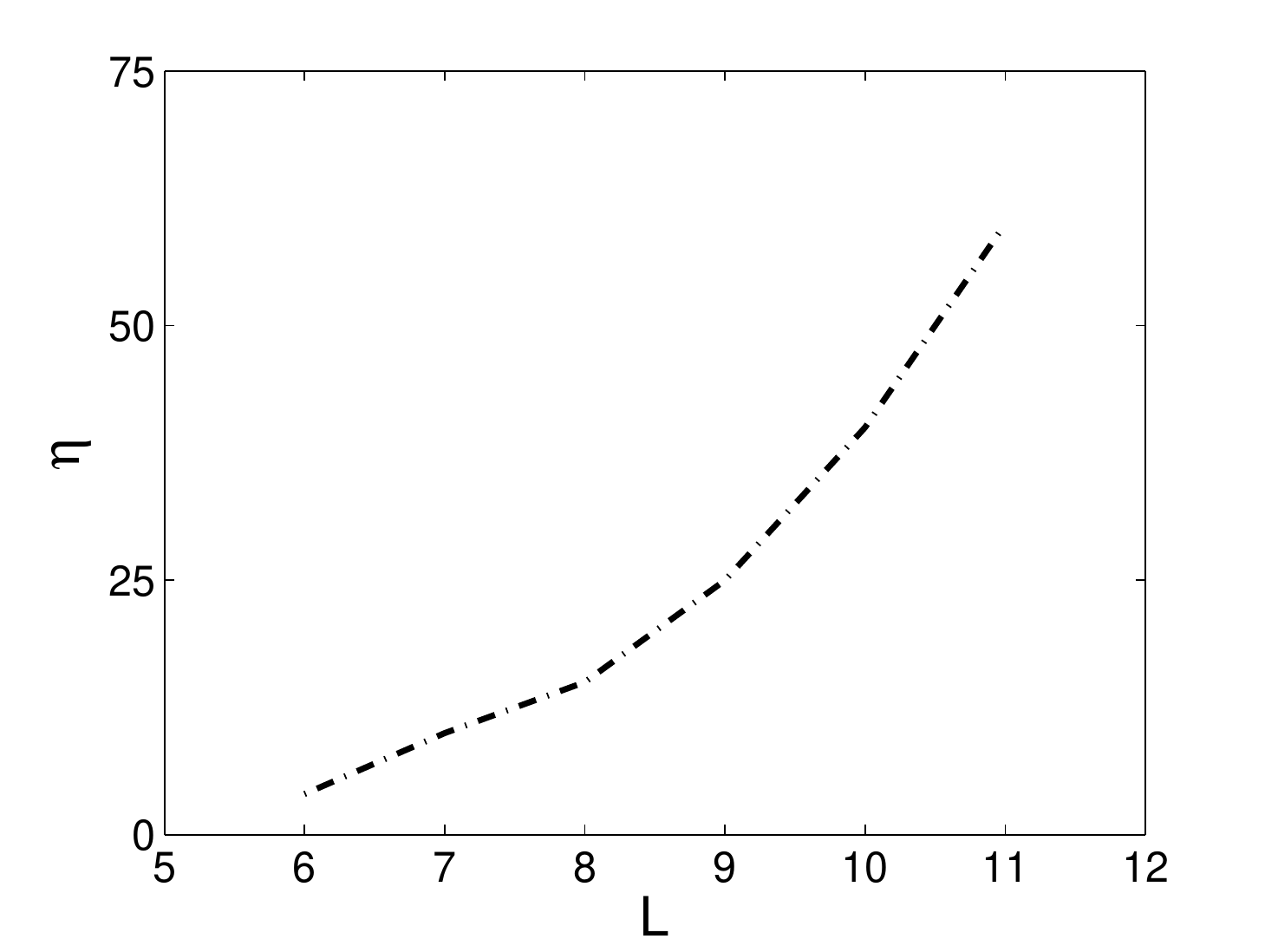}\\
(c) & (d) \\
\includegraphics[width=0.4\textwidth,height=0.35\textwidth]{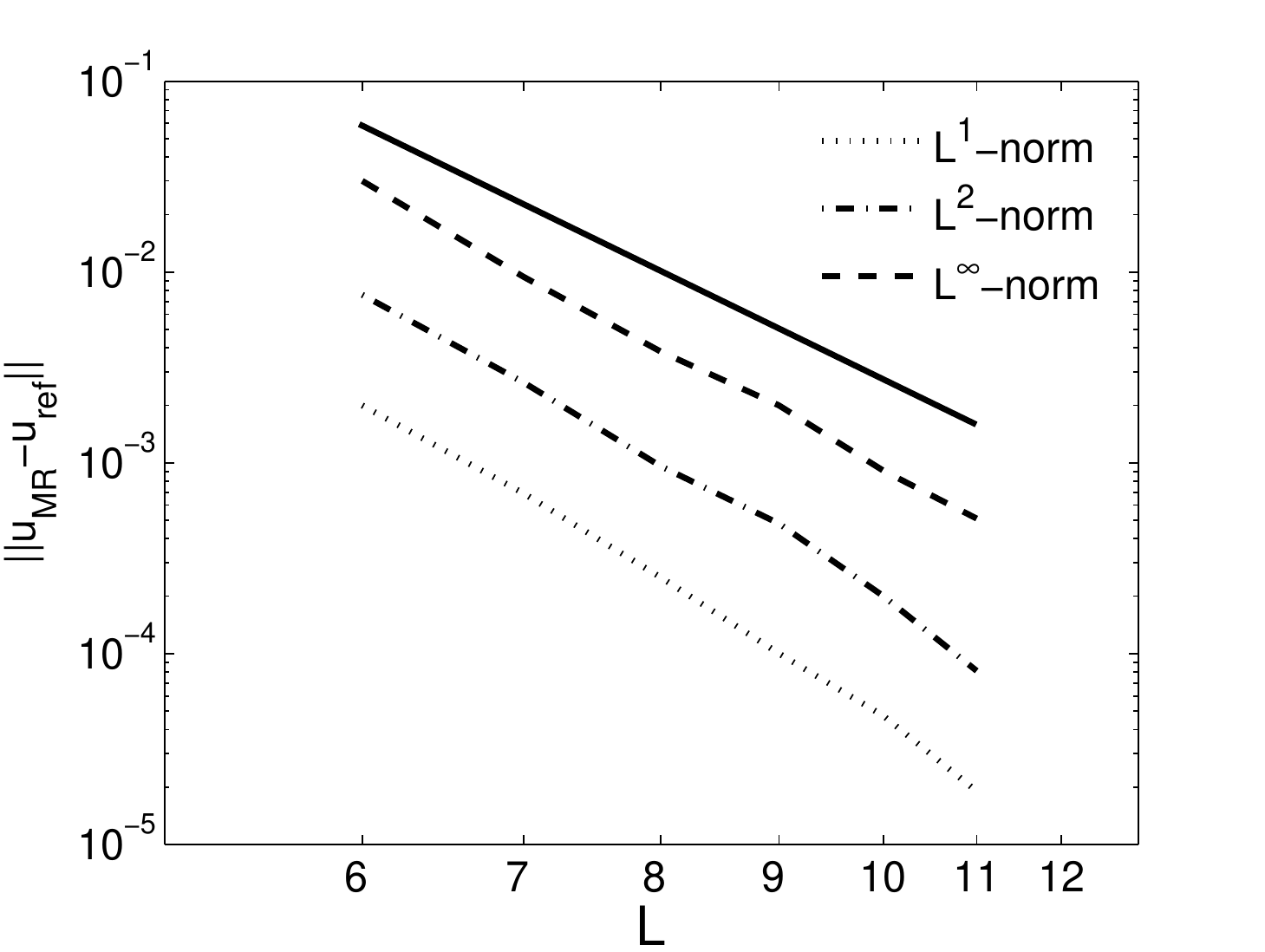}&
\includegraphics[width=0.4\textwidth,height=0.35\textwidth]{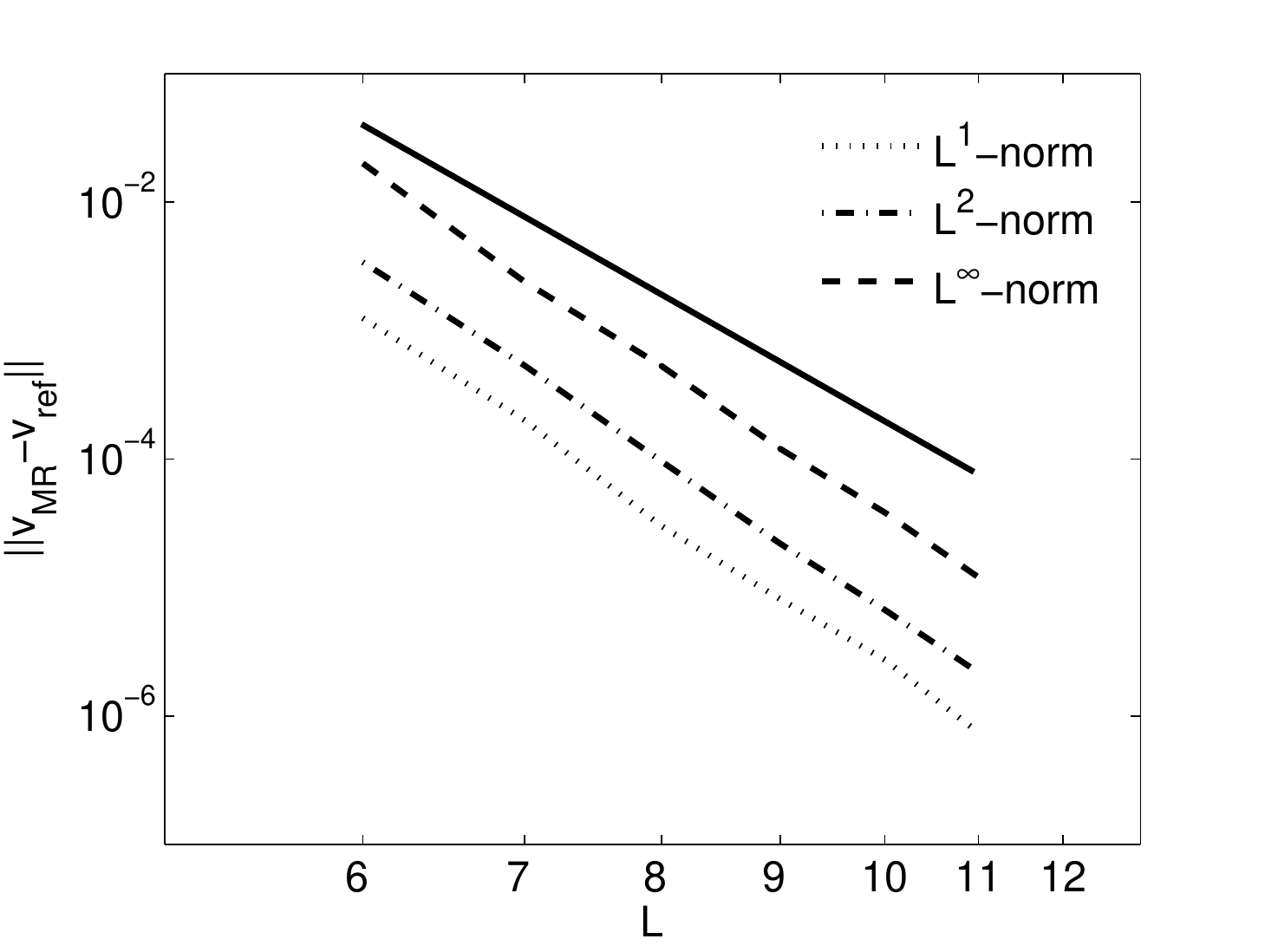}
\end{tabular}
\caption{\it Example~4 (Model~2 with non-degenerate diffusion):
(a) CPU acceleration $V$, (b)
 data compression rate $\eta$,
for different levels, at time $t=0.75\,\mathrm{s}$; (c)
errors $\|\bar{u}_{\mathrm{MR}}-\bar{u}_{\mathrm{ref}}\|_1$,
$\|\bar{u}_{\mathrm{MR}}-\bar{u}_{\mathrm{ref}}\|_2$,
$\|\bar{u}_{\mathrm{MR}}-\bar{u}_{\mathrm{ref}}\|_\infty$
 and  (d)  $\|\bar{v}_{\mathrm{MR}}-\bar{v}_{\mathrm{ref}}\|_1$,
$\|\bar{v}_{\mathrm{MR}}-\bar{v}_{\mathrm{ref}}\|_2$
and $\|\bar{v}_{\mathrm{MR}}-\bar{v}_{\mathrm{ref}}\|_\infty$
respectively for different levels $L$, at time $t=0.75$. }
\label{bbrs_fig:ex3a_cpu_eta}
\end{center}
\end{figure}

In Example~4 we choose a maximal resolution level of $N_L=256^2=65536$ control
volumes in the finest grid and a reference tolerance given by
$\varepsilon_{\mathrm{R}} = 2.6\times10^{-3}$.
The time step is the maximum allowed by the CFL condition \eqref{bbrs_cfl-turing}.
Table~\ref{bbrs_table:example_3a} summarizes the speed-up rate, compression rate
and errors in different norms between the numerical solution by multiresolution
and the fine-mesh finite volume reference solution for different times.
We depict errors between our multiresolution scheme and a reference FV solution with
$N_L=1024^2=1048576$ control volumes in the finest grid, for different multiresolution
levels $L$ in Figure~\ref{bbrs_fig:ex3a_cpu_eta} (c) and (d). In this case,
the slopes  equally indicate a rate of convergence slightly larger than two.
 Concerning the computation of errors,
in Examples~4, 5 and~6 the system is evolved until the ``random
noise'',  which is imposed as an initial condition on the finest grid,
  has been smoothed sufficiently; then, this solution is  projected  on
coarser levels to obtain auxiliary initial conditions for all the
 needed levels.

\begin{table}[t]
\begin{center}
\begin{tabular}{|l|c|c|c|c|c|c|c|}
\hline
Time & $V$  & $\eta$&  Species &  $L^1-$error  & $L^2-$error    &$L^\infty-$error  $\vphantom{\int_X^X}$   \\
\hline
\hline
$t=0.05 $
& 7.16 & 11.3783 & $u$ & $6.81\times10^{-4}$&$4.76\times10^{-5}$&$3.46\times10^{-3}$ $\vphantom{\int_X^X}$\\
&      &         & $v$ & $4.09\times10^{-4}$&$3.92\times10^{-4}$&$5.38\times10^{-4}$ $\vphantom{\int_X^X}$\\
\hline
$t=0.25$
& 9.29 &11.9756 & $u$& $8.37\times10^{-4}$&$6.94\times10^{-5}$&$9.93\times10^{-3}$ $\vphantom{\int_X^X}$\\
&     &         & $v$& $4.22\times10^{-4}$&$5.43\times10^{-4}$&$8.48\times10^{-4}$ $\vphantom{\int_X^X}$\\
\hline
$t=1.50$
& 11.87 & 14.4739 & $u$& $9.26\times10^{-4}$&$2.71\times10^{-4}$&$2.44\times10^{-2}$ $\vphantom{\int_X^X}$\\
&       &         & $v$& $4.30\times10^{-4}$&$9.77\times10^{-5}$&$8.39\times10^{-3}$ $\vphantom{\int_X^X}$\\
\hline
\end{tabular}
\end{center}

\vspace*{2mm}

\caption{\it Example~4 (Model~2 with non-degenerate diffusion):  Corresponding simulated time,
CPU ratio~$V$, compression rate~$\eta$  and componentwise errors.} \label{bbrs_table:example_3a}
\end{table}

\begin{figure}[t]
\begin{center}
\begin{tabular}{ccc}
\includegraphics[width=0.32\textwidth,height=0.3\textwidth]{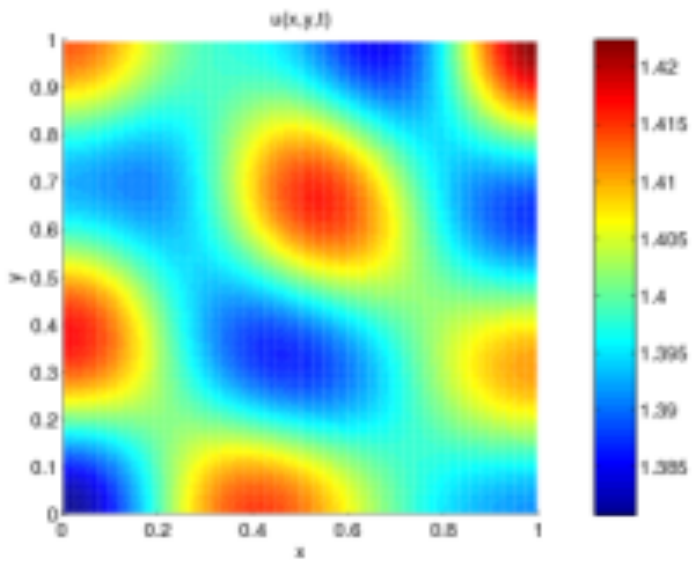}&
\includegraphics[width=0.32\textwidth,height=0.3\textwidth]{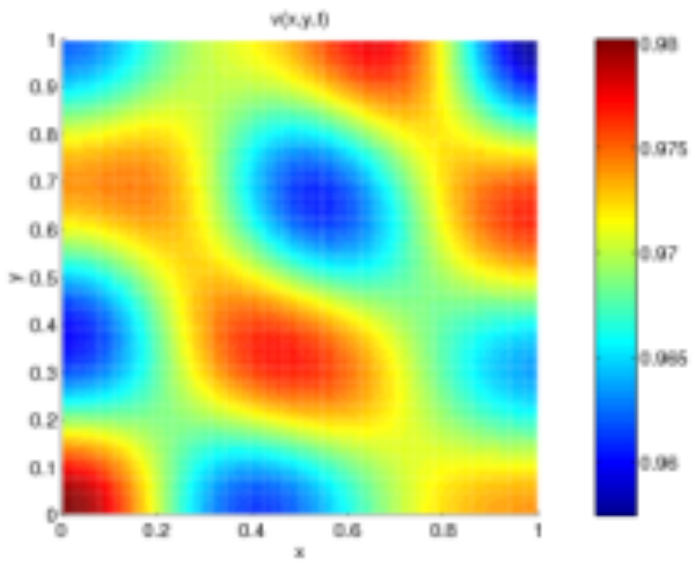}&
\includegraphics[width=0.28\textwidth,height=0.3\textwidth]{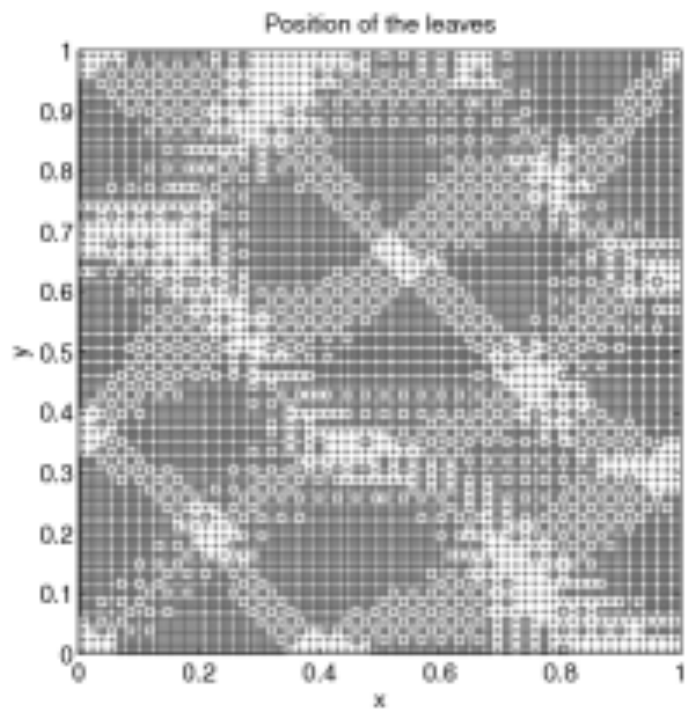}\\
\includegraphics[width=0.32\textwidth,height=0.3\textwidth]{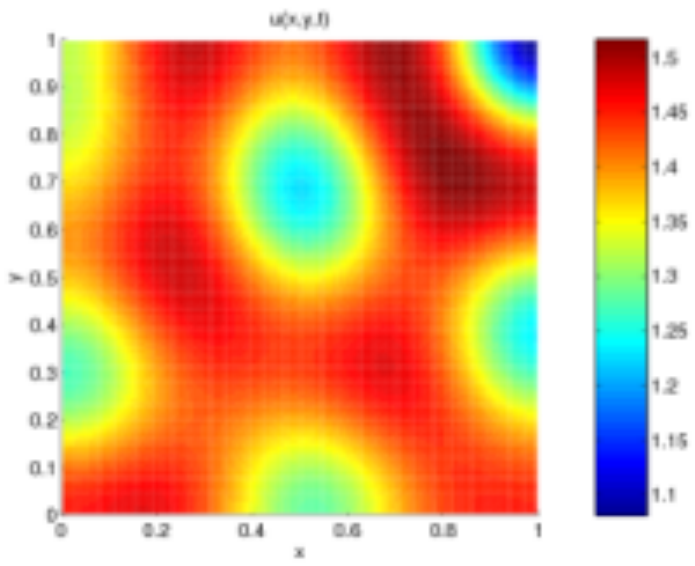}&
\includegraphics[width=0.32\textwidth,height=0.3\textwidth]{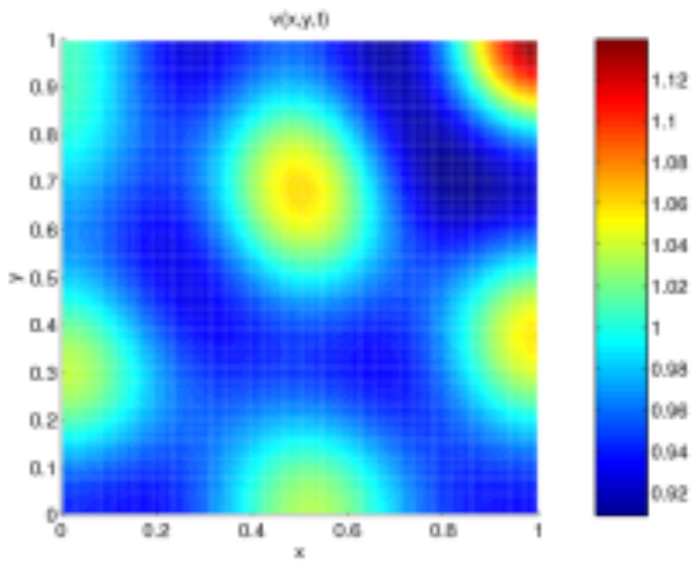}&
\includegraphics[width=0.28\textwidth,height=0.3\textwidth]{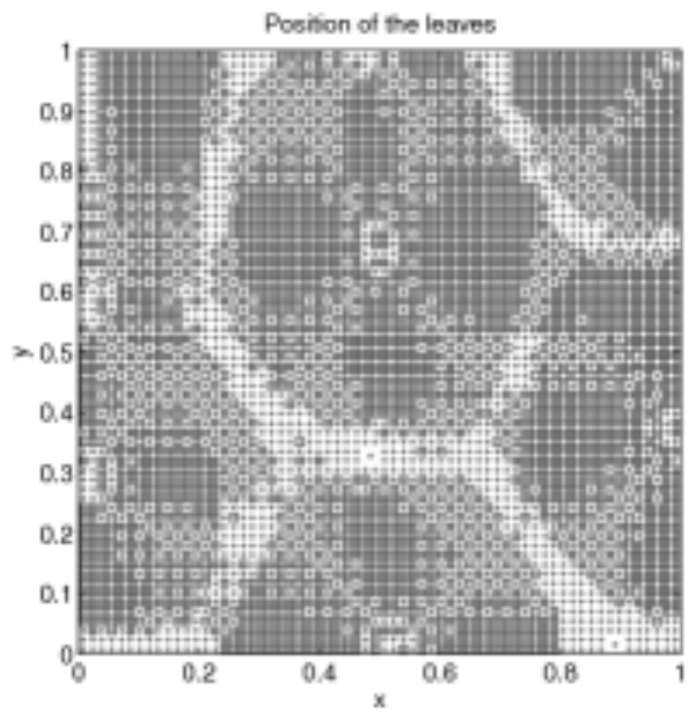}\\
\includegraphics[width=0.32\textwidth,height=0.3\textwidth]{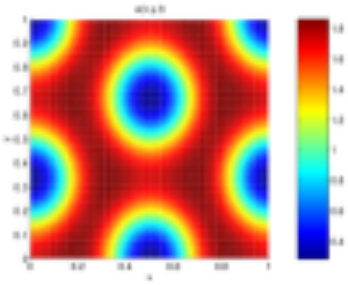}&
\includegraphics[width=0.32\textwidth,height=0.3\textwidth]{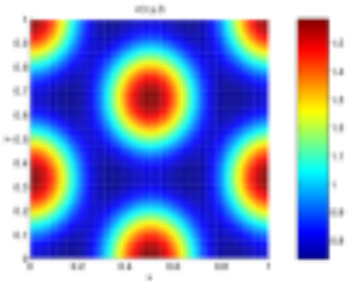}&
\includegraphics[width=0.28\textwidth,height=0.3\textwidth]{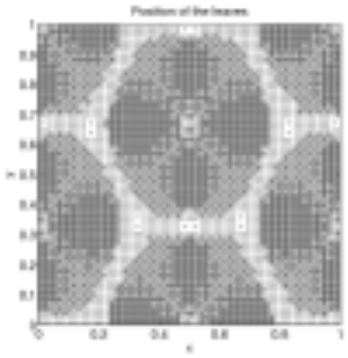}
\end{tabular}
\caption{\it Example~4 (Model~2 with non-degenerate diffusion):
Numerical solution for species $u$ (left) and $v$ (right),
 and leaves of the corresponding
tree data structure at times $t=0.05$ (top), $t=0.25$ (middle)
and $t=1.5$ (bottom).} \label{bbrs_fig:non-deg-dif}
\end{center}
\end{figure}

\begin{figure}[t]
\begin{center}
\begin{tabular}{ccc}
\includegraphics[width=0.32\textwidth,height=0.3\textwidth]{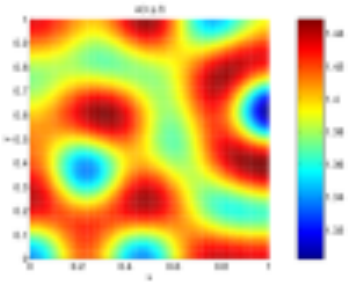}&
\includegraphics[width=0.32\textwidth,height=0.3\textwidth]{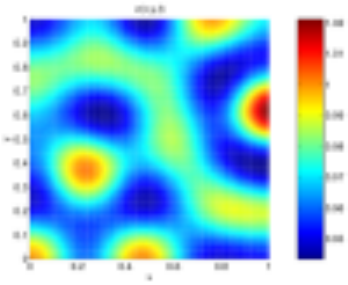}&
\includegraphics[width=0.28\textwidth,height=0.3\textwidth]{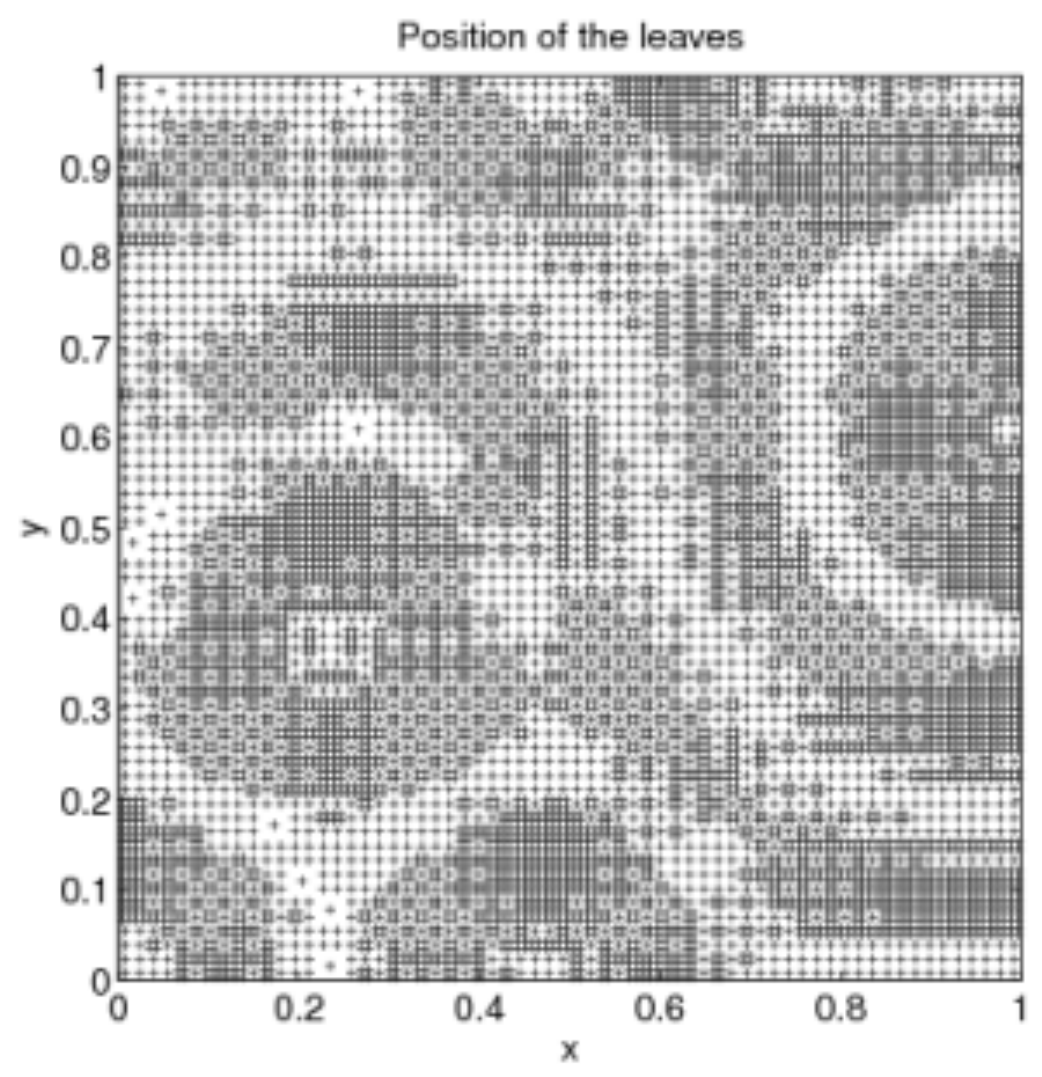}\\
\includegraphics[width=0.32\textwidth,height=0.3\textwidth]{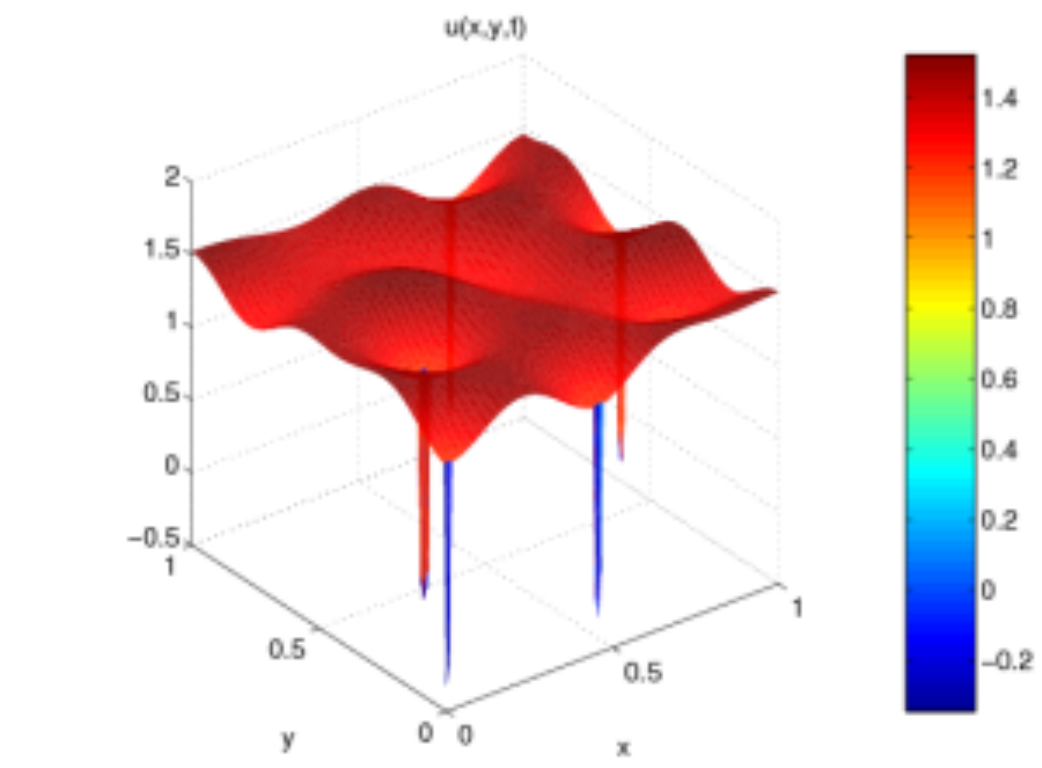}&
\includegraphics[width=0.32\textwidth,height=0.3\textwidth]{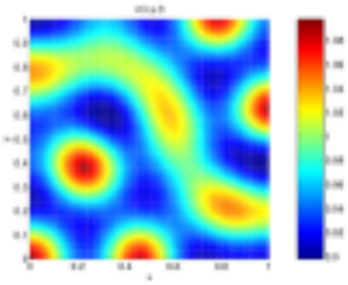}&
\includegraphics[width=0.28\textwidth,height=0.3\textwidth]{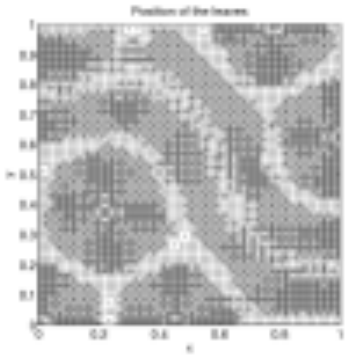}\\
\includegraphics[width=0.32\textwidth,height=0.3\textwidth]{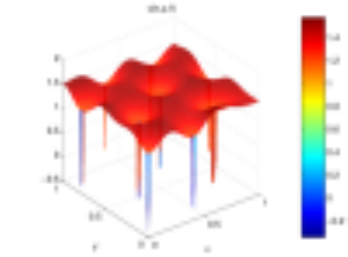}&
\includegraphics[width=0.32\textwidth,height=0.3\textwidth]{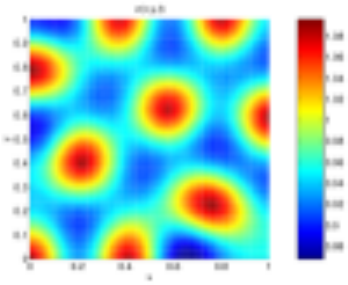}&
\includegraphics[width=0.28\textwidth,height=0.3\textwidth]{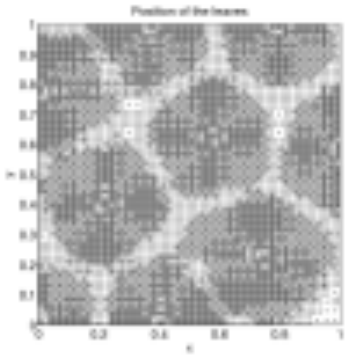}
\end{tabular}
\caption{\it Example~5 (Model~2 with degenerate diffusion):   Numerical solution
for species $u$ (left) and $v$ (middle), and leaves of the corresponding
tree data structure (right) at times $t=0.1$ (top), $t=0.25$ (middle) and
$t=1.5$ (bottom).}
\label{bbrs_fig:deg-dif}
\end{center}
\end{figure}

For Example~5, we use the degenerate diffusion coefficients
\eqref{bbrs_degen_param} with $ u_{\mathrm{c}}=1.2$ and
$v_{\mathrm{c}}=0.7$, and employ again the kinetics
 \eqref{bbrs_fgkinetics}, but this time we choose the parameters
$a=-0.5$, $b=1.9$, $d=4.8$ and $\gamma=395$.
We select a maximal resolution level of $N_L=256^2=65536$
control volumes in the finest grid, with a reference tolerance given by
$\varepsilon_{\mathrm{R}}=3.59\times10^{-4}$. From Table~\ref{bbrs_table:example_3b} we
see that  the multiresolution algorithm
 allows significant acceleration and
 data  compression rate are significantly increased by
the multiresolution algorithm with very good accuracy.
Figure~\ref{bbrs_fig:deg-dif} indicates that due to the
degeneracy of the diffusion given by \eqref{bbrs_degen_param},
 and in contrast to Example~4, species~$u$ exhibits
patterns with steeper gradients, and especially at $t=0.25$ and
$t=1.5$,  singularities appear.

\begin{table}[ht]
\begin{center}
\begin{tabular}{|l|c|c|c|c|c|c|c|}
\hline
Time & $V$  & $\eta$&  Species &    $L^1-$error     & $L^2-$error        &$L^\infty-$error  $\vphantom{\int_X^X}$   \\
\hline
\hline
$t=0.10$& 6.32 & 12.5438 & $u$& $6.31\times10^{-4}$&$5.82\times10^{-4}$&$2.72\times10^{-3}$ $\vphantom{\int_X^X}$\\
        &      &         & $v$& $4.98\times10^{-4}$&$5.37\times10^{-4}$&$9.46\times10^{-4}$ $\vphantom{\int_X^X}$\\
\hline
$t=0.25$& 9.79 & 10.3457 & $u$& $6.12\times10^{-4}$&$2.46\times10^{-5}$&$3.03\times10^{-3}$ $\vphantom{\int_X^X}$\\
        &      &         & $v$& $3.91\times10^{-4}$&$9.22\times10^{-4}$&$9.92\times10^{-4}$ $\vphantom{\int_X^X}$\\
\hline
$t=1.50$& 11.60& 10.1984 & $u$& $3.42\times10^{-4}$&$7.34\times10^{-4}$&$3.40\times10^{-3}$ $\vphantom{\int_X^X}$\\
        &      &         & $v$& $2.63\times10^{-4}$&$4.98\times10^{-4}$&$2.81\times10^{-3}$ $\vphantom{\int_X^X}$\\
\hline
\end{tabular}
\end{center}

\vspace*{2mm}

\caption{\it Example~5 (Model~2 with degenerate diffusion):
Corresponding simulated time,
speed-up rate~$V$, compression rate~$\eta$ and componentwise errors.} \label{bbrs_table:example_3b}
\end{table}

\subsection{Example~6: Chemotaxis-growth system}\label{bbrs_sec:int:ks}
\begin{figure}[t]
\begin{center}
\begin{tabular}{ccc}
\includegraphics[width=0.32\textwidth,height=0.3\textwidth]{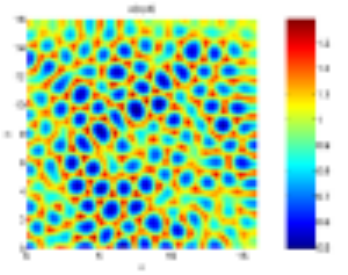}&
\includegraphics[width=0.32\textwidth,height=0.3\textwidth]{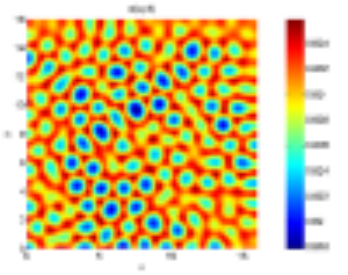}&
\includegraphics[width=0.28\textwidth,height=0.3\textwidth]{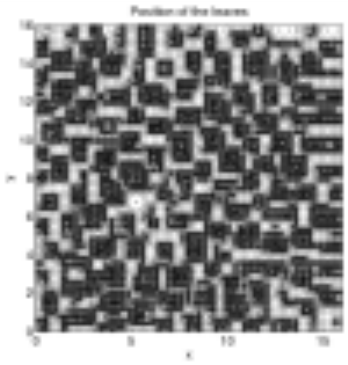}\\
\includegraphics[width=0.32\textwidth,height=0.3\textwidth]{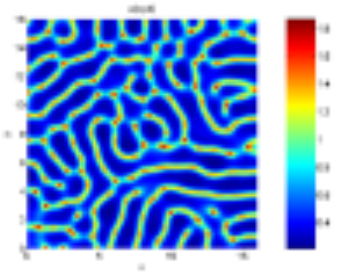}&
\includegraphics[width=0.32\textwidth,height=0.3\textwidth]{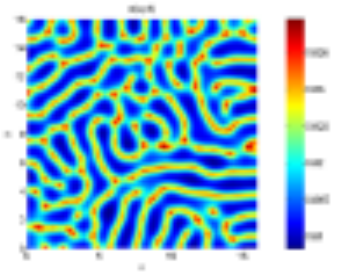}&
\includegraphics[width=0.28\textwidth,height=0.3\textwidth]{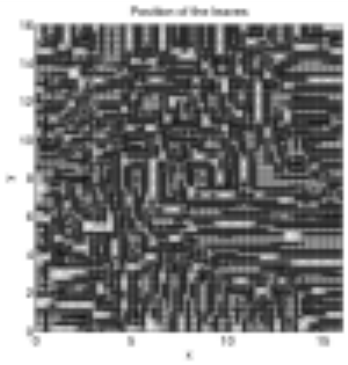}\\
\includegraphics[width=0.32\textwidth,height=0.3\textwidth]{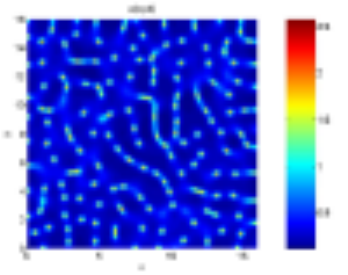}&
\includegraphics[width=0.32\textwidth,height=0.3\textwidth]{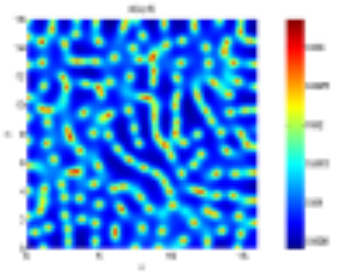}&
\includegraphics[width=0.28\textwidth,height=0.3\textwidth]{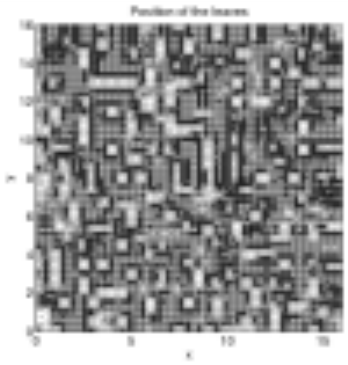}
\end{tabular}
\caption{\it Example~6 (chemotaxis with growth):
Numerical solution for species $u$ (left) and $v$ (middle), and leaves of the corresponding
tree data structure (right) for $\nu=6.04$ (top), $\nu=7$ (middle)  and
  $\nu=10$ (bottom).} \label{bbrs_fig:ex4-chemo}
\end{center}
\end{figure}

\begin{figure}[t]
\begin{center}
\begin{tabular}{cc}
(a) & (b) \\
\includegraphics[width=0.4\textwidth,height=0.35\textwidth]{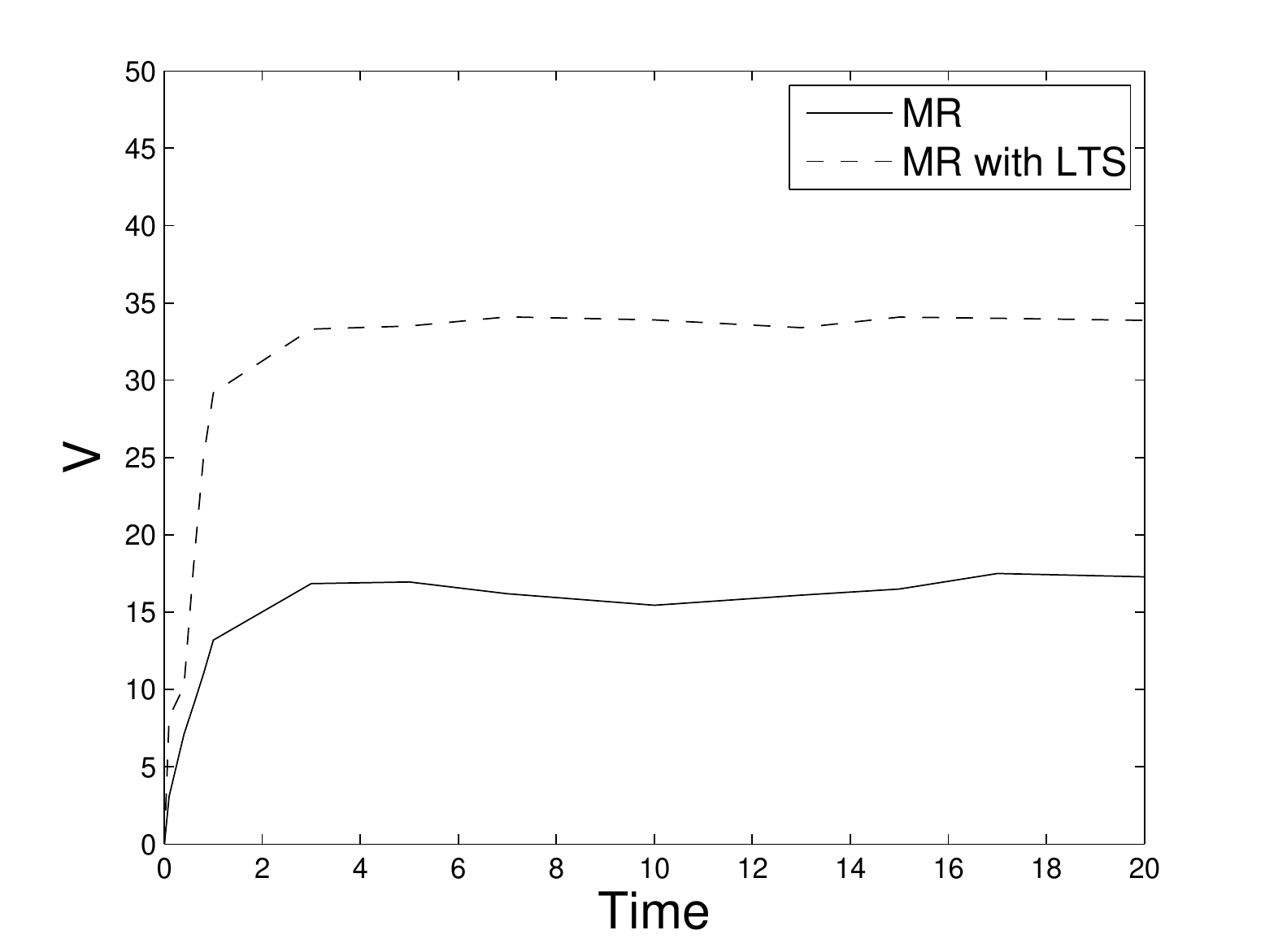}&
\includegraphics[width=0.4\textwidth,height=0.35\textwidth]{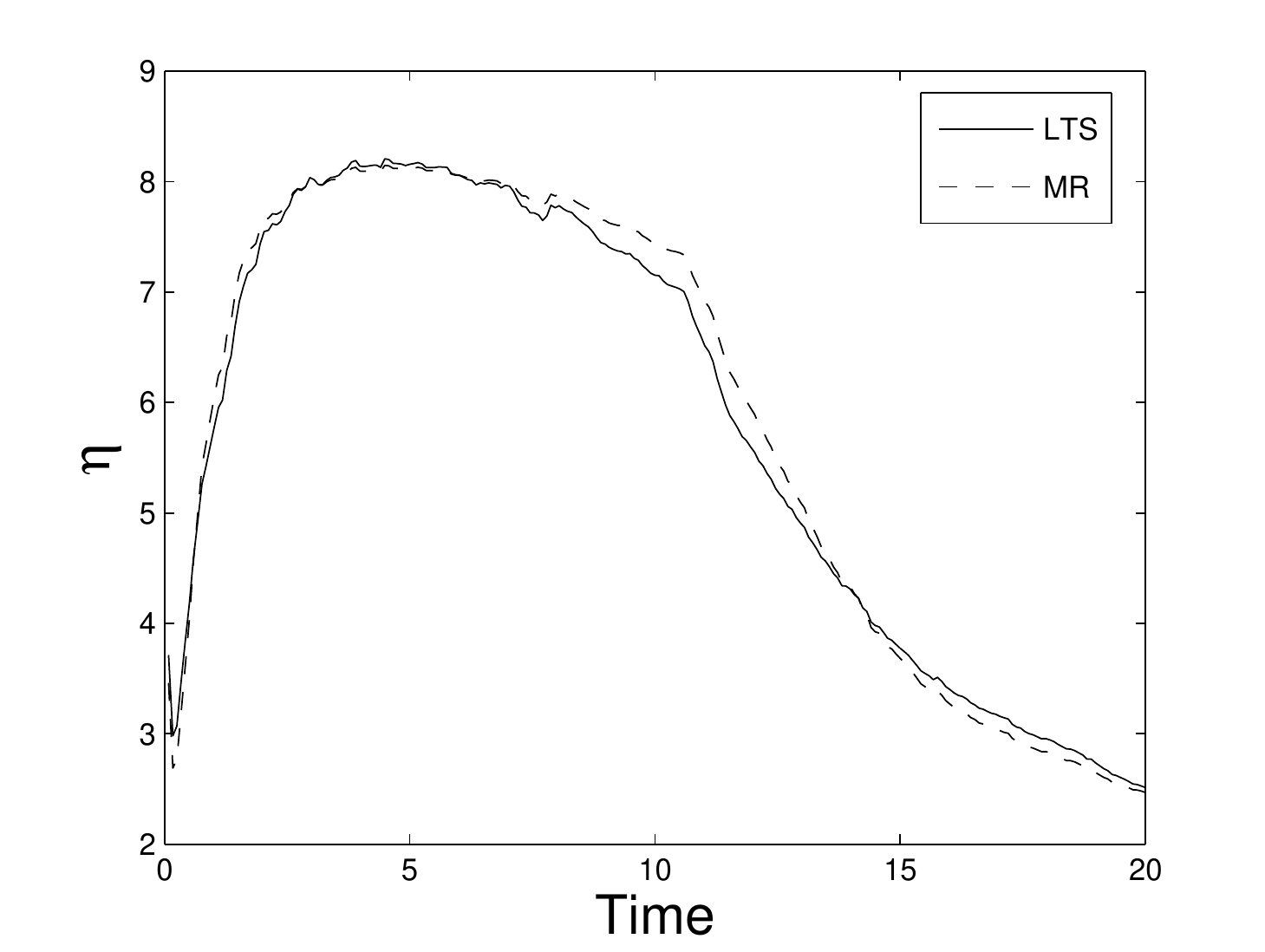}\\
(c) & (d) \\
\includegraphics[width=0.4\textwidth,height=0.35\textwidth]{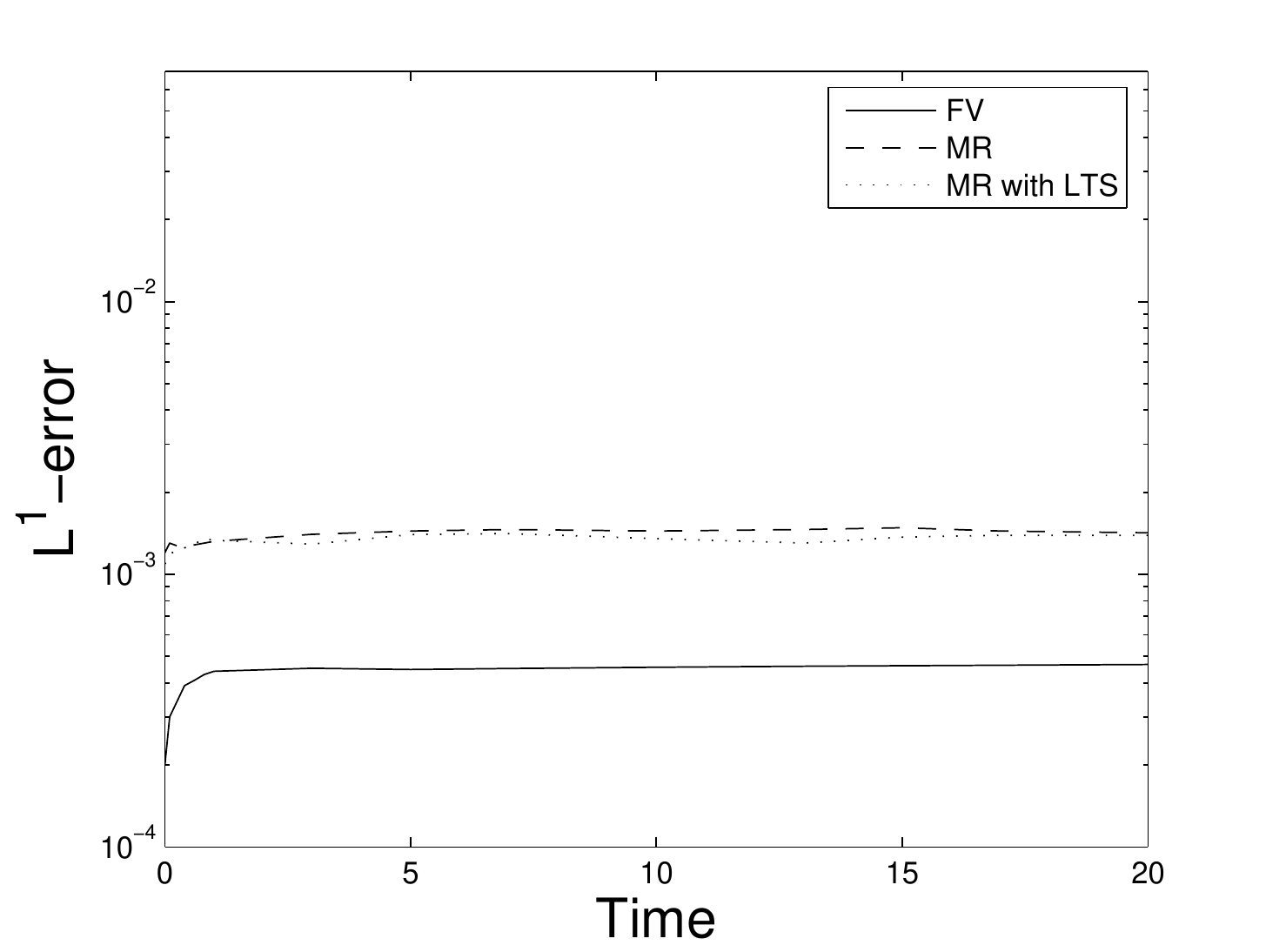}&
\includegraphics[width=0.4\textwidth,height=0.35\textwidth]{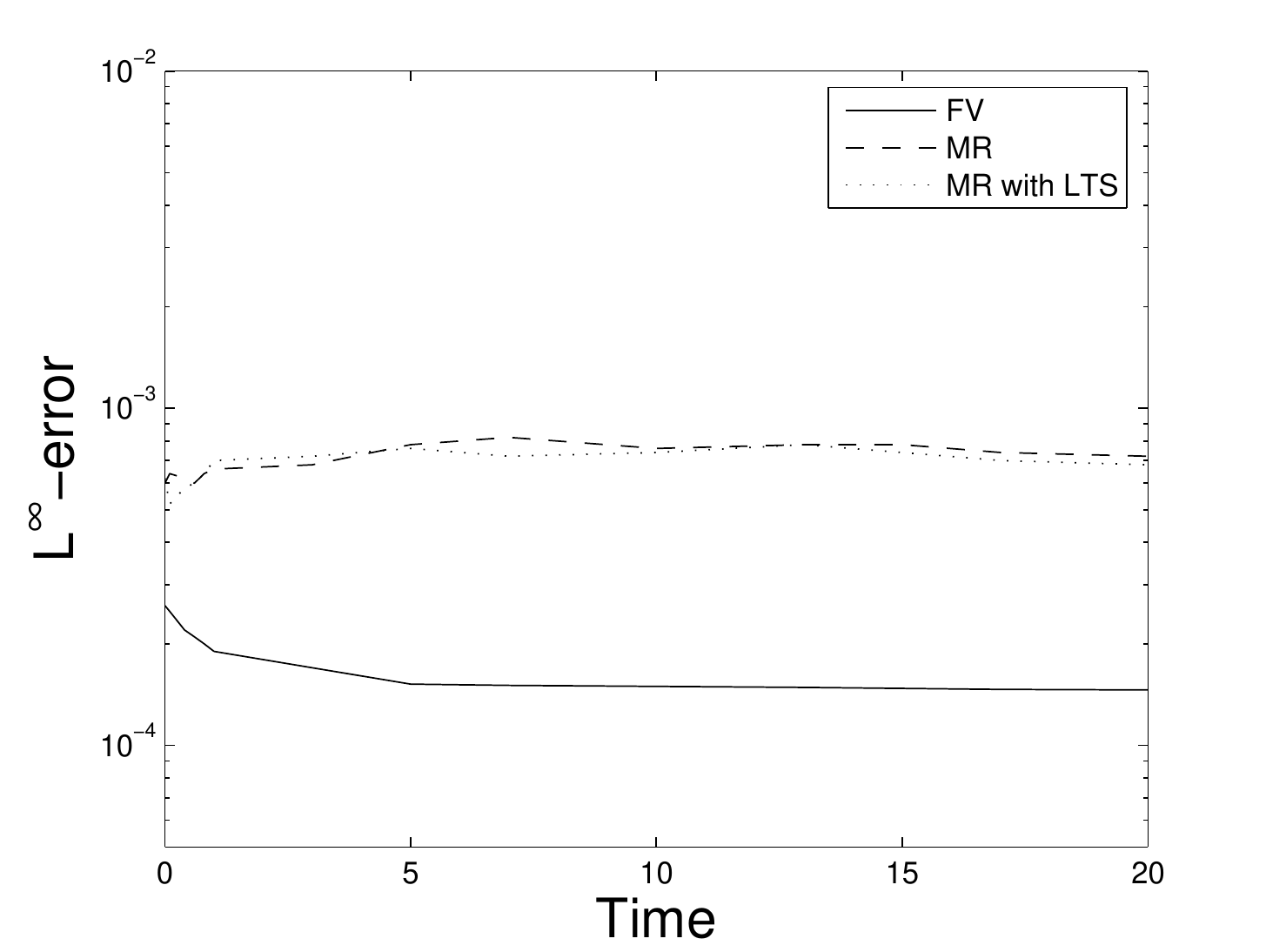}
\end{tabular}
\caption{\it Example~6 (chemotaxis with growth, $\nu=7$):  Time evolution for speed-up
rate~$V$, data compression rate~$\eta$  and  errors for different methods:
Multiresolution scheme with global time step, and multiresolution with
locally varying time step.} \label{bbrs_fig:ex4-varios}
\end{center}
\end{figure}

For Example~6, in \eqref{bbrs_eq:Kel-Seg} we consider a square domain $\Omega=[0,16]^2$
and fix the parameters $\sigma=0.0625$ and  $d=1$. The function $h(u,v)$ is given by
\eqref{bbrs_huvdef} with  $\alpha=1$ and $\beta=32$. The growth  function~$g(u)$
  for the species $u$ is given by \eqref{bbrs_gudef},  and the chemotactical sensitivity is
given by \eqref{bbrs_chinuv}. This configuration corresponds to the model of chemotaxis and
 growth presented in   \cite{mimura},  which is further analyzed
in   \cite{efendiev}. Similarly to \cite{efendiev}, the initial datum  is
$(u_0,v_0)=(1+\epsilon(\x),1/32)$,  where $\epsilon(\x)$ is a particular smooth
perturbation which goes to zero near $(8,8)$. We simulate the process until the
solution reaches inhomogeneous stationary states, and we present three cases
corresponding  to different values of~$\nu$, which is responsible
for the complexity of the spatial patterns. For example, for $\nu=7$
Figure~\ref{bbrs_fig:ex4-chemo} (middle) shows labyrinth-shaped patterns and for
$\nu=10$ (bottom), single filaments and spots. The corresponding adaptive
meshes were generated with $N_L=512^2=262144$ control volumes in the finest
grid, with $\varepsilon_{\mathrm{R}}=8.43\times10^{-4}$. For all these cases we
 implement  locally varying time stepping, so we will choose the maximum CFL
number allowed by \eqref{bbrs_cfl-chemo}, $\mathsf{CFL}_0=1$ for the coarsest level and
$\mathsf{CFL}_l=2^l\mathsf{CFL}_0$ for finer levels. From Figure~\ref{bbrs_fig:ex4-varios} we can observe
that if we incorporate the local time stepping strategy, a substantial gain
(a factor slightly lower than~2,  which is consistent with the results
 by Lamby, M\"{u}ller, and Stiriba
 \cite{LMS}) is obtained in speed-up rate when comparing with
a multiresolution calculation using global time stepping. The errors
are  computed  from a reference solution in a grid with $N_L=2048^2=4194304$
control volumes. We conclude that the errors are kept of the same slope
that the errors obtained with a global time step.

The compression rate
$\eta$ for both methods is lower than  in the previous examples, which
could be  explained by the  complexity and  density   of the spatial patterns in
this particular example.

\section{Conclusions} \label{bbrs_sec:conc}
This  paper describes an adaptive multiresolution scheme combined with a
locally varying time stepping used to approximate solutions of a
  class of two-dimensional
  reaction-diffusion systems in Cartesian geometry. Several numerical examples show
that the adaptive multiresolution mechanism with a suitable choice of the threshold
value represents a gain in CPU time while the errors are kept of the same order as the
reference finite volume method. In Examples~3 and~6,  we also see that the local time
stepping strategy is responsible for a gain in CPU time speed-up for a factor of
about 2. Also,  the errors between the solution using local time stepping and a reference
solution are of the same order that the solution obtained by the adaptive
multiresolution with global time stepping.

The motivation to employ explicit schemes only is the following.
 Even though implicit methods allow larger time steps, we need to
iteratively   solve a nonlinear system in each time step,
using e.g. Newton-Raphson method. The number of iterations is usually
controlled by  measuring the residual error, and cannot be controlled
a priori. Thus,  it appears  difficult to assess the true benefits of a
time-stepping strategy if the basic time discretization is an implicit one.

On the other hand, of course, for the Turing-type pattern formation problem,
it is conceded that patterns appear when one eigenvalue goes from
negative to positive. At steady state (when the pattern is visible) all
the eigenvalues again have negative real part. Thus to converge to steady
state once the domain of attraction of the pattern is reached, implicit
methods offer significant advantages since they can use larger and larger
time steps.

 We remark that for hyperbolic problems,
the incorporation of an implicit time discretization to the MR-LTS strategy
 can possibly form
a substantial improvement in the speed-up rate,  as presented in
\cite{MS}.

\section*{Acknowledgments}
MB acknowledges support by Fondecyt project 1070682, RB acknowledges
support by Fondecyt project 1050728 and Fondap in Applied Mathematics,
 project 15000001, RR acknowledges
support by Conicyt Fellowship and Mecesup projects UCO0406 and
UCO9907, and KS acknowledges support by the
 Agence Nationale de la Recherche, project M2TFP.
This work was partially done while RR visited
the Laboratoire de Mod\'{e}lisation et Simulation Num\'{e}rique en
 M\'{e}canique du CNRS and the Centre de Math\'{e}matiques et d'Informatique at the
 Universit\'{e} de Provence in Marseille, France.


\end{document}